\documentclass[12pt, openany]{memo-l}

\usepackage{amsmath, amsthm, amssymb, amscd, amsfonts}

\newtheorem{theoremmain}{Theorem}
\newtheorem{theorem}[equation]{Theorem}
\newtheorem{lemma}[equation]{Lemma}
\newtheorem{corollary}[equation]{Corollary}
\newtheorem{proposition}[equation]{Proposition}
\newtheorem{claim}[equation]{Claim}

\newtheorem{question}{Question}

\theoremstyle{definition}

 \newtheorem{example}[equation]{Example}

\theoremstyle{remark}
\newtheorem{remark}[equation]{Remark}

\numberwithin{section}{chapter}
\numberwithin{equation}{section}

\newcommand\Hom{\operatorname{Hom}}
\newcommand\id{\operatorname{id}}
\newcommand\ad{\operatorname{ad}}
\newcommand\Tr{\operatorname{Tr}}

\newcommand\co{\operatorname{co}}

\newcommand\cop{\operatorname{cop}}
\newcommand\End{\operatorname{End}}

\newcommand\Aut{\operatorname{Aut}}
\newcommand\op{\operatorname{op}}

\newcommand\Ind{\operatorname{Ind}}

\newcommand\Supp{\operatorname{Supp}}

\renewcommand{\subjclassname}{%
      \textup{2000} Mathematics Subject Classification}

\begin{document}
\frontmatter

\title[Semisimple Hopf Algebras of Low Dimension]{On the Semisolvability of Semisimple Hopf Algebras of
Low Dimension}
\author{Sonia Natale}
\address{Facultad de Matem\' atica, Astronom\' \i a y F\' \i sica,
Universidad Nacional de C\' ordoba, CIEM--CONICET, Ciudad Universitaria, (5000) C\' ordoba, Argentina}
\email{natale@mate.uncor.edu
\newline
\indent \emph{URL:}\/ http://www.mate.uncor.edu/natale}
\thanks{This work was partially supported by CONICET, Agencia C\'ordoba Ciencia,
ANPCyT, Fundaci\'on Antorchas and Secyt (UNC)} 
\date{September  2006}
\dedicatory{To the memory of my father}
\begin{abstract} We prove that every semisimple Hopf algebra of dimension less than
$60$ over an algebraically closed field $k$ of characteristic zero
is either  upper or lower semisolvable up to a cocycle twist.
\end{abstract}

\maketitle

\setcounter{page}{4}
\tableofcontents

\mainmatter

\chapter*{Introduction and Main Results}

 In recent papers several notions and results from the  theory of finite groups
have been generalized or adapted to the context of (semisimple)
Hopf algebras. Simultaneously, many results on the classification
of semisimple Hopf algebras have also  appeared. Some conjectural
analogies, such as Kaplansky's conjecture about the  dimensions of
the irreducible modules, still remain an open problem.

Let $H$ be a finite dimensional Hopf algebra over a field $k$.
 A Hopf subalgebra $A \subseteq H$ is called  \emph{normal} if $h_1 A \mathcal S(h_2) \subseteq A$, for all $h \in H$.  If $H$ does not contain proper normal Hopf subalgebras then it is called  \emph{simple}.

If $A \subseteq H$ is a normal Hopf subalgebra then the structure of $H$ can be reconstructed from $A$ and the quotient Hopf algebra $\overline H = H/ H A^+$; more precisely, it is known that in this case $H$ is isomorphic to a bicrossed product
$H \simeq A \#_{\rightharpoonup, \sigma}^{\rho, \tau} \overline H$, where  $(\rightharpoonup, \sigma, \rho, \tau)$ is a  \emph{compatible datum}; see for instance \cite{and-ext, Maext, ma-newdir}.
This fact implies that, when trying to classify Hopf algebras of a given finite dimension, it is an important problem to decide whether the Hopf algebra is simple or not.

\medbreak
\emph{We shall assume from now on that the field $k$ is algebraically closed of characteristic zero.}

\medbreak
We say that a finite dimensional Hopf algebra $H$ is \emph{trivial} if it is isomorphic to a group algebra or to a dual group algebra. Then, $H$ is trivial if and only if it is commutative or cocommutative.

The notions of upper and lower semisolvability for finite-dimensional Hopf algebras have been introduced in \cite{MoW}, as generalizations of the  notion of solvability for finite groups.
By definition, $H$ is called \emph{lower
semisolvable} if there exists a chain of Hopf subalgebras $H_{n+1}
= k \subseteq H_n \subseteq \dots \subseteq H_1 = H$ such that
$H_{i+1}$ is a normal Hopf subalgebra of $H_i$, for all $i$, and
all \emph{factors} ${\overline H}_i : = H_{i+1}/H_{i+1}H_i^+$ are
trivial. Dually, $H$ is called \emph{upper
semisolvable} if there exists a chain of quotient Hopf algebras
$H_{(0)} = H \to H_{(1)} \to \dots
\to H_{(n)}= k$ such that each of the maps $H_{(i-1)} \to
H_{(i)}$ is normal, and all
\emph{factors} $H_i : = H_{(i-1)}^{{\rm co } \pi_i}$ are trivial; here, $H_{(i-1)}^{{\rm co } \pi_i}$ is the space of coinvariants of the map $\pi_i$, see Section \ref{co-q}.

 We have that $H$ is upper semisolvable if and
only if $H^*$ is lower semisolvable \cite{MoW}. If this is the case, then
$H$ can be obtained from group algebras and their duals by means
of (a finite number of) extensions; in particular, $H$ is semisimple.

\medbreak The smallest non-solvable group is the simple
alternating group ${\mathbb A}_5$ of order 60. It is thus natural
to ask if an analogous statement is true for semisimple Hopf
algebras. A version of the following question was posed by S.
Montgomery in \cite[Question, pp. 269]{Mo1}.

\begin{question}\label{que}  Let $H$ be a semisimple
Hopf algebra of dimension less than $60$. Is $H$ necessarily  upper or  lower semisolvable? \end{question}

\medbreak
Let $H$ be a semisimple Hopf algebra over $k$. If $\dim H = p^n$, where $p$ is a prime number, then $H$ has a nontrivial central group-like element
\cite{masuoka-p^n}; inductively, one can see that $H$ is both upper and lower
semisolvable \cite{MoW}. Also, if $\dim
H = pq^2$, where $p \neq q$ are prime numbers, then it was shown
in \cite{pqq, pqq2, clspqq} that, under the assumption that $H$
and $H^*$ are both of Frobenius type, either $H$ or $H^*$ contains
a nontrivial central group-like element. This implies that these
Hopf algebras are also semisolvable, since semisimple Hopf algebras of dimension $p$, $pq$ and $q^2$ are trivial. In \cite{clspqq} we showed
that all semisimple Hopf algebras of dimension $pq^2 < 100$ are of Frobenius type (some instances of this fact, e.g., dimension 44, appeared in \cite{comp-k}); so that these are all semisolvable.

However, not every nontrivial semisimple Hopf algebra $H$ is
semisolvable. An example of a  simple nontrivial semisimple Hopf
algebra $H$ of dimension 60 and was constructed by D. Nikshych in
\cite{nik}: in this case $H$ is a cocycle twist of the group
algebra of the simple group ${\mathbb A}_5$. Moreover, it was
shown in \cite{nik} that if $G$ is a finite simple group and $\phi
\in kG \otimes kG$ is a nontrivial invertible pseudo $2$-cocycle,
then the twisted group algebra $(kG)_{\phi}$ is a nontrivial
semisimple Hopf algebra, which is simple as a Hopf algebra.

The smallest example of a semisimple Hopf algebra which is not
semisolvable is a cocycle twist of a group of order 36
\cite{twist-simple}. So the answer to Question \ref{que} is
negative, and it can only  be expected to be affirmative up to a
cocycle twist. The dimensions where the problem remains open are
24, 30, 36, 40, 42, 48, 54 and 56. We refer the reader to
\cite{andrusk, Mo1} for an account of previous results on the
problem of classification.

We also point out that, in the related context of Kac algebras, several  classification
results in low dimension were obtained by Izumi and Kosaki in
their work \cite{IK}; in that paper, the authors classify all Kac
algebras of dimensions 16, 24, $pq^2 < 60$ and $pqr < 60$.

 Our main result is the following theorem, giving an affirmative answer to Question \ref{que} up to cocycle twists.

\begin{theoremmain}\label{th-ss} Let $H$ be a semisimple Hopf algebra of dimension
$< 60$.   Then $H$ is either upper or lower semisolvable up to a
cocycle twist. \end{theoremmain}

We prove that a semisimple Hopf algebra of dimension  24, 30, 40,
42, 48, 54 or 56 is \emph{not simple}, and moreover, in dimension
36 the only simple example is a twisting of a finite group. This
is equivalent to the statement in the theorem, in view of previous
results. Indeed, if $H$ is a nontrivial semisimple Hopf algebra of
dimension $< 60$, then the following are equivalent (see
\cite{ama-sonia}):

\begin{itemize} \item $H$ is not simple;

\item  $H$ is either upper or lower semisolvable. \end{itemize}

Our approach to the problem is the  following: for each fixed
dimension, we first consider the possible algebra and coalgebra
decompositions (which turn out to be of Frobenius type). Then, for each possible type, we derive the existence of proper normal Hopf subalgebras.

In order to do this, we discuss some general result on semisimple Hopf algebras. Many of these results are new.
We discuss some properties of irreducible characters of low degree
which allow, in most cases, to prove the existence  of quotient
Hopf algebras or Hopf subalgebras, for each fixed algebra or
coalgebra structure, respectively. One of the main tools for this is the use of the Nichols-Richmond theorem on irreducible
characters of degree $2$ \cite{NR} and some of its consequences,
which we develope in Chapter \ref{dos}.

Using the character theory of $H$, we get a general result on Hopf subalgebras with index 3, which is often of use in low dimensions. See Theorem \ref{coinvariantes}.

\medbreak
In Chapter \ref{tres} we consider an inclusion $A \subseteq H$ of semisimple Hopf algebras. We show that if $C$ is a simple  subcoalgebra of $H$ such that $C
a \subseteq C$, for all $a \in A$, then the dual of the quotient coalgebra $C
/ CA^+$ and the crossed product $A_{\alpha}$,
where $\alpha : A \otimes A \to k$ is a certain $2$-cocycle, constitute a commuting pair in $C^*$.

This relates the corepresentation theory of $C / CA^+$ with the
representation theory of the Galois object $A_{\alpha}$ of $A$.
Applied in combination with   Masuoka's main result in
\cite{masuoka-cont} on deformations of cosemisimple Hopf algebras,
the result allows us to prove that some coalgebra decompositions
are impossible. The result is also useful to get some information
on the structure of the Hopf subalgebra $A$, especially when $A$
is a group algebra; see Section \ref{str-G}.

\medbreak
We also discuss braided Hopf algebras in relation with the Radford-Majid biproduct construction; see Chapter \ref{biprod}. It often happens, mainly because of the 'self-dual' nature of the assumption of simplicity, that if a given semisimple Hopf algebra $H$ is simple, then $H$ must have the structure of a biproduct $H = R \# A$, where $A$ is a semisimple Hopf subalgebra. We get several results on existence of proper (normal) Hopf subalgebras in biproducts; in particular, we give such a result in Corollary \ref{bip-8} for the case where $A$ is not cocommutative of dimension $p^3$, $p$ a prime number.

\medbreak
Let $p \neq q$ be prime numbers. We describe the known families in dimension $pq^2$ as cocycle twists of group algebras; we also prove that other families cannot be obtained in this fashion. See Chapter \ref{twist}.

 The classification of semisimple Hopf algebras of dimension $pqr$,
where $p$, $q$ and $r$ are distinct prime numbers, was given in
\cite{pqq} under the  assumption that $H$ admits an extension with
commutative 'kernel' and cocommutative 'cokernel' (a so called
\emph{abelian} extension).

In this paper we also prove that
semisimple Hopf algebras of dimension $30$ and $42$ admit abelian
extensions. This allows us to give the complete classification of
semisimple Hopf algebras of these dimensions. We obtain the
following theorems. See Chapters \ref{30}, \ref{42}.

\begin{theoremmain}\label{cls30} Let $H$ be a semisimple Hopf algebra of dimension
$30$ over $k$. Then $H$ is isomorphic to a group algebra $kG$ or to a dual
group algebra $k^G$, where $G$ is a group of order $30$. \end{theoremmain}

\medbreak
The known nontrivial examples in dimension $42$ were constructed in
\cite{examples}: these are denoted $\mathcal A_7(2, 3)$ and
$\mathcal A_7(3, 2) \simeq \mathcal A_7(2, 3)^*$.  We prove in this paper the following result.

\begin{theoremmain}\label{cls42} Let $H$ be a nontrivial semisimple Hopf algebra of dimension
$42$ over $k$. Then $H$ is isomorphic to one of the Hopf algebras
$\mathcal A_7(2, 3)$ or $\mathcal A_7(3, 2)$. \end{theoremmain}

\label{tabla} In Table \ref{tabl} we resume some known
facts about the classification of semisimple Hopf algebras of
dimension less than $60$.  In the first column, $p$, $q$ and $r$
are distinct prime numbers. The table is organized as follows: the
first column indicates the factorizations of the dimensions $\leq 60$ in terms of prime numbers; the second
and third column contain the references where nontrivial
examples were constructed and where the classification was given,
respectively. The fourth column gives additional information for
each specific type.

\begin{table}[t]
\begin{center}
\tiny{\begin{tabular}{|p{1,5cm}|p{3cm}|p{2,8cm}|p{3,8cm}|}
\hline  {\bf $\dim H$} & {\bf Nontrivial
\newline examples} & {\bf Classification} & {\bf
Remarks} \\ \hline $p$ & -- &   \cite{kac-p,
Z}. See also \cite{ZS}. & $H \simeq k\mathbb Z_p$. \\
\hline  $p^2$ & -- &  \cite{ma-pp}. & $H \simeq
kG$, $G$ abelian of order $p^2$ \\ \hline  $ pq$ & --
&  \cite{ma-6-8, EG, GW}. &
$H \simeq kG$, or $k^G$. Other proofs in \cite{So,
pqq}. \\ \hline
$p^3$ & \cite{k-p} for dim $8$.
\newline \cite{ma-pp} for dim $27$. & \cite{ma-6-8,
ma-pp}. & All
semisolvable by
\newline \cite{masuoka-p^n, MoW}. \\ \hline $pq^2$ &
\cite{fukuda, G, masuoka-further, pqq}. & \cite{fukuda}
for dim
$12$, \newline  \cite{masuoka-further} for dim $18$,
\newline
\cite{pqq, pqq2, clspqq}. & All semisolvable. \newline
\cite[Chapter X]{IK} for Kac algebras.\\ \hline $ p^4$ ($16$) &
\cite{kashina}. & \cite{kashina}. & All semisolvable
by
\newline \cite{masuoka-p^n, MoW}. \\ \hline
$ p^3q$ ($24$, $40$, $54$, $56$) & \cite{IK} for dim
$24$.
& \cite[Chapter XIV]{IK} for Kac algebras of dim
$24$. &
Semisolvable; \newline Chapters \ref{24}, \ref{40}, \ref{54}, \ref{56}.
\\ \hline $ pqr$ ($30$, $42$) & \cite{examples} for
dim
$42$. No nontrivial example in dim $30$. &  Chapters
\ref{30},
\ref{42}. & Abelian extensions in \newline dim $pqr$
classified in
\cite{pqq}. \cite[Chapter X]{IK} for Kac algebras. \\
\hline
$ p^5$ ($32$) & Yes. &  -- & All semisolvable by
\newline \cite{masuoka-p^n,
MoW}. \\ \hline $ p^2q^2$ ($36$) & \cite{eg-triangular};
$D(\mathbb S_3)$; \newline $(kD_3 \times D_3)_{\phi}$ simple
\cite{twist-simple}; more. & -- & First non-semisolvable
\newline example. Chapter \ref{36}.
\\ \hline $ p^4q$ ($48$) & Yes. & -- & Semisolvable; Chapter \ref{48} \\ \hline
$ p^2qr$ ($60$) & $(k\mathbb A_5)_{\phi}$ (simple) \cite{nik};
$(kD_3 \times D_5)_{\phi}$ (simple) \cite{twist-simple}. & -- &
--
\\ \hline
\end{tabular}}
\end{center}

\

\caption{Semisimple Hopf algebras of dimension $\leq
60$}\label{tabl}
\end{table}

 We include an appendix where we describe the structure of the Drinfeld doubles of the three non-commutative semisimple Hopf algebras of dimension 8. Tambara and Yamagami show in \cite{ty} that the categories of representations of these three Hopf algebras are not equivalent as monoidal categories.
The results in this appendix have been motivated by the paper \cite{Mo-8}, where the Schur indicators for the three Hopf algebras are compared: this gives  evidence that the representation theory of $H_8$ is in some sense closer to that of $kD_4$ than to that of $kQ$.  Here, $H_8$ denotes the unique nontrivial $8$-dimensional semisimple Hopf algebra over $k$ \cite{k-p, ma-6-8}, while $D_4$ and $Q$ denote, respectively, the dihedral and quaternionic groups of order $8$.

On the other hand, note that $H_8$ is an extension of $k^{\mathbb
Z_2 \times \mathbb Z_2}$ by $k\mathbb Z_2$, such that the
bicrossproduct group $(\mathbb Z_2 \times \mathbb Z_2) \bowtie
\mathbb Z_2$ associated to the matched pair of the extension is
$D_4$. As a consequence of \cite[Theorem 1.3]{gp-ttic}, we know
that $D(H_8)$ is a cocycle twist of the Dijkgraaf-Pasquier-Roche
quasi-Hopf algebra $D^{\omega}(D_4)$,  where $\omega \in H^3(D_4,
k^{\times})$ is the 3-cocycle associated to the extension
corresponding to $H_8$ via the Kac exact sequence \cite{kac}. We
present more evidence of this facts involving the Drinfeld
doubles. More precisely, we prove on the one hand that  $D(H_8)$
has no quotient Hopf algebra isomorphic to $kQ$. We also show that
$D(H_8)$ is a (central) extension of $k^G$ by $kG$, where $G =
G(D(H_8)^*) \simeq \mathbb Z_2 \times \mathbb Z_2 \times \mathbb
Z_2$. See Theorems \ref{th1}, \ref{th2}.

\subsection*{Acknowledgements.} Most of the results in this paper have been announced in \cite{ama-sonia}. They were also communicated in the conferences \emph{"Hopf algebras in Noncommutative Geometry and Physics"} (Brussels, May 2002) and \emph{"Colloquium on Homology Theories, Representations and Hopf Algebras"} (Luminy, June 2002).

The author is grateful to N. Andruskiewitsch, S. Montgomery, L.
Vainerman and  H.-J. Schneider  for interesting
discussions, comments and references. She also thanks  Y. Kashina and Y. Sommerh\" auser for helpful remarks on a previous version of this paper, and  the referee for many valuable comments.

This research has been done during a postdoctoral stay at the Department of Mathematics of the \' Ecole Normale Sup\' erieure, Paris. The author is grateful to Marc Rosso for his kind hospitality.

\bigbreak
\chapter*{Conventions and Notation.} \label{conv}

Throughout,  $k$ will denote an
algebraically closed field of characteristic zero.  The symbols
$\Hom$, $\otimes$, etc., will mean $\Hom_k$, $\otimes_k$, etc. Our
references for the theory of Hopf algebras are \cite{Mo, Sch}. The notation for Hopf algebras is standard; for
instance, the group of group-like elements in $H$ is denoted by
$G(H)$.  For an algebra $A$ (respectively, for a coalgebra $C$)
the notation $\Hom_A$ (resp. $\Hom^C$)  is used to indicate the
$\Hom$ bifunctor in the category of (left) $A$-modules (resp.
$C$-comodules).

A Hopf  algebra $H$ is called  \emph{semisimple} (respectively,
\emph{cosemisimple}) if it is semisimple as an algebra
(respectively, if it is cosemisimple as a coalgebra). Let $H$ be a
finite-dimensional Hopf algebra over $k$. By a result of Larson
and Radford, it is known that $H$ is semisimple if and only if $H$
is cosemisimple, if and only if ${\mathcal S}^2 = \id$. See
\cite{LR, LR2}. The \emph{character algebra} of $H$, denoted $R(H)$, is the subalgebra of $H^*$ spanned by the irreducible characters of $H$; if $H$ is semisimple, $R(H)$ coincides with the subalgebra of cocommutative elements in $H^*$.

Suppose $H$ is finite dimensional. For a Hopf subalgebra $A
\subseteq H$, the \emph{index} of $A$ in $H$ is defined by $[H: A] :
= \dfrac{\dim H}{\dim A}$; it is an integer by \cite{NZ}. Suppose  $q: H \to B$ is a surjective Hopf algebra map, and identify $B^*$ with a Hopf subalgebra of $H^*$ via $q^*$; by abuse of terminology, the index $[H^*: B^*]$ will be also called the index of $B$ in $H$.

\chapter{Semisimple Hopf Algebras}\label{uno}
Along this chapter, $H$ will be a semisimple (thus
finite-dimensional) Hopf algebra over $k$.

\section{Algebra structure}\label{alg-struct} As an algebra, $H$ is isomorphic to a
direct product of full matrix algebras
\begin{equation}\label{estructura} H \simeq k^{(n)} \times \prod_{d_i > 1}
M_{d_i}(k)^{(n_i)},\end{equation} where $n = |G(H^*)|$.   It
follows from the Nichols-Zoeller Theorem \cite{NZ}, that $n$
divides both $\dim H$ and $n_i d_i^2$, for all $i$. Moreover, by
\cite{NR} if $d_i = 2$ for some $i$, then the dimension of $H$ is
even.

By \cite{ZS}, if $n = 1$, then $\{ d_i : d_i > 1 \}$ has at least three elements.

Suppose that $A \subseteq H$ is a Hopf subalgebra. Then $A$ is also semisimple. Assume that $A$ is commutative. Then it follows from the Frobenius Reciprocity that  $d_i \leq [H:A]$, for all $i$. See \cite[Corollary 3.9]{harmonic}.

If $H$ is as in \eqref{estructura} as an algebra, we shall say
that $H$ is \emph{of type} $(1, n; d_1, n_1; \dots; d_r, n_r)$
\emph{as an algebra}. In this case, the dimension of the character
algebra of $H$ is $n + n_1 + \dots + n_r$.

If $H^*$ is of type $(1, n; d_1, n_1; \dots )$ as an algebra, we
shall say that $H$ is  \emph{of type} $(1, n; d_1, n_1; \dots )$
\emph{as a coalgebra}.

So that $H$ is of type $(1, n; d_1, n_1; \dots; d_r, n_r)$ as a
(co-)algebra if and only if $H$ has $n$ non-isomorphic
one-dimensional (co-)representations,  $n_1$ non-isomorphic
irreducible (co-)representations of degree $d_1$, etc. Sometimes,
we shall use the notation $X_d$ or $X_d(H)$ to indicate the set of
irreducible characters of $H$ of degree $d$.

\begin{example} The above arguments can be used to get the possible algebra structures
for a given finite dimension. For instance, suppose that $B$ is a
semisimple Hopf algebra of dimension $60$ such that $G(B^*) = 1$.
Then, as an algebra, $B$ is of one of the following types:
$$(1, 1; 3, 2; 4, 1; 5, 1), \quad  (1, 1; 2, 4;  3, 2; 5, 1) \quad
\text{or} \quad (1, 1; 2, 4;  3, 3; 4, 1).$$  \end{example}

\section{Irreducible characters}\label{irr-char}
Let $V$ be an $H$-module. The \emph{character} of $V$ is the
element $\chi = \chi_V \in H^*$ defined by $\langle\chi, h \rangle
= \Tr_V(h)$, $h \in H$. The \emph{degree} of $\chi$ is the integer
$\deg \chi = \chi (1) = \dim V$. If $U$ is another $H$-module, we
have $$\chi_{V \otimes U} = \chi_V \chi_U, \qquad \chi_{V^*} =
\mathcal S(\chi_V).$$ Thus the irreducible characters, {\it i.e.},
the characters of the irreducible $H$-modules, span a subalgebra
$R(H)$ of $H^*$, called the \emph{character algebra} of $H$. The
antipode induces an anti-algebra involution $*: R(H) \to R(H)$,
$\chi \mapsto \chi^*: = \mathcal S(\chi)$. The degree defines an
augmentation $R(H)  \to k$.

We first resume some of the basic properties of $R(H)$ that will
be often used in the rest of this paper.  Proofs of these facts
can be found in \cite{NR}.

Let $\chi_V, \chi_W \in R(H)$ be the characters of the $H$-modules $V$ and $W$, respectively. The integer $m(\chi_V, \chi_W) = \dim \Hom_H(V, W)$ will be called \emph{multiplicity} of $V$ in $W$. This extends to a bilinear form $m: R(H) \times R(H) \to k$.

Let $\widehat H$ denote the
set of irreducible characters of $H$. If $\chi \in R(H)$, then we may write $\chi = \sum_{\mu \in \widehat H} m(\mu, \chi) \mu$. Let $\chi$, $\psi$ and $\lambda$ be characters of $H$-modules; we have
\begin{equation}m(\chi, \psi \lambda) = m(\psi^*, \lambda \chi^*) = m(\psi, \chi \lambda^*). \end{equation}

Let $\chi$  be an irreducible character of $H$. Denote by
$G[\chi]$ the subgroup of $G(H^*)$ consisting of all those
elements $g$ such that $g \chi = \chi$. We have
\begin{equation}\label{mult-g}g \in G[\chi] \,
\Longleftrightarrow \, m(g, \chi \chi^*) > 0 \,
\Longleftrightarrow \, m(g, \chi \chi^*) = 1.\end{equation} See \cite[Theorem
10]{NR}. In particular,
\begin{equation}\label{des-ss}\chi  \chi^* = \sum_{g \in G[\chi]} g +
\sum_{\mu \in \widehat H,  \deg \mu > 1} m(\mu, \chi \chi^*)
\mu.\end{equation} Note that $G[\chi^*] = \{ g \in G(H^*): \chi g
= \chi \}$. Also, if $a \in G(H^*)$, we have  $$G[\chi a] =
G[\chi], \qquad  G[a\chi] = aG[\chi]a^{-1}.$$ As a consequence of
\cite{NZ}, we have that $|G[\chi]| / (\deg \chi)^2$. Moreover,  it
follows from the results in \cite[Section 2]{masuoka-cont}, that
the order of any element $g \in G[\chi]$ (hence also the exponent
of $G[\chi]$) divides $\deg \chi$.

The following lemma will be applied later; see Lemma
\ref{conteo54}.  It serves here to illustrate the (mostly
well-known) fact that not any algebra type can arise as the
structure of a semisimple Hopf algebra.

\begin{lemma}\label{alg-type} There is no semisimple Hopf algebra with algebra type $(1, 2; 2, 1; 4, m)$, $m \geq 1$. \end{lemma}

\begin{proof} Suppose on the contrary  that  $H$ is a semisimple Hopf algebra with this
algebra type. Let $G(H^*) = \{ \epsilon, g \}$  and let $\chi$ be
the unique irreducible character of degree $2$ of $H$; so that we
have $g \chi = \chi = \chi g$, $\chi^* = \chi$ and $\chi^2 =
\epsilon + g + \chi$.  In particular, $m(\chi, \zeta \chi) =
m(\zeta, \chi^2) = 0$, for all irreducible character $\zeta$ of
degree $4$. Then also $m(\chi, \chi \zeta) = m(\chi, \zeta^* \chi)
= 0$, for all such $\zeta$.

Let $\deg \zeta = 4$  and write $\zeta \zeta^* = \epsilon + g + n
\chi + \lambda$, where $\lambda$ is a character of $H$ such that $m (\chi, \lambda) = 0$. By taking degrees,
we get that $n > 0$ and $n$ is odd. Since $n = m(\zeta, \chi
\zeta)$, it follows that $n = 1$. Therefore, $\chi \zeta  = \zeta
+ \psi$, where $\deg \psi = 4$, $\psi \neq \zeta$. Since $m(\chi, \chi \zeta) = 0$, then $\psi$ is irreducible.

Then $m(\epsilon, \chi \zeta \psi^*)  = m (\epsilon, \psi \psi^*)
= 1$, and we have $m(\chi, \zeta \psi^*) = 1$. Thus $\zeta \psi^*
= \chi + \sum_l \psi_l$, where $\deg \psi_l = 4$. Taking degrees
we get a contradiction. This shows that this type is not possible.
\end{proof}

Let $n \geq 1$. The group $G(H^*) \times G(H^*)$ acts on the set
$X_n$ via \begin{equation}\label{accion}(g, h) . \chi = g \chi
h^{-1}, \qquad g, h \in G(H^*), \quad \chi \in X_p.\end{equation}
We have $G[(g, h).\chi] = gG[\chi]g^{-1}$, for all $g, h \in
G(H^*)$, and for all $\chi$.

Using this action we can get some information on the structure of the group $G(H^*)$. The following proposition gives an example of this fact; see also the proof of Lemma \ref{gl-z2}.

\begin{proposition}\label{nab-pq} Let $p < q$ be prime numbers. Suppose that $G(H^*)$ is nonabelian of order $pq$. Assume in addition that $G[\chi] \neq 1$, for all $\chi \in X_p$. Then $q^2$ divides $|X_p|$. \end{proposition}

\begin{proof} We may assume that $X_p \neq \emptyset$. Let $T \subseteq G(H^*)$ be a subgroup of order $q$.
Consider the action $(T \times T) \times X_p \to X_p$ obtained by restriction of the action \eqref{accion}. We shall show that the stabilizer $(T \times T)_{\chi}$ is trivial, for all $\chi \in X_p$, which will imply the proposition.

Let $g, h \in T$, $\chi \in X_p$, and suppose that $(g, h).\chi = g \chi h^{-1} = \chi$. Then $G[\chi] = G[(g, h).\chi] = gG[\chi]g^{-1}$. By assumption, $G[\chi] \neq 1$, hence $|G[\chi]| = p$, because of the assumption $|G(H^*)| = pq$. Then $g = 1$; otherwise,  $G[\chi]$ would be a normal subgroup of order $p$ in $G(H^*)$, implying that $G(H^*)$ is abelian against the assumptions.

Hence we have $(g, h).\chi = \chi h^{-1} = \chi$, and $h^{-1} \in G[\chi^*]$. Thus also $h = 1$. This finishes the proof of the proposition. \end{proof}

\section{Coinvariants of Hopf algebra maps}\label{co-q} Let $q: H \to B$ be a Hopf algebra
map and consider the subspaces of coinvariants
\begin{align*}H^{\co q} & = \{ h
\in H: (\id \otimes q) \Delta (h) = h \otimes 1 \}, \quad \text{and}\\
{}^{\co q}H & = \{ h \in H: (q \otimes \id) \Delta (h) = 1 \otimes h \}. \end{align*}

Then $H^{\co q}$ (respectively, ${}^{\co q}H$) is a left
(respectively, right) coideal subalgebra of $H$. We shall also use
the notation $H^{\co B} : = H^{\co q}$.
By \cite{sch-transitivity},
\begin{equation}\label{nl-bs}\dim H = \dim H^{\co q} \dim q(H) = \dim {}^{\co q}H \dim q(H).\end{equation}

The left coideal subalgebra $H^{\co q}$ is stable under the left
adjoint action of $H$. Moreover $H^{\co q} = {}^{\co q}H$ if and
only if $H^{\co q}$ is a (normal) Hopf subalgebra of $H$. If this is
the case, we shall say that the map $q: H \to B$ is normal.

\begin{remark}\label{inv-reg} The Hopf algebra $B^*$ acts on $H$ on the
left and on the right by $f \rightharpoonup h = \langle f, q(h_2)
\rangle h_1$, and $h \leftharpoonup f = \langle f, q(h_1) \rangle
h_2$, respectively.  Suppose that $q$ is surjective. Then we have
\begin{align*}H^{\co q} & = \{ h \in H: f \rightharpoonup h =
\epsilon (f) h, \, \forall f \in B^* \}, \quad \text{and} \\
{}^{\co q}H & = \{ h \in H: h \leftharpoonup f = \epsilon (f) h,
\forall f \in B^* \}. \end{align*}

In particular, let $\eta \in G(H^*)$ and consider the Hopf algebra
map $q: H \to k^{\langle \eta \rangle}$ obtained by transposing the inclusion
$k\langle \eta \rangle \subseteq H^*$. Then
\begin{equation*}H^{\co q}  = \{ h \in H: \eta \rightharpoonup h =
h\}, \qquad {}^{\co q}H  = \{ h \in H: h \leftharpoonup \eta = h\},
\end{equation*} where $\leftharpoonup$ and $\rightharpoonup$ are the regular
actions of $H^*$ on $H$.
\end{remark}

By \cite{NZ}, if $A$ is a Hopf  subalgebra of $H$ such that $A
\subseteq H^{\co q}$, then $H^{\co q}$ is free as an $A$-module,
with respect to the action by left multiplication of $A$. In
particular, $\dim A$ divides $\dim H^{\co q}$. The same holds true with  ${}^{\co q}H$ instead of $H^{\co q}$. Indeed, with this $A$-module structure and the coaction given by the comultiplication of $H$, $H^{\co q}$ is a left $(A, H)$ Hopf module.
We note the following consequence of this fact:

\begin{lemma} Suppose that $|G(H)|$ and $[H^* : G(H^*)]$ are relatively prime.
Then the group $G(H)$ is abelian and isomorphic to a subgroup of
$\widehat{G(H^*)}$.
\end{lemma}

\begin{proof} There is a surjective Hopf algebra map
$\pi: H \to k^{G(H^*)}$, and we have $\dim H^{\co \pi} = [H^* :
G(H^*)]$. If $1 \neq g \in G(H)$, then $\pi (g) \neq 1$,  since
otherwise $g$ would belong to $H^{\co \pi}$ implying that the
order of $g$ divides $\dim H^{\co \pi}$, contradicting the
assumption. Therefore the restriction of $\pi$ to $G(H)$ is
injective, and $G(H)$ is thus isomorphic to a subgroup of
$\widehat{G(H^*)} = G(k^{G(H^*)})$. This implies that $G(H)$ is abelian.
\end{proof}

\begin{lemma}\label{restriction} Suppose that $A \subseteq H$ is a Hopf subalgebra.
Then $A^{\co q\vert_A} = A \cap H^{\co q}$. In particular,  $\dim A = \dim
(A \cap H^{\co q}) \ \dim q(A)$.  \end{lemma}

\begin{proof} The first claim is evident. The second follows from \eqref{nl-bs}. \end{proof}

\section{Yetter-Drinfeld modules}\label{ss-yd} Let ${}_H^H\mathcal{YD}$ denote the  category of
(left-left) Yetter-Drinfeld modules over $H$. Objects of this category
are vector spaces $V$ endowed with an $H$-coaction $\rho: V \to H
\otimes V$ and an $H$-action $. : H \otimes V \to V$, subject to
the compatibility condition $\rho (h . v) = h_1 v_{-1} \mathcal
S(h_3) \otimes h_2 . v_0$, $v \in V$, $h \in H$; morphisms are
$H$-linear and colinear maps.

The category ${}_H^H\mathcal{YD}$ is a modular category,  which
coincides as such with the category of modules over the Drinfeld
double of $H$; see \cite{EG}.
It is shown in \cite{EG} that if $V$ is a simple  Yetter-Drinfeld
module over $H$, then $\dim V$ divides $\dim H$.

With respect to the left adjoint action $\ad: H \otimes H \to
H$,  $(\ad h) (a) = h_1 a \mathcal S(h_2)$ and the left regular coaction
$\Delta: H \to H \otimes H$, $H$ becomes an object of
${}_H^H\mathcal{YD}$.

The Yetter-Drinfeld submodules $V \subseteq H$ are exactly the
left coideals $V$ of $H$ such that $h_1 V \mathcal S(h_2)
\subseteq V$, for all $h \in H$. Thus, a one-dimensional
Yetter-Drinfeld submodule of $H$ is exactly the span of a central
group-like element of $H$.

It is well-kown that the space of (right) coinvariants of a Hopf
algebra map is a left  coideal stable under the left adjoint
action. We thus obtain the following lemma:

\begin{lemma}\label{YD-coinv} Let $H \to B$ be a Hopf algebra map.
Then $H^{\co B}$ is a Yetter-Drinfeld submodule of $H$. \qed \end{lemma}

\begin{remark} Suppose that $H \to B \to B'$ is a sequence of Hopf algebra maps.
Then $H^{\co B}$ is a Yetter-Drinfeld
submodule of $H^{\co B'}$.  In particular, since ${}_H^H\mathcal{YD}$ is a
semisimple category, there exists a Yetter-Drinfeld submodule $W
\subseteq H^{\co B'}$ such that $H^{\co B'} = H^{\co B} \oplus W$.
\end{remark}

We recover the following result, due to Kobayashi and Masuoka \cite{kob-mas}; see also \cite[Theorem 2.1.1]{pqq}. Our alternative proof is based on \cite{EG}.

\begin{corollary}\label{kob-mas} Suppose that $B \subseteq H$ is a
Hopf subalgebra such that $[H : B] = p$ is the smallest prime
number dividing $\dim H$.  Then $B$ is a normal Hopf subalgebra
and $H$ fits into a (co-)central extension $1 \to B \to H \to
k\mathbb Z_p \to 0$. \qed \end{corollary}

Our argument proves indeed the following more precise statement:  if $A \subseteq H$ is a \emph{normal}  Hopf subalgebra such that
$\dim A$ is the smallest prime number dividing $\dim H$, then $A$
is \emph{central} in $H$.

\begin{proof} Consider  the dual projection $H^* \to B^*$.
So that we have $\dim (H^*)^{\co B^*} = [H : B] = p$.   Let $V$ be an
irreducible Yetter-Drinfeld submodule of $(H^*)^{\co B^*}$. Since
the dimension of $V$ divides $\dim H$ and is less than $p$, we
find that $\dim V  = 1$; therefore $V = kg$, for some $g \in
Z(H^*) \cap G(H^*)$.

Decomposing $(H^*)^{\co B^*}$ into irreducible  Yetter-Drinfeld
modules, we see that $(H^*)^{\co B^*}$ is a central group-like
Hopf subalgebra of $H^*$ of dimension $p$. This implies that $H^*$
fits into a central extension $0 \to k\mathbb Z_p \to H^* \to B^*
\to 1$. The lemma follows after dualizing this extension.
\end{proof}

\section{Yetter-Drinfeld modules and the character algebra}In the paper \cite{zhu2} Y. Zhu establishes a bijective correspondence between primitive idempotents in the character algebra $R(H) \subseteq H^*$ and irreducible Yetter-Drinfeld submodules of $H$.
Indeed, it is shown in \cite{zhu2} that $D(H)$ and $R(H)$ form a commuting pair in $\End H$, with respect to the $D(H)$-action corresponding to the Yetter-Drinfeld module structure in $H$ considered in Section \ref{ss-yd} (or a version of it thereof) and the $R(H)$-action $\rightharpoonup: R(H) \otimes H \to H$, $f \rightharpoonup h : = \langle f, h_2\rangle h_1$.

Let $q: H \to B$ be a Hopf algebra projection. Then $H^{\co B}$ is a Yetter-Drinfeld submodule of $H$.
Consider the dual inclusion of Hopf algebras $B^* \to H^*$ and let $e_0 \in B^*$ be the normalized integral; $e_0$ is the primitive idempotent in $B^*$ corresponding to the trivial representation. Since $e_0$ is a cocommutative element, we have $e_0 \in R(H)$. Hence we may write
\begin{equation}\label{idpts}e_0 = \Lambda + E_1 + \dots + E_n,
\end{equation}
where $\Lambda, E_1, \dots, E_n$ are orthogonal primitive idempotents in $R(H)$, such that $\Lambda$ is the normalized integral in $H^*$.

The following proposition gives a refinement of the result in \cite{zhu2}.

\begin{proposition}\label{ref-corr-z2} The idempotents $\Lambda, E_1, \dots, E_n$ in \eqref{idpts} correspond bijectively with the irreducible $H$-Yetter-Drinfeld submodules of $H^{\co B}$.
\end{proposition}

\begin{proof}  We saw in Remark \ref{inv-reg} that $H^{\co B}$ coincides with the subalgebra of $B^*$-invariant elements of $H$ under the left regular action $\rightharpoonup: B^* \otimes H \to H$. Hence, since $e_0 \in B^*$ is the primitive idempotent corresponding to the trivial representation, we have $H^{\co B} = e_0 \rightharpoonup H$. This implies the proposition. \end{proof}

\section{One dimensional Yetter-Drinfeld modules}
Let   $g \in G(H)$, $\eta \in G(H^*)$. Let $V_{g, \eta}$ denote
the one dimensional vector space endowed with the action $h.1 = \eta (h)1$,
$h \in H$, and the coaction $1 \mapsto g \otimes 1$. The following
 is a restatement of a result due to Radford \cite[Proposition 10]{R} that describes the group-like elements in the dual of $D(H)$.

\begin{lemma}\label{yd-1}
The one-dimensional Yetter-Drinfeld modules of $H$ are exactly
of the form  $V_{g, \eta}$, where $g \in G(H)$ and $\eta \in
G(H^*)$ are such that $(\eta \rightharpoonup h) g = g (h
\leftharpoonup \eta)$, for all $h \in H$. \qed  \end{lemma}

\begin{remark}\label{yd-*} Let $g \in G(H)$, $\eta \in G(H^*)$.
Then $V_{g, \eta}$ is a Yetter-Drinfeld module of $H$ if and only if $V_{\eta, g}$ is a Yetter-Drinfeld module of $H^*$.

\begin{proof} We use Lemma \ref{yd-1}. We have that $V_{g, \eta}$ is a Yetter-Drinfeld module of $H$ if and only if $(\eta \rightharpoonup h) g = g (h
\leftharpoonup \eta)$, for all $h \in H$.
This equivalent to $(g \rightharpoonup f) \eta = \eta (f
\leftharpoonup g)$, for all $f \in H^*$.  Hence the claim follows.
\end{proof} \end{remark}

\begin{lemma}\label{inv-eta} Let  $g \in G(H)$ and $\eta \in G(H^*)$ such that $V_{g,
\eta}$ is a Yetter-Drinfeld module of $H$. Let also $q: H \to
k^{\langle \eta \rangle}$ be the Hopf algebra map obtained by transposing the
inclusion $k\langle \eta \rangle \subseteq H^*$.

Suppose that $V \subseteq H^{\co q}$ is a subspace of $H$ such
that $g^{-1} v g = v$, for all $v \in V$. Then $V \subseteq
{}^{\co q}H$.
\end{lemma}

\begin{proof} Using Lemma \ref{yd-1}, we have for all $v \in V$,
\begin{equation*} v \leftharpoonup \eta  = g^{-1} (\eta \rightharpoonup v)
g = g^{-1}vg = v,
\end{equation*}
the second equality because $v \in H^{\co q}$, see Remark
\ref{inv-reg}. This shows that $v \in {}^{\co q}H$ and finishes
the proof of the lemma. \end{proof}

The following result will be used in Chapters \ref{24} and
\ref{54}.

\begin{theorem}\label{A-comm-g} Let $g$, $\eta$ and $q$ be as in Lemma \ref{inv-eta}.
Let $A \subseteq H$ be a Hopf subalgebra such that $g^{-1} a g =
a$, for all $a \in A^{\co q}$. Then  the restriction $q\vert_A: A \to
k^{\langle \eta \rangle}$ is normal.
\end{theorem}

\begin{proof} We shall prove that $A^{\co q}$ is a Hopf
subalgebra of $A$. By Lemma \ref{inv-eta}, we have $A^{\co
q} \subseteq {}^{\co q}A$. Thus ${}^{\co q}A = A^{\co q}$, since
they have the same finite dimension. This implies the theorem.
\end{proof}

\section{$H^{\co B}$ as a left coideal of $H$}
Let $q: H \to B$ be a surjective Hopf algebra map. Identify $B^*$ with its image under the transpose map $q^*: B^* \to H^*$; so that $B^*$ is a Hopf subalgebra of $H^*$.

For each left $B^*$-module $W$, we may consider the \emph{induced} left $H^*$-module $V = \Ind_{B^*}^{H^*}W : = H^*\otimes_{B^*}W$. Most basic properties of the induction functor are discussed, for instance, in \cite{harmonic}. Proposition \ref{ind-norm} below establishes a relationship between the decomposition of the induced module $\Ind_{B^*}^{H^*}\epsilon$ and the normality of the map $q$.
Here, $\Ind_{B^*}^{H^*}\epsilon$ indicates the representation induced from the \emph{trivial} one-dimensional representation $W = k_{\epsilon}$ of $B^*$.
The key ingredient is the following identification.

Recall that $H^{\co B}$ is a left coideal of $H$, hence a right $H^*$-module. Thus $(H^{\co B})^*$ is naturally a left $H^*$-module through the action given by $\langle p . f, x \rangle : = \langle p, x_1 \rangle  \langle f, x_2 \rangle$,  for $p \in H^*$, $f \in (H^{\co B})^*$, $x \in H^{\co B}$.

\begin{lemma}\label{idf-coinv} $(H^{\co B})^* \simeq \Ind_{B^*}^{H^*}\epsilon$ as left $H^*$-modules. \end{lemma}

As a consequence of this lemma, we see that  for every irreducible left coideal $V$ of $H$, $V^*$ appears in $H^{\co B}$ with the same multiplicity as $V$ does.

\begin{proof} As left $H^*$-modules, $\Ind_{B^*}^{H^*}\epsilon = H^*\otimes_{B^*}k_{\epsilon} \simeq H^*/H^*(B^*)^+$. On the other hand, the evaluation map $\langle \, , \, \rangle: H^* \otimes H \to k$ induces a left $H^*$-linear isomorphism $H^*/H^*(B^*)^+  \simeq (H^{\co B})^*$.
 \end{proof}

\begin{proposition}\label{ind-norm} The map $q: H \to B$ is normal if and only if every irreducible $H^*$-module $V$ appears with multiplicity $\dim V$ or 0 in $\Ind_{B^*}^{H^*}\epsilon$. \end{proposition}

\begin{proof} We know that $q$ is normal if and only if $H^{\co B}$ is a subcoalgebra of $H$. In turn, the last holds if and only if for every irreducible $H$-coideal $V \subseteq H^{\co B}$, $H^{\co B}$ contains the simple subcoalgebra corresponding to $V$; that is, if and only if, every irreducible left coideal of $H$ appears with multiplicity $\dim V$ or 0 in $H^{\co B}$. The proposition follows from Lemma \ref{idf-coinv}. \end{proof}

\chapter{The Nichols-Richmond Theorem}\label{dos}

Recall from \ref{irr-char} that the \emph{character algebra} of $H$, denoted $R(H)$, is the subalgebra of $H^*$ spanned by the irreducible characters of $H$.

A subalgebra $M$ of $R(H)$ is called a \emph{standard}
subalgebra if $M$ is spanned by irreducible characters of $H$. So
that if $X$ is a subset of $\widehat H$, $X$ spans a standard
subalgebra of $R(H)$ if and only if the product of characters in
$X$ decompose as a sum of characters in $X$.

By \cite[Theorem 6]{NR} there is a bijection between standard
subalgebras of $R(H)$ and quotient Hopf algebras of $H$. Under
this bijection, the quotient $H \to B$ corresponds to the
character algebra of $B$: $R(B) \subseteq R(H)$.

\section{Irreducible characters of degree 2}\label{6-8}
In this section we collect some facts related to the fusion
rules of irreducible characters of degree $2$.

Suppose that $H$ contains an irreducible character $\chi$ of
degree $2$, such that
\begin{equation}\label{dim8} \chi^2 =  \sum_{g \in G[\chi]} g;\end{equation}
in particular, $\chi^* = \chi$ and $|G[\chi]| = 4$. The set
$G[\chi] \cup \{ \chi \}$ spans a standard subalgebra of $R(H)$.
This subalgebra corresponds to a quotient Hopf algebra $H \to
\overline H$, where $\overline H$ is a semisimple non-commutative
Hopf algebra of dimension $8$ with algebra type $(1, 4; 2, 1)$.
The classification of semisimple Hopf algebras of dimension $8$ implies that the group $G[\chi] = G({\overline H}^*)$ is not cyclic; a more general
picture will appear later in \ref{str-G}.

Conversely, every semisimple non-commutative Hopf algebra of
dimension $8$ has four linear characters, which constitute a group
isomorphic to $\mathbb Z_2 \times \mathbb Z_2$, and one
(self-dual) irreducible character of degree $2$, satisfying the
relations \eqref{dim8}.

Suppose now that $H$ has an irreducible character $\chi$ of degree
$2$ such that \begin{equation}\label{dim6} \chi \chi^* = \epsilon
+ g + \chi,\end{equation} for some $g \in G(H^*)$. Then $G[\chi] =
\{ \epsilon, g \}$, $\chi^* = \chi$ and $G[\chi] \cup \{ \chi \}$
spans a standard subalgebra of $R(H)$, which corresponds to a
quotient Hopf algebra $H \to \overline H$, where $\overline H
\backsimeq k \mathbb S_3$ is the unique non-commutative semisimple
Hopf algebra of dimension $6$.

\begin{proposition}\label{cociente8} Suppose that the following conditions are
fulfilled:

 (i) $\vert \{ \chi \in \widehat H: \chi(1) = 2 \} \vert$ is odd;

 (ii) $G(H^*)$ contains a subgroup $\Gamma$ of order $4$.

Then  there is a quotient Hopf algebra $H \to \overline H$, where
$\overline H$ is a semisimple non-commutative Hopf algebra of dimension
$8$ such that $\Gamma \subseteq {\overline H}^*$. \end{proposition}

\begin{proof} The group $\Gamma$ acts on the set
$X_2 =  \{ \chi \in \widehat H: \chi(1) = 2 \}$ in the form $g .
\chi = g \chi g^{-1}$. Let $X_2' \subseteq X_2$ be the set of
fixed points under this action. Since  $|X_2|$ is odd by
assumption, then $|X_2'|$ is also odd. Moreover, since $\Gamma$ is
abelian, $\Gamma$ acts on $X_2'$ by left multiplication.

Let $Y \subseteq X_2'$ be the set of fixed  points of $X_2'$ under
left multiplication by $\Gamma$. Once again we find that $|Y|$ is
odd and, in particular, that $Y$ is not empty. It is easy to see that
$\mu^* \in Y$ for all $\mu \in Y$, whence there must exist $\chi
\in Y$ such that $\chi^* = \chi$.

By construction,  $\chi^2 = \chi \chi^* = \sum_{g \in \Gamma} g$; see \eqref{des-ss}.
Hence the set $\Gamma \cup \{ \chi \}$ spans a standard subalgebra
of $R(H)$ which corresponds to a Hopf algebra quotient of
dimension $8$ as claimed. \end{proof}

\begin{lemma}\label{dual-estable}
Suppose that  $\lambda \lambda^* \in kG(H^*)$ for some irreducible
character $\lambda$. Assume in addition that $\deg \lambda \leq
\deg \mu$ for all irreducible character $\mu$ with $\deg \mu > 1$.
Then also $\lambda^* \lambda \in kG(H^*)$. \end{lemma}

\begin{proof} By assumption we have $\lambda \lambda^* = \sum_{g \in G[\lambda]} g$.
On the other hand, we may write $\lambda^* \lambda = \sum_{g \in
G[\lambda^*]} g + \sum_{\deg \mu > 1} n_{\mu} \mu$, where $n_{\mu}
= m(\mu, \lambda^* \lambda)$. Then we have $$\lambda \lambda^*
\lambda = \sum_{g \in G[\lambda^*]} \lambda g + \sum_{\deg \mu >
1} n_{\mu} \lambda \mu  = |G[\lambda^*]| \lambda + \sum_{\deg \mu
> 1} n_{\mu} \lambda \mu,$$ and also $$\lambda \lambda^* \lambda =
\sum_{g \in G[\lambda]} g \lambda = |G[\lambda]| \lambda.$$
Comparing the multiplicity of $\lambda$ in both expressions, we
find that $\lambda \mu = (\deg \mu) \lambda$, for all $\mu$ such
that $n_{\mu} \neq 0$. Thus $\deg \mu = m(\lambda, \lambda \mu) =
m(\mu, \lambda^* \lambda) = n_{\mu}$. But then $n_{\mu} \deg \mu \geq
(\deg \lambda)^2$ for all $\mu$ such that $\deg \mu > 1$ and $n_{\mu} \neq 0$; so we
see that $n_{\mu} = 0$, for all $\mu$ such that $\deg \mu > 1$.
This finishes the proof of the lemma. \end{proof}

\section{The Nichols-Richmond theorem} The following theorem is due to Nichols and
Richmond. See \cite[Theorem 11]{NR}. We state here a version
convenient to the finite dimensional context.

\begin{theorem}\label{thm-nr} Suppose that $H$ has an irreducible character
$\chi$ of degree $2$. Then at least one of the following conditions holds:

 (i) $G[\chi] \neq 1$;

(ii) $H$ has a Hopf algebra quotient  of dimension
$24$, which has a degree one character $g$ of order $2$ such that
$g \chi \neq \chi$;

 (iii) $H$ has a Hopf algebra quotient of dimension
$12$ or $60$. \qed \end{theorem}

Notice that, as a consequence of the theorem, if $\dim H < 60$ and
$H$ has an irreducible character of degree $2$, then $G(H^*) \neq
1$.

\begin{remark}\label{el-nr} (i) Suppose that $H$ has an irreducible character $\chi$ of
degree  $2$ such that  $G[\chi] = 1$. Then $H$ must also contain
an irreducible character $\psi$ of degree $3$, which is
necessarily self-dual, and such that  $\chi \chi^* = \epsilon +
\psi$.

Assume that  $H$ has no irreducible character of degree $4$. Then $|G[\psi]|
= 3$ and $\psi^2 = \sum_{g \in G[\psi]}g + 2 \psi$; so that  $G[\psi] \cup \psi$ spans a standard subalgebra of $R(H)$,
which corresponds to a quotient Hopf algebra of type $(1, 3; 3, 1)$.

\begin{proof} We follow the lines of the proof of Case 1 in \cite[Theorem 11]{NR}. Since $H$ has no irreducible character of degree $4$ and $m(\chi, \psi \chi) = 1$, a counting argument implies that $\psi \chi$ is a sum of irreducible characters of degree 2. Moreover, it is easy to see that all these characters are conjugated to $\chi$. By \cite[Theorem 10 (3)]{NR}, $\psi \chi = \sum_{g \in G[\psi]} g \chi$, so in particular $|G[\psi]| = 3$. Multiplying both sides of this equation  by $\chi^*$ on the right, we get $\psi^2 = \sum_{g \in G[\psi]}g + 2 \psi$. This proves the claim. \end{proof}

(ii)  The assumption $G[\chi] = 1$, for an irreducible character
$\chi$ of degree $n$, implies also that $\chi$ has exactly
$|G(H^*)|$ distinct conjugates under the action of $G(H^*)$ by
left multiplication. Thus, in this case, we must have an
inequality $|G(H^*)| \leq |X_n|$.
\end{remark}

The following corollary of the Nichols-Richmond Theorem gives some restrictions for the possibility  $G(H^*) = 1$.

\begin{corollary}\label{cor-g-1} Suppose that $H$ has an
irreducible character of degree 2.
If  $G(H^*) = 1$, then $H$ has a  Hopf algebra quotient of
dimension 60.  In particular, $60 / \dim H$. \end{corollary}

\begin{proof}  Any semisimple Hopf algebra of dimension $12$
contains nontrivial group-like elements \cite{fukuda}.
By Theorem \ref{thm-nr}, if $G(H^*) = 1$,  $H$ must have a  Hopf algebra quotient of dimension 60 as claimed. \end{proof}

\section{An application to $D(H)$} In this section, we consider an application of Theorem \ref{thm-nr} to the Drinfeld double of
$H$.

\begin{proposition}\label{doble} Suppose that $H$ is a  semisimple
Hopf algebra such that $\dim H < 60$.  If the Drinfeld double
$D(H)$ has an irreducible character of degree $2$, then $G(D(H)^*)
\neq 1$.   \end{proposition}

In particular,  $D(H)$ is not simple; see \cite[Corollary 2.3.2]{pqq}.

\begin{proof} If $D(H)$ has an irreducible character of degree $2$, then by
Corollary \ref{cor-g-1}, $G(D(H)^*)
\neq 1$, unless there is a quotient Hopf algebra $q: D(H) \to B$,
with  $\dim B = 60$ such that $G(B^*) = 1$.

If this were the case, then by \cite{NZ} $\dim H$ is divisible by
$30$  (since it must be divisible by $2$, $3$ and $5$); thus by
assumption, $\dim H = 30$. The proposition follows from Theorem
\ref{cls30}, which says that $H$ is necessarily trivial.
\end{proof}

In our applications we shall combine the preceding proposition with the following fact.

\begin{lemma}\label{index-3}
Suppose that $H$ has a Hopf subalgebra or quotient Hopf algebra of
index $3$.  If $G(H^*) \cap Z(H^*) = 1$ and $G(H) \cap Z(H) = 1$,
then the Drinfeld double of $H$ contains an irreducible character
of degree $2$. \end{lemma}

\begin{proof} We may assume that $H$ has a quotient Hopf algebra $H \to B$ of index $3$;
the other case is dual,  once we notice that $D({H^*}^{\cop})
\simeq D(H)^{\op} \simeq D(H)$. By assumption $\dim H^{\co B} = 3$, and by
Lemma \ref{YD-coinv} $H^{\co B}$ is a Yetter-Drinfeld submodule of $H$. Decomposing $H^{\co B}$ as a direct sum of irreducible Yetter-Drinfeld submodules implies the lemma, since the trivial appears with multiplicity 1. \end{proof}

\section{Existence of  proper Hopf subalgebras}
In this section we apply the Nichols-Richmond theorem in order to assure,
in certain cases, the existence of proper Hopf subalgebras.

\begin{lemma}\label{productodesimples} Let $\chi$ and $\psi$  be irreducible
characters of $H$ such that the product $\chi \psi$ is
irreducible. Then for all irreducible character  $\mu \neq
\epsilon$ we must have $m(\mu, \psi \psi^*) = 0$ or $m(\mu, \chi^*
\chi) = 0$.

In particular $G[\psi] \cap G[\chi^*] = 1$. \end{lemma}

\begin{proof} Let $\zeta = \chi \psi$. By Schur's Lemma, $\zeta$
is irreducible if and only if $m(\epsilon, \zeta \zeta^*) = 1$.

On the other hand, we have
\begin{align*} \zeta \zeta^*  & = \chi \psi \psi^* \chi^*
= \chi  \left( \sum_{\mu \in \widehat H} m(\mu, \psi \psi^*) \mu
\right) \chi^*   \\ & = \chi \chi^* +  \sum_{\mu \neq \epsilon} m(\mu,
\psi \psi^*) \chi \mu \chi^*. \end{align*} Therefore,  $m
(\epsilon, \zeta \zeta^*) = 1$ if and only if  for all $\mu \neq
\epsilon$, with $m(\mu, \psi \psi^*) > 0$, we have
$m(\epsilon, \chi \mu \chi^*) = 0$ or equivalently, $m(\mu, \chi^*
\chi) = 0$. This proves the lemma. \end{proof}

\begin{theorem}\label{corolario} Suppose that  $1 \neq G[\chi] \cap G[\psi]$,
for all irreducible characters $\chi$ and $\psi$
of degree $2$. Then there is a quotient Hopf algebra $\pi: H  \to \overline H$,
such that $\overline H$ is of type $(1, |G(H^*)|; 2, |X_2|)$ as an algebra. \end{theorem}

\begin{proof} It follows from Lemma \ref{productodesimples} and the assumptions, that
if $\chi$ and $\psi$ are irreducible characters of degree $2$,
then their product $\chi \psi$ is not irreducible. Thus $\chi \psi$
decomposes as a sum of irreducible characters of degree at most
$2$; indeed, the one-dimensional characters appearing in $\chi
\psi$ with positive multiplicity form a coset of the stabilizer
$G[\chi]$ in $G(H^*)$ and thus there are $0$, $2$ or $4$  of them, so that the
other irreducible summands must be of degree $2$. Therefore, the
set $\{ \chi \in \widehat H: \chi(1) \leq  2 \}$ spans a standard
subalgebra of the character algebra of $H$, implying  the claim.
\end{proof}

We next collect some useful variations of this theorem.

\begin{remark}\label{nr} (i) If $H$ has no irreducible character of degree $4$ and
the  dimension of $H$ is not divisible by $12$, then the set
$G(H^*) \cup X_2$ spans a standard subalgebra of $R(H)$.
Therefore,
$H$ has a quotient Hopf algebra of type $(1, |G(H^*)|; 2,
|X_2|)$.

\begin{proof} In this case the product of two characters of
degree $2$ cannot  be irreducible. By Theorem \ref{thm-nr}, since $12$ does not
divide $\dim H$, then $G[\chi] \neq 1$ for all $\chi \in X_2$,
implying that for all $\psi \in X_2$, $\chi \psi$ decomposes as a
sum of irreducible characters of degree $\leq 2$ (arguing as in
the proof of Theorem \ref{corolario}). \end{proof}

(ii) Suppose that $G[\chi] \neq 1$ for all $\chi \in
X_2$.  Assume in addition that $G(H^*)$ has a unique subgroup $F$
of order $2$. Then all irreducible characters of degree $2$ are
stable under multiplication by $F$.

In particular $F \subseteq
G[\psi] \cap G[\chi^*]$, for all irreducible characters $\psi$ and
$\chi$ of degree $2$.

(iii)  Suppose that the action by right multiplication
of $G(H^*)$ on $X_2$ is transitive and $|G[\chi]| \neq 1$ for some
irreducible character of degree $2$. Then the set $G(H) \cup X_2$
spans a standard subalgebra of $R(H^*)$ which corresponds to a
quotient Hopf algebra $H \to B$, such that $B$ is of type $(1,
|G(H^*)|; 2, |X_2|)$ as an algebra.

\begin{proof} The assumption implies that $|G[\chi]| \neq 1$, for all irreducible character $\chi$ with $\deg \chi = 2$: indeed, the  action of $G(H^*)$ by right multiplication, which is transitive, preserves
$G[\chi]$ and this group is not trivial at least for one $\chi$. Therefore $\chi \chi^*$ belongs to the span of $G(H^*)$ and $X_2$.

Since the action of $G(H^*)$ on the set $X_2$ by right
multiplication is transitive, for all $\chi' \in X_2$, there
exists $h \in G(H^*)$ such that $\chi' = \chi^* h$; then $\chi \chi'
= \chi \chi^*h$ belongs to the span of $G(H^*) \cup X_2$.

This shows
that the set $G(H^*) \cup X_2$ spans a standard subalgebra of $R(H)$ and
implies the claim. \end{proof} \end{remark}

\section{Hopf subalgebras of index 3} Using character theory, we shall show that for some Hopf algebra inclusions $A \subseteq H$, with $[H:A] = 3$, there
is a quotient Hopf algebra $H \to \overline H$, such that the simple $\overline H$-modules have dimension 1 or 2.

Observe that, by Corollary \ref{kob-mas}, if $H$ has a Hopf subalgebra of index 3 which is not normal, then the dimension of $H$ must be even; hence $6 / \dim H$. The following theorem and its corollary give more precise information.

\begin{theorem}\label{coinvariantes} Let $A \subseteq H$ be a Hopf subalgebra such that $[H:A] = 3$.
Suppose  that $A$ is not normal in $H$.
Then there exist a Hopf algebra $\overline H$ of algebra type
$(1, |G(H^*)|; 2, \dfrac{|G(H^*)|}{2})$ together with surjective Hopf algebra maps $H \to \overline H \to B$, where $B \simeq k \mathbb S_3$, and a surjective coalgebra map $B \to H/HA^+$ such that the following diagram is commutative:
$$\begin{CD}H @>>> \overline H \\
@VVV @VVV \\
H/HA^+ @<<< B. \end{CD}$$ \end{theorem}

\begin{proof} Since $A$ is not normal in $H$, by Proposition \ref{ind-norm}, $\Ind_A^H\epsilon = \epsilon + \chi$, where $\chi$ is an irreducible character of degree 2. In particular, $\chi^* = \chi$,  by Lemma \ref{idf-coinv}.

\begin{claim}\label{af-s3} $|G[\chi]| = 2$ and $G[\chi] \cup \{ \chi \}$ spans a standard subalgebra of $R(H)$.  \end{claim}

\begin{proof} It is enough to show that the product $\chi \chi^* = \chi^2$ admits a decomposition as in \eqref{dim6}.

Since $m(\chi, \Ind_A^H\epsilon) = 1$, we have $\chi\vert_A = \epsilon + x$, where $\epsilon \neq x \in G(A^*)$, by Frobenius reciprocity. Hence $\chi^2\vert_A = 2(\epsilon + x)$.

In the character algebra of $H$ one of the following decompositions must hold:

(a) $\chi^2 = \epsilon + a + b + c$, where $a, b, c \in G(H^*) \backslash \{ \epsilon \}$ are pairwise distinct;

(b) $\chi^2 = \epsilon + \psi$, where $\psi$ is an irreducible character of degree 3;

(c) $\chi^2 = \epsilon + a + \lambda$, where $\epsilon \neq a \in G(H^*)$ and  $\lambda$ is an irreducible character of degree 2.

We shall show that the cases (a) and (b) are impossible, thus proving that $|G[\chi]| = 2$. Suppose that (b) holds. Restricting to $A$, we find that necessarily $m(\epsilon, \psi\vert_A) > 0$. But by Frobenius reciprocity $m(\epsilon, \psi\vert_A) = m(\psi, \Ind_A^H \epsilon) = 0$; this contradiction discards case (b). Case (a) is similarly discarded.

Hence (c) holds. It remains to show that $\lambda = \chi$. For this, we restrict the equation $\chi^2 = \epsilon + a + \lambda$ to $A$, and apply the Frobenius reciprocity to find that $m(\lambda, \Ind_A^H \epsilon) > 0$, whence $\lambda = \chi$. This proves the claim. \end{proof}

Clearly, the standard subalgebra in Claim \ref{af-s3} corresponds to a quotient Hopf algebra $H \to B$, such that $B \simeq k\mathbb S_3$. Moreover, by construction, and using Lemma \ref{idf-coinv}, $(H^*)^{\co A^*} \subseteq B^* \subseteq H^*$. Hence there is a surjective coalgebra map $B \to H/HA^+$ which factorizes the canonical map $H \to H/HA^+$.

Lemma \ref{idf-coinv} implies also that $g \chi g^{-1} = \chi$,
for all $g \in G(H^*)$.  Therefore, $G(H^*) \cup G(H^*)\chi$ spans
a standard subalgebra of $R(H)$, which corresponds to a quotient
Hopf algebra $H \to \overline H$ with the desired properties. This
finishes the proof of the theorem.  \end{proof}

Note that $\dim \overline H = 3 \vert G(H^*) \vert$. The following corollary is an immediate consequence of the theorem.

\begin{corollary}\label{cor-ind3} Suppose that $A \subseteq H$ is a Hopf subalgebra such that $[H:A] = 3$ and $A$ is not normal in $H$.
Then we have

(i) $|G(H^*)|$ is even;

(ii) $3|G(H^*)|$ divides $\dim H$. \qed \end{corollary}

\chapter{Quotient Coalgebras}\label{tres}

Let $H$ be a semisimple Hopf
algebra and let $A \subseteq H$ be a Hopf subalgebra.  Consider
the quotient coalgebra $p: H \to \overline H : = H / HA^+$. By
\cite[3.4]{masuoka-freeness}, $\overline H$ is a cosemisimple
coalgebra. In this chapter we aim to relate the corepresentations
of $H$ and $\overline H$. We discuss the
corepresentation theory of $\overline H$ in relation with that of
$H$ and the corestriction functor ${}^H\mathcal M \to {}^{\overline
H}\mathcal M$.

We show that if $C$ is a simple  subcoalgebra of $H$ such that $C
a \subseteq C$, for all $a \in A$, then the dual of the quotient
coalgebra $C / CA^+$ and the crossed product $A_{\alpha}$, where
$\alpha : A \otimes A \to k$ is a certain $2$-cocycle, constitute
a commuting pair in $C^*$. This is applied in combination with
Masuoka's main result in \cite{masuoka-cont}, in some instances of
the proof of Theorem \ref{th-ss}.

In particular, when $A = kG$ is the group algebra of a subgroup
$G$ of $G(H)$ and  $V$ is a simple $H$-comodule, we deduce that
${\rm End }^{\overline H} (V)$ is isomorphic as an algebra to a twisted
group algebra $k_{\alpha}\Gamma$, where $\Gamma \subseteq G$ is
the stabilizer of $V$,  {\it i.e.} $\Gamma = \{ g \in G: V \otimes
g \simeq V \}$, and $\alpha : \Gamma \times \Gamma \to k^{\times}$
is a $2$-cocycle. This result implies that the multiplicity of an
irreducible $\overline H$-comodule in $V$ is a divisor of the
order of $\Gamma$. In particular, when the group $\Gamma$ is
abelian,  all irreducible $\overline H$-comodules in the
restriction of $V$ to $\overline H$ appear with the same
multiplicity $d$, where $d$ divides the order of $\Gamma$. This
allows us to recover the result in \cite[Proposition 2.4]{masuoka-cont}.

Some of these results are applied to the case
when $H$ is a biproduct in the  sense of Radford: $H \simeq R \#
A$. Indeed, in this case $R$ is isomorphic as a
coalgebra to the quotient $H / HA^+$.

We shall denote by ${}^H\mathcal M$ the category of left
$H$-comodules; for a left $H$-comodule algebra $A$, the category
of (left-right) $(A, H)$-Hopf modules will be indicated by ${}^H\mathcal M_A$.

\section{A multiplicity formula} Let $\rho: V \to H \otimes V$,
$\rho (v) = v_{-1} \otimes v_0$, be a left $H$-comodule.  Consider
the  $\overline H$-comodule structure on $V$ obtained by
corestriction along the coalgebra map $p$; we shall sometimes use
the notation $\overline V$ for this comodule structure, as well as
$\overline{\rho} = (p \otimes \id) \rho: V \to \overline H \otimes
V$ for the structure map. We obtain in this way a functor
${}^H\mathcal M \to {}^{\overline H}\mathcal M$, $V \mapsto
\overline V$.

By a result of Schneider  \cite[Theorem II]{sch-ppal}, there is a
category equivalence $\omega: {}^H\mathcal M_A \to {}^{\overline
H}\mathcal M$ between the category of $(A, H)$-Hopf modules and
the category of $\overline H$-comodules. The equivalence is given
by $\omega: M \mapsto M/ MA^+$, for any  object $M$ of $ {}^H\mathcal
M_A$.

Consider now the functor  $F: {}^H\mathcal M \to
{}^H\mathcal M_A$, $F(V) = V \otimes A$; where  $H$ coacts
diagonally on $V \otimes A$ and the right $A$-module structure is
given by $v \otimes a . b = v \otimes ab$.

\begin{lemma}\label{adjointness}
(i) The functor $F$ is a left adjoint of the forgetful  functor
$U: {}^H\mathcal M_A \to {}^H\mathcal M$;

 (ii) for all left $H$-comodules $V$
there is an isomorphism $\omega F(V) \simeq \overline V$.
\end{lemma}

\begin{proof} (i) We define natural maps $$\psi: \Hom^H (X, U(Y)) \to
\Hom^H_A(X \otimes A, Y),$$
$$\phi:\Hom^H_A(X \otimes A, Y) \to
\Hom^H (X, U(Y)),$$ as follows:
\begin{equation*} \psi (f) (x \otimes a) = f(x).a, \quad
\phi(g) (x) = g(x \otimes 1), \end{equation*} for all $x \in X$,
$a \in A$. It is not hard to check that $\psi$  and $\phi$ are
well defined and are indeed inverse isomorphisms.

(ii)  The maps $\mu: (V \otimes A) / (V \otimes A^+) \to \overline
V$, $\mu ([v \otimes a]) : = \epsilon (a) v$, and $\eta: \overline
V \to (V \otimes A) / (V \otimes A^+)$, $\eta (v) : = [v \otimes
1]$, define inverse $\overline H$-colinear  isomorphisms.
\end{proof}

As a consequence we get the following proposition.

\begin{proposition}\label{adj} Suppose that  $U$ and $V$ are
finite-dimensional left $H$-comodules. There is a natural linear
isomorphism
$$\Hom^{\overline H} (\overline U, \overline V) \simeq
\Hom^H(V^* \otimes U, A).$$
Suppose that $V$ is a simple $H$-comodule.  Then $\overline V$ is simple if and only if
$$\Hom^H(V, V \otimes W) = \Hom^H(V^* \otimes V, W) = 0,$$ for all
simple left $A$-comodule $k1 \neq W$. \end{proposition}

\begin{proof} We  have isomorphisms
\begin{align*}\Hom^H (U, V \otimes A) & \simeq \Hom^H_A(U
\otimes A, V \otimes A) \\ & \simeq \Hom^{\overline H}(\omega (U \otimes A), \omega (V \otimes A))  \simeq \Hom^{\overline H} (\overline U,
\overline V);\end{align*}
the first isomorphism by Lemma \ref{adjointness} (i), the second since $\omega$ is a category equivalence and the third by Lemma \ref{adjointness} (ii).
This proves the first statement since $\Hom^H (U, V
\otimes A) \simeq \Hom^H(V^* \otimes U, A)$.

As left $H$-comodule, $A$ decomposes in the form $A \simeq \oplus_{W} (\dim W) W$, where the sum runs over the set of isomorphism classes of irreducible left coideals  of $A$, which coincides with the set of isomorphism classes of irreducible left coideals $W$ of $H$ which are contained in $A$. Hence we have $$\Hom^{\overline H} (\overline V, \overline V) \simeq \Hom^H(V^* \otimes V, A) \simeq \oplus_{W} (\dim W) \Hom^H(V^* \otimes V, W);$$ this implies the last statement, in view of Schur's Lemma. \end{proof}

\begin{remark}\label{dimension} (i) Suppose that $U$ and  $V$ are finite-dimensional left $H$-comodules,
and let $\chi_U$ and $\chi_V \in H$ be the corresponding
characters. As a left $H$-comodule, $A \simeq \oplus_{\lambda \in
\widehat{A^*}}  \deg \lambda \ W_{\lambda}$.
Proposition \ref{adj} implies the following multiplicity formula:
$$\dim \Hom^{\overline H} (\overline U,
\overline V) = \sum_{\lambda \in \widehat{A^*}} \deg \lambda \
m(\lambda, \chi_V^* \chi_U).$$

 (ii)  Suppose that $A = kG$, where $G$ is a subgroup
of $G(H)$. Let $G[V^*]$ denote the subgroup of $G$ consisting of all elements $g$ for which $Vg \simeq V$; that is, $G[V^*] = G \cap G[\chi_V^*]$. Recall from \eqref{mult-g} that for $g \in G$, we have $\dim \Hom^H (V^* \otimes V, g) = 1$ if and only if $g \in G[V^*]$, and $\dim \Hom^H (V^* \otimes V, g) = 0$ otherwise.

It follows from Proposition \ref{adj} that $$\End^{\overline H}V \simeq \Hom^H
(V^* \otimes V, kG) = \bigoplus_{g \in G[V^*]}
\Hom^H (V^* \otimes V, g).$$ We thus get $\dim \End^{\overline H} V = |G[V^*]|$. \end{remark}

\section{Stable subcoalgebras}
We keep the notation  in the previous section. Let $C \subseteq
H$ be a simple subcoalgebra and let $V \subseteq C$ be an
irreducible left coideal of $H$. Suppose that  $Ca \subseteq  C$,
for all $a \in A$. That is, $C$ is a  right $A$-module with
action given by right multiplication; by \cite{NZ} $C$ is a free
right  $A$-module.  Let $\overline C : = p(C)$ and let $t \in A$
be the normalized integral. We have $H = H (1 - t) \oplus Ht$, and
$H A^+ = H (1 - t)$. In particular,  $p(h) = p(ht)$, for all $h
\in H$.

Notice also that  $\End^H V = \End^CV$ and $\End^{\overline H} V =
\End^{\overline C}V$.

\begin{lemma}\label{dim-cbar}
(i) The map $p : H \to \overline H$  induces  an identification
$C / C(1 - t) = \overline C$;

 (ii) $\dim \overline C = (\dim A)^{-1} \dim C$.
\end{lemma}

\begin{proof} (i) The map $p$ induces an identification $\overline C = C / C
\cap HA^+ = C / C \cap H(1-t)$. We claim that $C \cap H(1-t) = C
(1-t)$. Indeed, let $c \in C \cap H(1-t)$; then  $c = h (1-t)$
implying, since $(1-t)$ is an idempotent, that $c = c(1-t) \in
C(1-t)$. The other inclusion is immediate from the fact that $Ca
\subseteq C$, for all $a \in A$.

 (ii) We have $C = C(1-t) \oplus Ct$, since $t$ and
$1-t$ are orthogonal idempotents, and $C$ is stable under right
multiplication by $A$. By part (i), $\dim \overline C = \dim C t$.
But $Ct = C^A$ is the space of $A$-invariant elements in $C$ under
the action by right multiplication. Since $C$ is a free right
$A$-module of rank  $(\dim A)^{-1} \dim C$, then $\dim C^A = (\dim
A)^{-1} \dim C$.  \end{proof}

\begin{remark}\label{fusion-rules} Let  $\chi \in C$ be the irreducible character of $H$ corresponding to $C$.
The simple subcoalgebra $C$ satisfies  $Ca \subseteq C$ for all $a
\in A$ if and only if $\chi \psi \in \mathbb Z \chi$, for all
$\psi \in \widehat{A^*}$. This follows from the fact that the
multiplication map $m: H \otimes H \to H$ is a left (and right)
$H$-comodule map. \end{remark}

Since $C$ is a right $A$-module
coalgebra under right multiplication, then  $C^*$ is a
left $A$-module algebra  under the action $(a . f) (c) = f (ca)$,
$f \in C^*$, $c \in C$, $a \in A$.

By the Skolem-Noether Theorem for Hopf algebras \cite{ma-sk-nt},
since $C^*$ is a simple algebra over $k$, there exists a
convolution invertible map $\psi: A \to C^*$ such that $a. f =
\psi (a_1) f \psi^{-1} (a_2)$, for all $f \in C^*$, $a \in A$.
This gives  rise to an algebra map $\psi: A_{\alpha} \to C^*$,
where $\alpha \in Z^2(A, k)$ is the $2$-cocycle associated to
$\psi$ in the form $\alpha (a, b) = \psi^{-1}(a_1b_1) \psi(a_2)
\psi(b_2)$. Here, and elsewhere,  $A_{\alpha} : = A
\#_{\alpha} k$ denotes the associated crossed product with respect
to the cocycle $\alpha$; that is, $A_{\alpha} = A$ as vector spaces, with the multiplication $a.b = \alpha(a_2, b_2) a_1b_1$, $a, b \in A$.
Note that, when $A = kG$ is a group
algebra, $(kG)_{\alpha} = k_{\alpha}G$ is the \emph{twisted}
group algebra.

\begin{proposition}\label{comm-pair} As algebras, $\left(\End^{\overline C} V\right)^{\op} \simeq A_{\alpha}$. Moreover,
$A_{\alpha}$ and ${\overline C}^*$  form a commuting pair in
$(\End V)^{\op} \simeq C^*$. \end{proposition}

\begin{proof} The coalgebra projection $p: C \to \overline C$ induces by
transposition an algebra inclusion ${\overline C}^* \subset C^*$.
We have that ${\overline C}^*$ coincides
with the subalgebra of invariants $(C^*)^A$ under the action of
$A$. Indeed, $f \in {\overline C}^*$ if and only if $f(c) =
f(p(c))$, for all $c \in C$, if and only if $f(c) = f(ct) = (t.f)
(c)$, for all $c \in C$; whence, ${\overline C}^* = (C^*)^A$.

Hence, by definition of $\psi$, ${\overline C}^*$ coincides with
the commutant of $\psi(A_{\alpha})$ in $C^*$. Since
${\overline C}^*$ is semisimple, the double commutant theorem
implies that $\psi(A_{\alpha})$ is the commutant of ${\overline
C}^*$ in $C^*$.

It follows from Burnside's Density Theorem \cite{CR} that there
is an anti-isomorphism of algebras $I: C^* \to \End V$, given by
\begin{equation}\label{c^*} I(f) (v) : = \langle f, v_{-1} \rangle v_0, \quad v
\in V, f \in C^*.\end{equation}
Under this identification, the commutant of ${\overline
C}^*$ in $C^*$ coincides with the subalgebra $\left(\End^{\overline C} V\right)^{\op}$ of $(\End V)^{\op}$.
Therefore $\psi$ defines a surjective algebra map $\psi: A_{\alpha} \to \psi(A_{\alpha}) = \left(\End^{\overline C} V\right)^{\op}$.

Since $C A \subseteq C$, we have $\chi_V \lambda = (\deg \lambda)
\chi_V$, for all irreducible character $\lambda \in
\widehat{A^*}$; that is, $m(\lambda, \chi_V^* \chi_V) = \deg \lambda$, for all $\lambda \in \widehat{A^*}$. By Remark \ref{dimension} (i), $\dim
\End^{\overline C} V = \sum_{\lambda \in \widehat{A^*}} (\deg \lambda)^2  = \dim A$. Therefore $\psi: A_{\alpha} \to
C^*$ is an injective algebra map and determines an isomorphism
$A_{\alpha} \simeq \left(\End^{\overline C} V\right)^{\op}$. This finishes the
proof of the proposition.  \end{proof}

\begin{corollary}\label{comm-pair-cor} There exists a bijective
correspondence between irreducible
$A_{\alpha}$-modules and irreducible $\overline C$-comodules.

There is an isomorphism $V \simeq \bigoplus_i U_i \otimes W_i$, where $W_i$ runs
over a system of representatives of the isomorphism classes of irreducible
$A_{\alpha}$-modules and $U_i$ is the  irreducible $\overline C$-comodule
corresponding to  $W_i$.
 \qed \end{corollary}

In particular, if $A = kG$ where $G$ is a finite group, then the multiplicity of $U_i$ in $V$
divides the order of $G$, for all $i$.

If $G$ is an abelian group, then all the irreducible $k_{\alpha}G$-modules
have the same dimension $d$. Therefore all irreducible $\overline
C$-subcomodules of  $V$ appear with the same multiplicity $d$.

As an application of the methods of this section, we have the following proposition.

\begin{proposition}\label{dima=dimc} Suppose that $AC = C = CA$. Assume in addition that $\dim A = \dim C$. Then $A$ is normal in $k[C]$. \end{proposition}

Here, $k[C]$ denotes the subalgebra generated by $C$; this is a Hopf subalgebra of $H$ containing $A$.

\begin{proof} By Lemma \ref{dim-cbar}, $\dim C t = \dim tC
= 1$, where $t \in A$ is the normalized integral. Therefore
$Ct = tC = k\psi$, where $\psi \in C$ is the
corresponding irreducible character. Hence, for all $c \in C$,
$tc = l(c) \psi$ and $c t = r(c) \psi$, which implies that $tc
= ct$, after applying the counit.

Hence  $t$ commutes pointwise with $C$; thus it is central
in $k[C]$, and {\it a fortiori}, $A$ is normal in $k[C]$. \end{proof}

\begin{remark}\label{particular} It follows from Corollary \ref{comm-pair-cor} that, if
$G$  is cyclic and $|G| = \dim V$, then $\overline C$ is a
cocommutative coalgebra. We thus recover a fact in the proof of
\cite[Proposition 2.4]{masuoka-cont}. \end{remark}

\section{Quotients modulo group-like Hopf subalgebras}\label{general}
Let $G$ be a subgroup of $G(H)$ and let $A = kG$.
We shall now specialize the description in the previous sections.

Let $t : = \dfrac{1}{|G|} \sum_{g \in G} g$ be the normalized integral in $kG$.

\begin{lemma}\label{intersection} Let $C$ and $D$ be simple subcoalgebras of
$H$. Then the following are equivalent:

 (i) $p(C) \cap p(D) \neq 0$;

 (ii) $p(C) = p(D)$;

 (iii) There exists $g \in G$ such that $Cg = D$.\end{lemma}

\begin{proof} (ii) $\Longrightarrow$ (i). Clear, since $p(C) \neq 0$ for
all simple subcoalgebra $C$.

 (iii) $\Longrightarrow$ (ii). If $D = Cg$, then $Dt = Ct$ and therefore $p(D) =
p(Dt) = p(Ct) = p(C)$.

 (i) $\Longrightarrow$ (iii). Suppose that $c \in C$ and $d \in D$ are
such that $p(c) = p(d) \neq 0$. Then  $c-d$ belongs to the kernel of $p$, and there
exists $h \in H$ such that $c-d = h(1-t)$. Hence $$(c-d)t = h(1-t) t = 0,$$ implying that $ct = dt$. But $ct \in
\sum_{g \in G} C g$ and $dt \in \sum_{g \in G} Dg$. Therefore,
$Cg \cap Dg' \neq 0$, for some $g, g' \in G$; since both $Cg$ and
$Dg'$ are simple subcoalgebras, this implies that $Cg = Dg'$ and
$D = Cg(g')^{-1}$. \end{proof}

The group $G$ acts on the set of simple subcoalgebras of $H$ by right
multiplication. Let $C_1, \dots, C_n$ be a system of representatives of this
action, and let $G_i \subseteq G$ be the stabilizer of $C_i$.

\begin{corollary}\label{cor-cociente}  There is an isomorphism of  coalgebras $\overline H \simeq
\bigoplus_{i = 1}^n \overline{C_i}$, where
$\overline C_i = p(C_i) \simeq C_i / C_i (kG_i)^+$.  \end{corollary}

\begin{proof} It follows from Lemma \ref{intersection} that $\overline H \simeq
\bigoplus_{i = 1}^n p(C_i)$. Thus it remains to see that $p (C_i) \simeq C_i /
C_i (kG_i)^+$.

Fix $1 \leq i \leq n$ and let $C = C_i$, $G_C = G_i$.
We claim that $C \cap  H(kG)^+ = C \cap C(1-t) = C \cap C(1-t_C)$,
where $t_C = \dfrac{1}{|G_C|}\sum_{h \in G_C}h$ is the normalized
integral in $kG_C$.

Note first that if $c \in C \cap H(kG)^+ = C \cap H(1-t)$ then $c
= h(1-t)$  and thus $c(1-t) = h(1-t)^2 = c$; this shows that $C
\cap H(kG)^+ = C \cap C(1-t)$. On the other hand, we have $$t =
\dfrac{1}{|G|}\sum_{g \in G}g = \dfrac{1}{|G|}\sum_{g \in G_C
\backslash G} \sum_{h \in G_C} hg = \dfrac{|G_C|}{|G|}\sum_{g \in
G_C \backslash G} t_Cg.$$ Thus, for  $c \in C$, $ct =
\dfrac{|G_C|}{|G|}\sum_{g \in G_C \backslash G} ct_Cg$ belongs to
$\bigoplus_{g \in G_C \backslash G} Cg$, implying that $ct = 0$ if
and only if $ct_C = 0$. Then $c \in C \cap C(1-t)$ if and only if
$ct = 0$ if and only if $ct_C = 0$ if and only if $c \in C \cap
C(1-t_C) = C(1-t_C)$. This proves the claim and the corollary
follows.
\end{proof}

\section{On the structure of $G(H)$}\label{str-G}
In this section we  present some consequences of Proposition \ref{comm-pair}.
Some  of them are special cases of results in the papers \cite{ty, tambara}.
We keep the notation in the previous sections: $C \subseteq H$ is a simple subcoalgebra, $V \subseteq C$ is an irreducible left coideal, $G$ is a subgroup of $G(H)$ such that $Cg = C$, for all $g \in G$, and $\overline C = C / C(kG)^+$ is the quotient coalgebra.

\begin{proposition}\label{spl-qt}  Suppose that
$\overline C$ is a simple coalgebra.  Then there exists a
non-degenerate $2$-cocycle $\alpha: G \times G \to k^{\times}$. In
particular, the group $G$  is solvable and not cyclic.
\end{proposition}

\begin{proof}  By Corollary \ref{comm-pair-cor} the twisted group algebra
$k_{\alpha}G$ is simple. This implies the proposition. \end{proof}

\begin{corollary} Suppose that $|G| = \dim C = (\dim V)^2$.
Then $G$ is solvable and not cyclic. \qed \end{corollary}

\begin{remark} Suppose as above that $|G| = \dim C$. If we assume in addition that  $V^* \simeq V$, then $A = kG
\oplus C$ is a Hopf subalgebra whose category $\mathcal C$ of corepresentations has necessarily the fusion rules in \cite{ty}.
The results in {\it loc. cit.} imply
that $G$ is abelian. Since by definition $\mathcal C$ admits a fiber functor, the existence of a non-degenerate $2$-cocycle on $G$ is a consequence of
\cite{tambara}. \end{remark}

In the following proposition we give a  Hopf theoretical proof of
these facts concerning $G$, under  rather less restrictive
assumption.

\begin{proposition}\label{ty-hopf} Suppose that the following conditions hold:

(i) $|G| = \dim C$;

(ii) $gC = C = Cg$, for all $g \in G$.

Then the group $G$ is abelian and admits a non-degenerate $2$-cocycle.

If in addition $C = \mathcal S(C)$, then  $A = kG
\oplus C$ is a Hopf subalgebra of $H$ of dimension $2 \dim C$,
which fits into a cocentral extension $1 \to k^{\widehat G} \to A
\to k\mathbb Z_2 \to 0$. \end{proposition}

An analogous result, in the context of Kac algebras, appears in \cite[Theorem IX.8 (ii)]{IK}.

\begin{proof} Keep the notation in the proof of Proposition \ref{comm-pair}.
Let $X_g \in C^*$ denote the image of  $g \in G$ under the map
$\psi : G \to C^*$; so that $\{ X_g: g \in G \}$ is a basis of
$C^*$ and $X_aX_b = \alpha (a, b) X_{ab}$, for all $a, b \in G$.
Also, the action $(a . p) (c) = p(ca)$ is given by $a.p = X_a p
X_a^{-1}$, for all $p \in C^*$, $a \in G$.

In view of condition (ii),  the same arguments apply to the Hopf
algebra $H^{\op}$. Therefore there exists a basis $\{ Y_g: g \in G
\}$ of $C^*$ such that the action $(p . b) (c) = p(bc)$ is given
by $p.b = Y_b^{-1} p Y_b$, for all $p \in C^*$, $b \in G$. The
relation $(a.p).b = a.(p.b)$ implies that $X_aY_bX_a^{-1}Y_b^{-1}
= \zeta (a, b) 1$,  for some  map $\zeta: G \times G \to
k^{\times}$. The definition  of $\zeta$ implies that $a . Y_b =
\zeta (a, b) Y_b$ and similarly that $X_a^{-1} . b = \zeta(a,
b)X_a^{-1}$. Thus, by the associativity of the actions, we get
that $\zeta$ is a bicharacter on $G$.

Suppose that there exists $a \in G$ such that $\zeta(a, b) = 1$,
for all $b \in G$. This implies that $a . Y_b = Y_b$, for all $b
\in G$. Since $\{ Y_b \}_{b \in G}$ is a basis of $C^*$, we get
that the action of $a$ on $C^*$ (and thus on $C$) is trivial. This
implies that $a = 1$, because by \cite{NZ} $C^* \simeq kG$ as a
left and right $kG$-module. Therefore the bicharacter $\zeta: G
\times G \to k^{\times}$  is non-degenerate. This proves that $G$
is abelian as claimed. Finally, since $kG$,  which is isomorphic
to $k^{\widehat G}$, has index  $2$ in $A$, the last part of the
proposition follows; see Corollary \ref{kob-mas}. \end{proof}

\section{A criterion of normality}
We review in this section a result of A. Masuoka, which appears in \cite[Section 2]{masuoka-cont}.

\emph{We shall assume that $C \subseteq H$ is a simple subcoalgebra of dimension $n^2$ of $H$, and $g \in G(H)$ is a group-like element of order $n$ such that $gC = C = Cg$.} In particular, $g \in k[C]$, hence $k\langle g \rangle \subseteq k[C]$.

Let $V \subseteq C$ be an irreducible left coideal, so that $kg \otimes V \simeq V \simeq V \otimes kg$. Let also $\alpha: kg \otimes V \to V$ and $\alpha': V \otimes kg \to V$ be $H$-colinear isomorphisms.

\begin{lemma}\label{l-mas} $k\langle g \rangle$ is a normal Hopf subalgebra of $k[C]$ if and only if $\alpha$ and $\alpha'$ commute as endomorphisms of $V$.
\end{lemma}

\begin{proof} Let $t \in k\langle g \rangle$ be the normalized integral. Then $k\langle g \rangle$ is
normal in $k[C]$ if and only if $tc = ct$, for all $c \in C$. By
\cite[Lemma 2.1]{masuoka-cont} this is  in turn equivalent to
$[\alpha, \alpha'] = 0$ in $\End V$.
\end{proof}

\begin{corollary}\label{cor-mk} Suppose that $V$ is an irreducible left coideal of $C$ such that $gV = V = Vg$.
Then $k\langle g \rangle$ is  normal in $k[C]$.
\end{corollary}

\begin{proof} In this situation, we may take as $\alpha$ the left multiplication by $g$, and $\alpha'$ the right multiplication by $g$. Since this endomorphisms commute with each other, the corollary follows from Lemma \ref{l-mas}.
\end{proof}

\chapter{Braided Hopf Algebras}\label{biprod}

Let $A$ be a semisimple Hopf algebra and let ${}^A_A\mathcal{YD}$ denote the braided category of Yetter-Drinfeld modules over $A$. Let $R$ be a
semisimple braided Hopf algebra in ${}^A_A\mathcal{YD}$. The results in this chapter concern the \emph{biproduct} construction, as described in Section \ref{rm-bip}. This construction was introduced by Radford \cite{Rbip} and interpreted in categorical terms by Majid  \cite{majid-bip, majid-book}.

\section{Radford-Majid biproduct construction}\label{rm-bip}
Denote by $\rho: R \to A \otimes R$, $\rho(a) = a_{-1}
\otimes a_0$, and $. : A \otimes R \to R$, the coaction and
action of $A$ on $R$, respectively.
So that the Yetter-Drinfeld compatibility condition reads as follows:
\begin{equation}\label{yetter-drinfeld} \rho (h . a) =
h_1 a_{-1} \mathcal S(h_3) \otimes h_2 . a_0, \quad \forall a\in R, \ h \in A. \end{equation}
We shall use the notation
$\Delta_R(a) = a^1 \otimes a^2$  and $\mathcal
S_R$ for the comultiplication and the antipode of $R$, respectively.

Thus our assumption amounts to the following conditions:
\begin{flalign}\label{a}
& R \ \text{is an} \ A\text{-module and } \ A\text{-comodule
algebra;} & \\ & \label{b} R \ \text{is an} \ A\text{-module and }
\ A\text{-comodule coalgebra;} & \\ & \label{c} \Delta_R(ab) = a^1
((a^2)_{-1}.b^1) \otimes (a^2)_0b^2; & \\
& \label{d} \mathcal S_R(a^1) a^2 =
\epsilon_R(a)1_R =  a^1 \mathcal S_R(a^2). &
\end{flalign}

Let  $H = R \# A$ be the corresponding biproduct; so
that $H$ is a semisimple Hopf algebra with multiplication,
comultiplication and antipode given by
\begin{equation}\label{mca}
(a \# g) (b \# h)  = a (g_1 . b) \#
g_2h,  \quad \Delta(a \# g)   = a^1 \# (a^2)_{-1} g_1 \otimes
(a^2)_0 \# g_2, \end{equation}
\begin{equation*} \mathcal S (a \# g)  = (1 \# \mathcal S(a_{-1}g)) (\mathcal
S_R(a_0) \# 1), \end{equation*}
for all $g, h \in A$, $a, b \in R$; here we use the notation $a \# g$ to
indicate the element $a \otimes g \in R \# A$. See \cite{Rbip}.

Consider the natural maps $\pi : H \to A$, $\pi (r \# a) = \epsilon_R(r)a$,  and $\iota: A \to H$, $\iota (a) = 1 \otimes a$. Then $\pi$ is a Hopf algebra surjection and $\iota$ is a Hopf algebra injection.
Moreover we have  $\pi \iota = \id_A$ and
\begin{equation}\label{algebraR} R = H^{\co \pi} = \{ h \in H: \ (\id \otimes
\pi) \Delta(h) = h  \otimes 1 \},\end{equation}
coincides, as an $A$-module and $A$-comodule algebra  with the  left  coideal subalgebra of right $A$-coinvariants in $H$.
On the other hand, the map $\id \otimes \epsilon: H \to R$ induces an
isomorphism of left $A$-module and $A$-comodule coalgebras
\begin{equation}\label{coalgebraR} R \simeq H / HA^+.\end{equation}
Indeed, the biproduct construction for finite dimensional Hopf
algebras is   characterized by these properties. Namely, suppose
that there are Hopf algebra maps $\iota: A \to H$ and $\pi: H \to
A$ such that $\pi \iota: A \to A$ is an isomorphism. Then the
subalgebra $R : = H^{\co \pi}$ of right coinvariants of $\pi$  has
a natural structure of Yetter-Drinfeld Hopf algebra over $A$ such
that the multiplication map $R \# A \to H$ induces a Hopf algebra
isomorphism. This principle will be often used throughout this
paper.

Typically, and mainly following the lines described in Chapter \ref{dos}, we shall encounter a Hopf subalgebra $A \subseteq H$ and a surjective Hopf  algebra map  $\pi: H \to B$, where $H$, $A$ and $B$ are certain semisimple Hopf algebras such that $\dim A = \dim B$. After an analysis of the possible left coideal decompositions of $H^{\co B}$, we shall be able sometimes to deduce that $A \cap H^{\co B} = k1$: this property guarantees the injectivity of $\pi \vert_A$, and hence that $\pi\vert_A: A \to  B$ is an isomorphism. This will tell us that $H \simeq R \#A$ has the structure of a biproduct and will enable us to use the biproduct techniques that we discuss in the rest of this chapter.

One simple instance of this situation, frequently used along this paper, is described in the following lemma.

\begin{lemma}\label{spl-ins} Suppose $A \subseteq H$ is a cocommutative Hopf subalgebra and $\pi: H \to B$ is a surjective Hopf  algebra map,  such that $\dim A = \dim B$ and $\dim A$ is relatively prime to $[H: A]$. Then $\pi\vert_A: A \to  B$ is an isomorphism, and $H \simeq R \#A$ is a biproduct. \end{lemma}

\begin{proof} It is enough to show that $A \cap H^{\co \pi} = k1$. This follows from the assumptions, since $\dim A \cap H^{\co \pi}$ divides both $\dim A$ and $\dim H^{\co \pi} = [H: A]$. \end{proof}

\section{Coalgebra structure of $R$}
It turns out that  $R$ is a normal left coideal subalgebra of $H$
as well as a quotient left $H$-module coalgebra through the identification $R \simeq H/HA^+$.
Since $H$ is also cosemisimple, $R$ decomposes as a direct sum $R =
\oplus_{i} V_i$, where $V_i$ are irreducible left coideals of $H$.
It is clear that any left coideal $V$ of $H$ such that $V \subseteq R$
is an $A$-subcomodule of $R$.

The following lemma gives insight into the relationship between the $H$-comodule
structure and the coalgebra structure on $R$.

\begin{proposition}\label{R-coideal}
Let $V \subseteq R$ be a left coideal of $H$. Then the following hold:

\smallbreak
(i) $V$ is a left coideal of $R$;

 (ii) $\dim \End^R (V) \leq \dim V$,  and the equality
holds if and only if $V$ is a subcoalgebra of $R$. If this is the
case,  and if $V$ is an irreducible left coideal of $H$, then $V$
has multiplicity $1$ as a left $H$-subcomodule of $R$.
\end{proposition}

\begin{proof} (i) The map $p = \id \otimes \epsilon: H \to R$
is a coalgebra surjection, and $p\vert_R = \id_R$.  Since
$\Delta(V) \subseteq H \otimes V$, then $\Delta_R(V) = (p \otimes
p) \Delta(V) \subseteq R \otimes V$, showing that $V$ is a left
coideal of $R$.

 (ii) Let $V = \oplus_i m_i V_i$, where $V_i$ are
simple,  pairwise non-isomorphic, left coideals of $R$. Thus $m_i
\leq \dim V_i$, and $\dim \End^R (V) = \sum_i m_i^2 \leq \dim V$.
Moreover, $\dim \End^R (V) = \dim V$ if and only if $m_i = \dim
V_i$, for all $i$, if and only if $V$ is a subcoalgebra of $R$.

Suppose now that $V$ and $U$ are irreducible coideals  of $H$ such
that $V$ is a subcoalgebra of $R$ and $V \simeq U \subseteq R$.
Then $U \simeq V$ as left coideals of $R$.  Since $V$ is a
subcoalgebra of $R$, this implies that $U = V$. The proof of the
proposition is now complete.
\end{proof}

\begin{corollary}\label{formula} Suppose that $V \subseteq R$
is a left coideal of $H$ and let  $\chi$ be the character of $V$.
Then we have $$\sum_{\lambda \in \widehat{A^*}} \deg \lambda \
m(\lambda, \chi_V^*\chi_V) \leq \dim V,$$ and the equality holds if and
only if $V$ is a subcoalgebra of $R$.

In particular, $|G(A) \cap G[\chi^*]|  \leq \dim V$. \end{corollary}

Recall that $G[\chi^*] = \{ g \in G(H): \chi g = \chi \}$.

\begin{proof} Combine Proposition \ref{R-coideal} (ii) with Remark \ref{dimension}.
\end{proof}

\begin{remark} Suppose that $A = kG$ is a group algebra and
$V \subseteq R$ is an irreducible left  coideal
of $H$ such that $\dim V = |G|$ and $Vg \simeq V$, for all $g \in
G$, as in \cite[1.3]{pqq2}.

Then we have  $\dim \End^RV = |G|$ and, by Proposition \ref{comm-pair}, there is an isomorphism of algebras $(\End^RV)^{\op} \simeq k_{\alpha}G$ for some $\alpha \in Z^2(G, k^{\times})$.

On the other hand, by Proposition \ref{R-coideal}, $V$ is a  subcoalgebra of
$R$. So that $\End^RV = \End^VV \simeq (V^*)^{\op}$ as algebras.
Whence $V$ is isomorphic to a dual twisted group algebra as a
coalgebra. We thus recover the statement in \cite[Proposition
1.3.1]{pqq2}. \end{remark}

\section{Hopf subalgebras} In this section
we  discuss some results on the existence of proper Hopf
subalgebras in the biproduct $H \simeq R \# A$.

\begin{lemma}\label{util}  Suppose  that $\widetilde R \subseteq R$ is a subspace such that:

(1) $\widetilde R$ is a subalgebra and a subcoalgebra of $R$;

(2) there exists a Hopf subalgebra  $B \subseteq A$ such that
$\rho (\widetilde R) \subseteq B \otimes \widetilde R$ and  $B .
\widetilde R \subseteq \widetilde R$.

Then $\widetilde R$ is a braided Hopf  algebra over $B$ and the
biproduct $\widetilde R \# B$ is a Hopf subalgebra of $H$.
\end{lemma}

\begin{proof} It is not hard to see that  the conditions \eqref{a}, \eqref{b} and \eqref{c}
are verified, so that $\widetilde R$ is a braided Hopf subalgebra
of $R$. Using \eqref{mca}, we see that indeed $\widetilde R \# B$ is
a Hopf subalgebra of $H$.
\end{proof}

\begin{remark}\label{rmk-util} (i) Suppose that  $V \subseteq R$ is a subcoalgebra satisfying condition (2) in Lemma \ref{util}; that is,  suppose that there exists a Hopf subalgebra $B \subseteq A$ such that
$\rho (V) \subseteq B \otimes V$ and  $B .
V \subseteq V$.

Then the subalgebra $k[V]$ generated by $V$ in $R$   is both a
subalgebra and a subcoalgebra, by \eqref{c}; moreover $\rho (k[V])
\subseteq B \otimes k[V]$  and $B .  k[V] \subseteq k[V]$ because
the multiplication of $R$ is a comodule and module map. Therefore, Lemma \ref{util} implies that $k[V]$ is a braided Hopf algebra over $B$ and $k[V]
\# B$ is a Hopf subalgebra of $H$.

(ii) Suppose that there exists a Hopf subalgebra $B \subseteq A$
such that $\rho (R) \subseteq B \otimes R$. Then Lemma \ref{util}
applies with $\widetilde R = R$, and we find that $R \# B$ is a Hopf subalgebra of $H$.
\end{remark}

\begin{lemma}\label{bip-2-dim}
Suppose that $V \subseteq R$ is an irreducible left coideal of $H$. Then
we have

(i)  $V$ is not irreducible as an $A$-subcomodule of $R$.

(ii)  Assume in addition that $\dim V = 2$ and
$\sum_{g \in G(A)} g.V$ generates $R$ as an algebra.
Then $\rho (R) \subseteq kG(A) \otimes R$ and   $R \# kG(A)$ is a
Hopf subalgebra of $H$. \end{lemma}

\begin{proof} (i) Let $\pi := \epsilon \otimes \id: H \to A$
be the  canonical Hopf algebra projection. We have $\rho (V)
\subseteq A \otimes V$.  Suppose on the contrary that $V$ is an
irreducible $A$-subcomodule of $R$.

Then  $\rho (V) \subseteq C_0 \otimes V$, where $C_0 \subseteq A$
is a simple subcoalgebra of dimension $(\dim V)^2$. Let $C
\subseteq H$ be the simple subcoalgebra containing $V$. We have
$\rho (V) = (\pi \otimes \id) \Delta (V) \subseteq \pi (C) \otimes
V$, where $\pi : H \to A$ is the projection. Thus, $C_0 \subseteq
\pi (C)$, and since $\dim C = \dim C_0$, we find that $\pi
\vert_C$ is injective. This is absurd since $\pi\vert_V =
\epsilon\vert_V$ because $V \subseteq R$.

 (ii)  By part (i), $V$ is not irreducible as an
$A$-comodule and therefore $\rho (V) \subseteq kG(A) \otimes V$.
By \eqref{yetter-drinfeld}, $\rho (g.V) \subseteq kG(A) \otimes
g.V$, for all $g \in G(A)$. This implies, in virtue of the
assumption, that $\rho (R) \subseteq kG(A) \otimes R$. Thus, by
Remark \ref{rmk-util} (ii), $R$ is a braided Hopf algebra over
$kG(A)$ and the biproduct $R \# kG(A)$ is a Hopf subalgebra of
$H$.   \end{proof}

 Assume that $\dim A = p^3$, $p$  prime, and $A$ is not
cocommutative. Then the index of $G(A)$ in $A$ is $p$ and $kG(A)$
in normal in $A$. Moreover, the irreducible $A$-comodules have dimension $1$ or $p$. See \cite{ma-pp}.
We thus obtain the following corollary.

\begin{corollary}\label{bip-8} Let $H = R \# A$ be a biproduct, where $A$ is a non-cocommutative semisimple Hopf algebra of dimension $p^3$.
Suppose that  $V \subseteq R$ is an irreducible $p$-dimensional left coideal of
$H$ such that $\sum_{g \in G(A)} g.V$ generates $R$ as an
algebra.

Then $H$ contains a Hopf subalgebra of index $p$.
\end{corollary}

\begin{proof}
By Lemma \ref{bip-2-dim} (i) $V$ is not irreducible as an $A$-comodule;
therefore $\rho (V) \subseteq kG(A) \otimes V$. Then $\rho (R)
\subseteq kG(A) \otimes R$, since by assumption the sum $\sum_{g
\in G(A)} g.V$ generates $R$. Remark \ref{rmk-util} (ii) now
implies that $R \# k G(A)$ is a  Hopf subalgebra of $H$.
\end{proof}

\section{Biproducts over finite groups}
Assume that $A = kG$ is the group algebra of a finite group $G$.
Thus, $R$ is a $G$-graded algebra
$$R = \bigoplus_{g \in G}R_g, \quad R_g = \{ r \in R: \rho(a) = g \otimes a
\},$$ such that $\Delta_R(R_g) \subseteq \bigoplus_{st = g}R_s
\otimes R_t$. The action of $G$ is by algebra and coalgebra
automorphisms, and we have
\begin{equation}\label{accion-coaccion}h . R_g =
R_{hgh^{-1}},\end{equation}
for all $g, h \in G$.
The braiding $\tau_{R,R}: R \otimes R \to R \otimes R$ is given by
\begin{equation}\label{trenza} \tau_{R,R} (a \otimes b) = (g . b)
\otimes a, \quad \forall a \in R_g, \ b \in R.\end{equation}

Let $\Supp R \subseteq G$ denote the set of elements $g
\in G$ such that $R_g \neq 0$. Let also $G_R$ denote the subgroup of $G$
generated by $\Supp R$.

\begin{lemma}\label{soporte} $G_R$ is a normal
subgroup of $G$. Moreover, $R$ is a Yetter-Drinfeld Hopf algebra  over $G_R$
with respect to the coaction $\rho$ and the restricted action of $G_R$ on $R$,
and the biproduct $R \# kG_R$ is a normal Hopf subalgebra in $R \# kG$.
\end{lemma}

\begin{proof} It follows from \eqref{accion-coaccion} that $G_R$ is normal
in $G$. By Remark \ref{rmk-util} (ii), $R \# kG_R$ is a Hopf subalgebra of $R \# kG$. It is normal thanks to \eqref{mca}. See
\cite{mont-indec}. \end{proof}

The following lemma will help to find, in certain cases, normal Hopf subalgebras of Hopf algebras obtained as biproducts.

\begin{lemma}\label{nucleo-accion} (i) Suppose that $G$ contains a normal subgroup $N$ such that $N$ acts trivially on $R$. Then the group algebra $kN$ is a normal Hopf subalgebra in $H$.

(ii) Assume that $R$ is cocommutative and let $n = \dim R$. If $|G|$ does not divide $(n-1)!$ then there exists a subgroup $1 \neq N$ of $G$ such that the group algebra $kN$ is normal in $H$. \end{lemma}

\begin{proof} (i) Since $N$ acts trivially on $R$, then $ha = ah$ in $H$, for all $h \in N$, $a \in R$. On the other hand, the Yetter-Drinfeld condition \eqref{accion-coaccion} implies that $N \subseteq Z_G(G_R)$.

Note the following consequence of \eqref{mca}:
\begin{equation} \ad_{a\# g} (b \# h) =
\left( a^1 \# a^2_{-2} \right) \left( g . b \# g h g^{-1}\mathcal S (a^2_{-1}) \right) \left( \mathcal S_R(a^2_0) \# 1 \right).
\end{equation}
In particular, for all $h \in N$, we have
\begin{align*}\ad_{a\# g} (h) & = a^1 \left( a^2_{-2} (g h g^{-1}) S (a^2_{-1}) \right) S_R(a^2_0) \\ & = a^1 (g h g^{-1}) S_R(a^2) = \epsilon (a) ghg^{-1} \in kN; \end{align*}
the second equality because $g h g^{-1} \in N \subseteq Z_G(G_R)$.
This proves (i).

 (ii) Since $G$ acts on $R$ by coalgebra automorphisms fixing $1 \in G(R)$, then $G$ permutes the set $X = G(R) \backslash \{ 1 \}$. This corresponds to a group homomorphism $f: G \to \mathcal S(X)$, where $\mathcal S(X)$ is the group of all permutations of $X$.
Let $N$ be  the kernel of $f$. By assumption $N \neq 1$ and $N$ acts trivially on $R$. Then the claim follows from part (i). \end{proof}

\begin{proposition}\label{dimR-3-4} Suppose that $\dim R = 3$ or $4$.
Then  either $H$ or $H^*$ contains a proper normal cocommutative
Hopf subalgebra.  \end{proposition}

\begin{proof} In this case, $R$ is necessarily cocommutative.
By Lemma \ref{nucleo-accion}, we may assume that  $|G| = 2$ if
$\dim R = 3$ and $|G| = 6, 3$ or $2$ if $\dim R = 4$. In view of
\cite{ma-6-8} and \cite{fukuda}, it remains to consider the case
where $|G| = 6$, $\dim R = 4$. In this case, $X = G(R) \backslash
\{ 1 \}$ has three elements and we may assume that $G \simeq
\mathcal S(X)$ is not abelian; so that $\dim H = 24$ and $|G(H)|$
is divisible by $6$. We may then also assume that $H$ is not
trivial and $|G(H)| \neq 12$. Therefore $G = G(H)$ and $H$ is of
type $(1, 6; 3, 2)$ as a coalgebra. The proposition will follow
from the following claim.

{\bf Claim.} The group $G$ is abelian.

{\it Proof of the Claim.}
We have that  $[H: G] = 4$ and $\dim R(H^*) = 8$. Consider the inclusion $kG \subseteq R(H^*)$. If  $G$ is not abelian, then $R(H^*) \simeq M_2(k) \times k^{(4)}$ as an algebra, and a complete set of orthogonal primitive idempotents of $kG$ is of the form $e_0, e_1, f_0, f_1$, where $\dim He_i = 4$ and $\dim Hf_j = 8$. Moreover, $f_0, f_1$ are not central.

The Kac-Zhu Theorem implies that $R(H^*)$ has a primitive idempotent $e$ with $\dim He = 2$, and by \cite{zhu2} $D(H)$ has an irreducible character of degree $2$. By Proposition \ref{doble} $|G(D(H)^*)| \neq 1$. In view of \cite{R}, this implies that the subgroup of the group $G(D(H)) = G(H^*) \times G(H)$ consisting of elements which are central in $D(H)$ is not trivial. We may assume that the elements of this group are of the form $\eta \otimes g$, where $g \neq 1$ if $\eta \neq \epsilon$.
Thus $Z(G) \neq 1$ and the claim follows. \end{proof}

\section{Cocommutative braided Hopf algebras}
Suppose  that $R$ is a cocommutative coalgebra, that is, $a^1
\otimes a^2 = a^2 \otimes a^1$, for all $a \in R$; so that  all
irreducible $R$-comodules are one-dimensional.

\begin{lemma}\label{estabilizadores} Let $V \subseteq H = R \# kG$ be an irreducible left coideal of dimension $\dim V > 1$ and let $\chi = \chi_{V^*}$. Then $G[\chi] \cap G \neq 1$. \end{lemma}

\begin{proof} As coalgebras, $R \simeq H / H (kG)^+$. Let $C \subseteq H$ be the simple subcoalgebra containing $V$. Then the (right) stabilizer of $C$ in $G$ is $G[\chi] \cap G$.
Consider the corestriction $\overline V$ of $V$ to $R$ as in Chapter \ref{tres}.

By Corollary \ref{cor-cociente}, the image of $C$ under the canonical  projection $H \to H / H (kG)^+$ is isomorphic to $C / C (k(G[\chi] \cap G))^+$. Therefore, applying  Proposition \ref{comm-pair}
with $A = k(G[\chi] \cap G)$,
we find an isomorphism of algebras $(\End^R \overline V)^{\op} \simeq k_{\alpha}(G[\chi] \cap G)$, where $\alpha \in Z^2(G[\chi] \cap G, k^{\times})$. Since $\dim V > 1$, $\overline V$ cannot be irreducible. Therefore $|G[\chi] \cap G| \neq 1$ and the lemma follows. \end{proof}

\begin{proposition}
Let  $|G| = q^r$, where $q$ is  a prime number and $r > 0$.
 Then $q$ divides the dimension of
$V$, for  all irreducible left coideals $V$ of $H$ such that $\dim
V > 1$.

Assume in addition that $R \cap G(H) = 1$. Then $\dim R = 1 \mod q$. \end{proposition}

\begin{proof} Let $V$ be an irreducible left coideal of $H$ of dimension $\dim V > 1$. It follows from Lemma \ref{estabilizadores} that $G[\chi_V] \cap G \neq 1$ and since $G[\chi_V] \cap G \subseteq G$, we find that $q$ divides $|G[\chi_V] \cap G|$.

Since $|G[V]|$ divides the dimension of $V$, it follows that $q / \dim V$.

Suppose now that $R \cap G(H) = 1$.
Observe that $R \cap G(H)$ coincides with the set $G(R)^{\co \rho}$ of coinvariant group-like elements in $R$.

Since $R$ decomposes as a direct sum of irreducible left coideals of $H$: $R = k1 \oplus V_1 \oplus \dots \oplus V_m$, and by the assumption $\dim V_j > 1$, for all $j = 1, \dots, m$, the last claim follows. \end{proof}

\section{Cocommutative braided Hopf algebras over $\mathbb Z_p$} We begin this section by reviewing some of the results in \cite{So}, concerning the classification of cocommutative cosemisimple braided Hopf algebras over groups of prime order, that will be used later. We then show some applications.

\smallbreak The braided Hopf algebra $R$ is called \emph{trivial}
if the braiding $\tau_{R, R}: R \otimes R \to R \otimes R$ is the
canonical flip of vector spaces, {\it i.e.}, $\tau_{R, R} (a
\otimes b) = b \otimes a$, for all $a, b \in R$. See \cite[Definition
1.1]{So}.

By \cite{Scha} $R$ is trivial if and only if $R$ is a (usual) Hopf algebra;
that is, if and only if $\Delta_R(ab) = a^1 b^1 \otimes a^2 b^2$, for all $a, b \in R$.

  Let $p$ be a prime number and let $\mathbb Z_p$ denote
the cyclic group of order $p$. Let $R$ be a braided semisimple Hopf algebra over $\mathbb Z_p$. Suppose that $R$ is a cocommutative coalgebra.
It is shown in \cite[Proposition 7.2]{So} that if $R$ is nontrivial,
then $p$ divides the dimension of $R$. More precisely, we have the following proposition.

\begin{proposition}\label{som}  Let $H = R \# k\mathbb Z_p$ such that $H$ is not cocommutative. Suppose that $R$ is cocommutative. Then we have:

(i) Assume that $R$ is  trivial. Then $H$ fits into an
\emph{abelian} central extension
\begin{equation}\label{ex-ses0} 0 \to k\mathbb Z_p \to H \to R \to 1. \end{equation}

(ii) Assume that $R$ is not trivial. Then there are exact sequences of Hopf algebras
\begin{align}\label{ex-ses1}
& 1 \to k^{\mathbb Z_p} \to H \to k F \otimes k \mathbb Z_p \to 1, \\
\label{ex-ses2} & 1 \to k^{\mathbb Z_p}  \otimes k \mathbb Z_p \to H \to k F \to 1, \end{align}
for a certain group $F$. \end{proposition}

\begin{proof} (i) Since $R$ is  cocommutative and trivial, it is isomorphic as
a Hopf algebra to a group algebra $R \simeq k\Gamma$ and the
action and coaction of $\mathbb Z_p$ on $R$ correspond,
respectively, to actions by group automorphisms $\mathbb Z_p \to
\Aut \Gamma$. By \cite[Proposition 1.11]{So}, either the action or
the coaction of $\mathbb Z_p$ on $R$ must be trivial. If the
coaction is trivial, it follows from \eqref{mca} that $H$ is
cocommutative, against the assumption. Therefore the action is
trivial, and again by \eqref{mca}, we find that the coalgebra
surjection $\pi = \id \otimes \epsilon: H \to R$ is indeed a Hopf
algebra surjection; moreover, we have $k\mathbb Z_p = H^{\co
\pi}$. Therefore, $\pi$ gives rise to an abelian central extension
\eqref{ex-ses0}.

(ii) By \cite[Theorems 7.7 and 7.8]{So}, $R$ is isomorphic to a crossed product $k^{\mathbb Z_p} \#_{\sigma} kF$, where $F$ is a finite group, with 'diagonal' action and coaction of $k \mathbb Z_p$. These braided Hopf algebras fit into the construction described in Section 3 of \cite{braid-ext}. Therefore, by \cite[Proposition 3.7]{braid-ext}, there are exact sequences of Hopf algebras
\eqref{ex-ses1} and \eqref{ex-ses2}. \end{proof}

As an application, we shall give in the following theorem a classification result for certain semisimple Hopf algebras. A stronger result (without the assumption on $|G(H^*)|$) is proved in \cite[Corollary IX.9]{IK} in the context of Kac algebras.

\begin{theorem}\label{1-2-2-n} Suppose that $H$ is of type $(1, 2; 2, n)$ as a coalgebra.
Assume in addition that $|G(H^*)|$ is even. Then $H$ is
commutative.
\end{theorem}

\begin{proof} We have $\dim H = 2 (2n+1)$. In particular, since $4$ does not divide $\dim H$, it follows from Theorem \ref{thm-nr} that every irreducible character of degree $2$ is stable under left multiplication by $G(H)$. Therefore, by Remark \ref{particular}, $H/H(kG(H))^+$ is a cocommutative coalgebra.
By assumption, $G(H^*)$ contains a subgroup $\Gamma \simeq \mathbb Z_2$, and the projection $q: H \to k^{\Gamma}$ restricts injectively to $kG(H) \simeq k\mathbb Z_2$, because $\dim H^{\co q} = 2n+1$ is odd. Hence $H : = R \# k\mathbb Z_2$, where $R \simeq H/H(kG(H))^+$ is a cocommutative Yetter-Drinfeld Hopf algebra of odd dimension. By Proposition \ref{som}, since $H$ is not cocommutative,  $H$ fits into an exact sequence  $1 \to k\mathbb Z_2 \to H \to R \to 1$, where $R$ is a cocommutative Hopf algebra. Hence $R = kF$, where $F$ is a group of order $2n+1$.

\begin{claim}\label{f-abel} The group $F$ is abelian. \end{claim}

\begin{proof} By dualizing, we get an exact sequence  $1 \to k^F \to H^* \to k\mathbb Z_2 \to 1$. Therefore, as an algebra, $H^*$ is isomorphic to a crossed product $H^* \simeq k^F \#_{\tau}k\mathbb Z_2$, corresponding to an action
$\rightharpoonup : k\mathbb Z_2 \times k^F \to k^F$ and a
2-cocycle  $\tau: \mathbb Z_2 \times \mathbb Z_2 \to k^F$. By
\cite{pqq}, we may assume that $\tau = 1$.

The action is in this case of the form $g \rightharpoonup \delta_x
= \delta_{g \triangleright x}$, for all $g \in \mathbb Z_2$, $x
\in F$, where $\triangleright : \mathbb Z_2 \times F \to F$ is an
action by group automorphisms (because the extension is cocentral,
see \cite{pqq}). By Clifford theory, the irreducible $H^*$-modules
are exactly of the form $V_{x, U} = \Ind_{k^F \# kG_x}^{H^*} kx
\otimes U$, where $x \in F$ runs over a system of representatives
of the orbits of $\mathbb Z_2$ in $F$, $G_x$ is the stabilizer of
$x$ in $\mathbb Z_2$ and $U$ runs over all irreducible
nonisomorphic $kG_x$-modules. Note that $\dim V_{x, U} = [\mathbb
Z_2: G_x]$. The assumption on the coalgebra structure of $H$
implies that the action $\triangleright : \mathbb Z_2 \times F \to
F$ has exactly one fixed point $x = 1$.

Let $\phi \in \Aut F$, $\phi (x) = g \triangleright x$, where $1 \neq g \in \mathbb Z_2$. Then $g$ is an automorphism of order 2, with exactly one fixed point $x = 1$. This implies that $F$ is abelian and moreover $g \triangleright x = x^{-1}$, for all $x \in F$. \end{proof}

The proof of the theorem can be now concluded as follows: we have a Hopf algebra inclusion $k^F \subseteq H^*$. Since the group $F$ is abelian, $k^F = k\widehat F$ is cocommutative, and thus $2n+1$ divides the order of $G(H^*)$. Since by assumption the order of $G(H^*)$ is even, then $|G(H^*)| = \dim H$ and $H^* = kG(H^*)$ is cocommutative. One could alternatively use the more general statement on abelian exact sequences that we prove in Lemma \ref{sek} below.   \end{proof}

\begin{lemma}\label{sek} Suppose that $H$ fits into an exact sequence $$1 \to k^{\mathbb Z_2} \to H \to kF \to 1,$$ where $F$ is an abelian group of odd order. Assume in addition that $H$ is of type $(1, 2; 2, n)$ as a coalgebra.
Then $H$ is commutative. \end{lemma}

\begin{proof} Note that the exact sequence is necessarily central, by the dual version of Corollary \ref{kob-mas}.
Hence, in the associated matched pair $(\mathbb Z_2, F)$,  the
action $\triangleleft : \mathbb Z_2 \times F \to \mathbb Z_2$ is
trivial and the action $\triangleright : \mathbb Z_2 \times F \to
F$ is by group automorphisms.

As an algebra, $H \simeq k^{\mathbb Z_2} \#_{\sigma} kF$ is a crossed product for the trivial action corresponding to $\triangleleft$ and for some 2-cocycle $\sigma: F \times F \to k^{\mathbb Z_2}$, and as a coalgebra $H \simeq k^{F} \#_{\tau} k\mathbb Z_2$ with respect to the action of $\mathbb Z_2$ on $k^F$ corresponding  to $\triangleright$ and the \emph{trivial} cocycle $\tau$; see \cite{pqq}.

To establish the lemma, it will be enough to show that the cocycle $\sigma$ is a coboundary, that is, it is symmetric.

Let $g$ be the generator of $\mathbb Z_2$. The assumption on the
coalgebra structure of $H$ implies, as in the proof of Theorem
\ref{1-2-2-n}, that $g$ acts on $F$ by $g \triangleright x =
x^{-1}$ for all $x \in F$. Write $$\sigma(x, y) = \sigma_1(x, y)
\delta_1 + \sigma_g(x, y)\delta_g, \qquad x, y \in F,$$ where
$\delta_1, \delta_g \in k^{\mathbb Z_2}$ are the canonical
idempotents. The normalized cocycle condition for $\sigma$ implies
that $\sigma_g : F \times F \to k^{\times}$ is a 2-cocycle and
$\sigma_1 = 1$ \cite{Maext}. Moreover, the compatibility condition
in \cite[(4.8)]{Maext} between $\sigma$ and the trivial cocycle
$\tau$ implies that
$$\sigma_g(x, y) \sigma_g(x^{-1}, y^{-1}) = 1, \qquad \forall x, y \in F.$$
So that $\sigma_g (x^{-1}, y^{-1}) = \sigma_g(x, y)^{-1}$, for all $x, y \in F$. But the 2-cocycle $\beta: F \times F \to k^{\times}$, $\beta(x, y): = \sigma_g (x^{-1}, y^{-1})$ is cohomologous to $\sigma_g$, because of the injectivity of the map
$$\phi: H^2(F, k^{\times}) \to \Hom(\Lambda^2(F), k^{\times}), \qquad \phi ([\alpha]) (x, y): = \alpha(x, y)\alpha(y, x)^{-1}.$$

Hence, in $H^2(F, k^{\times})$, we have $[\sigma_g] = [\beta] = [\sigma_g]^{-1}$, implying that $[\sigma_g]^2 = 1$, and {\it a fortiori} that $[\sigma_g] = 1$ since the order of $F$ is odd. This finishes the proof of the lemma. \end{proof}

The  following  more general form of Lemma \ref{sek} was suggested by a comment of the referee.

\begin{corollary}\label{sek-rf} Suppose that $H$ is of type $(1, 2; 2, n)$ as a coalgebra. Assume in addition that $H$  fits into an exact sequence $$1 \to k^{\mathbb Z_2} \to H \to K \to 1,$$ where $K$ is a  Hopf algebra of odd dimension. Then $K$  is cocommutative and $H$ is commutative. \end{corollary}

\begin{proof} Note first  that $K$ is cocommutative: this follows
from Remark \ref{particular}, since  $K \simeq H/H(kG(H))^+$, and
all simple subcoalgebras  of dimension 4 are necessarily stable
under left and right multiplication by $G(H)$. Thus $K \simeq kF$,
for some group $F$ of odd order.  In view of Lemma \ref{sek} the
corollary  will follow if we show that $K$ is also commutative. We
do this as follows: consider the matched pair $(\mathbb Z_2, F)$
associated with $H$. As in the proof of Claim \ref{f-abel} we see
that the assumption on the coalgebra structure of $H$ implies that
the action $\triangleright: \mathbb Z_2 \times F \to F$ is by
group automorphisms and cannot have fixed points $x \neq 1$. Then
the group is abelian, and $K$ is commutative.
\end{proof}

\section{Transitive actions of central subgroups}
We  shall assume in this section that $R$ contains a
cocommutative subcoalgebra $V$ such that $\rho (V) \subseteq kG
\otimes V$. As before, we shall use the notation $G_V$ to indicate
the subgroup of $G$ generated by $\Supp V$. We have $\rho (V)
\subseteq kG_V \otimes V$ and there is  a basis $v_1, \dots, v_n$
of $V$ such that $\rho (v_i) = g_i \otimes v_i$, where $g_i \in
\Supp V$, $\forall i = 1, \dots, n$.

\begin{lemma}\label{tr-act} Suppose that there exists a central subgroup $S$ of $G$ with the property that $S$ permutes transitively the set $G(V)$. Then $G_V$ is an abelian subgroup of $G$ of order $\leq \dim V$. \end{lemma}

\begin{proof} Let $G(V) = \{ a_1, \dots, a_n \}$.
Since $G(V)$ is a basis of $V$, we may write $\rho (a_i) = \sum_j t_{ij} \otimes a_j$, for some $t_{ij} \in kG$.
As a consequence of the $kG$-colinearity of $\Delta_R$ we get the
following relations:
\begin{equation}\label{idemp}t_{ij} t_{il} = \delta_{jl}t_{ij}, \quad \forall 1 \leq i, j, l \leq n.\end{equation}
On the other hand, the coassociativity of $\rho$ implies that, for
all $p = 1, \dots, n$, the subspace $R_p$ spanned by the set $\{
t_{pj}: \ 1 \leq j \leq n \}$ is a right coideal of $kG$; by
definition, we have $kG_V = \sum_{p}R_p$.

The assumption on $S$ implies that $R_p = R_j$, for
all $1 \leq p, j \leq n$. Fix $p$, and let  $e_j : = t_{pj}$, $1
\leq j \leq n$. By \eqref{idemp} we get that $kG_V $, which is
spanned by $\{ e_j \}$, is a commutative subalgebra of $kG$ of
dimension $\leq n$. This proves the lemma. \end{proof}

\chapter[Cocycle Deformations]{Cocycle Deformations of Some Hopf Algebras}\label{twist}

In this chapter, we describe some known examples of semisimple Hopf algebras as cocycle twists of group algebras. We also show that some others cannot be obtained  in this fashion.

\section{Lifting from abelian groups} Let $G$ be a finite group and let $\phi \in kG \otimes kG$ be a normalized invertible $2$-cocycle. The following lemma is a consequence of \cite{eg-isocat}.

\begin{lemma} The Hopf algebra $(kG)_{\phi}$ is cocommutative if and only if there exist a normal abelian subgroup $A \subseteq G$ and an $\ad G$-invariant cohomology class $[\omega] \in H^2 (\widehat A, k^{\times})$ such that
\begin{equation}\label{lifted} \phi = \sum_{x, y \in \widehat A} \omega (x, y) \delta_x \otimes \delta_y \in A \otimes A,\end{equation} where, for $x \in \widehat A$, the element $\delta_x \in A$ is the primitive idempotent defined by $\delta_x : = \dfrac{1}{|A|}\sum_{a \in A} x(a) a$. \end{lemma}

\begin{proof} Note that $(kG)_{\phi}$
is cocommutative if and only if $\phi \phi_{21}^{-1}$ commutes
with  $\Delta (G)$.
The 'if' direction follows from this observation and the invariance of
$[\omega]$. The 'only if' part follows from \cite[Proof of Theorem
1.3]{eg-isocat}.
\end{proof}

\begin{remark}\label{rmk-lifted} Suppose that the $2$-cocycle $\phi$ has the form \eqref{lifted}.
Then we say that $\phi$ is \emph{lifted} from the (not necessarily
normal) abelian subgroup $A$; see \cite{nik}.

Suppose that $\phi$ is lifted from the normal abelian
subgroup $A$.
Recall the isomorphism  $$\theta: H^2(\widehat A, k^{\times}) \to
\Hom (\Lambda^2 (\widehat A), k^{\times}), \quad \theta ([\sigma]) (x,
y) = \sigma (x, y) \sigma (y, x)^{-1}.$$ This isomorphism commutes
with the $\ad G$-action (i.e., the action coming from the adjoint
action of $G$ on $A$).

 Then $(kG)_{\phi}$ is cocommutative if and only if
$[\omega]$ is $\ad G$-invariant in $H^2 (\widehat A, k^{\times})$.
\end{remark}

In the rest of this chapter $p$ and $q$ will be distinct prime numbers.

\section{Examples in dimension $pq^2$; $p = 1 \mod q$}\label{caso1} Suppose that $p = 1 \mod q$. Let $\mathcal A_i$, $0 \leq i \leq q-1$ be the nontrivial semisimple Hopf algebras constructed in \cite{G}.
The Hopf algebras $\mathcal A_i$ are self-dual and we have $G(\mathcal A_0) \simeq \mathbb Z_q \times \mathbb Z_q$, while $G(\mathcal A_i) \simeq \mathbb Z_{q^2}$ for all $i = 1, \dots , q-1$; see \cite{examples}.

Let $F = \langle a: \ a^p = 1 \rangle$ be the cyclic group of
order $p$ and let  $\Gamma = \langle s, t: \ s^q = t^q =
sts^{-1}t^{-1} = 1 \rangle \simeq \mathbb Z_q \times \mathbb Z_q$.

  Let $1 < m, l \leq p-1$ be units modulo $p$ such that
$m^q = l^q = 1 \mod p$.  Consider the semidirect product $G = F
\rtimes \Gamma$ corresponding to the action by group automorphisms
of $\Gamma$ on $F$ defined on generators by $s.a = a^m$ and $t . a
= a^l$.

Let $\omega \in Z^2(\widehat \Gamma, k^{\times})$ be a nontrivial $2$-cocycle; in particular, the cohomology class $[\omega]$ generates $H^2 (\widehat \Gamma, k^{\times}) = \Hom (\Lambda^2\widehat \Gamma, k^{\times}) \simeq \mathbb Z_q$.
Consider the invertible normalized $2$-cocycle $\phi = \sum_{x, y \in \Gamma} \omega (x, y) \delta_x \otimes \delta_y \in k\Gamma \otimes k\Gamma$.

The following proposition generalizes \cite[Example 2.9]{nik}. It also follows from the results in \cite[Theorem 4.8]{masuoka-cont} in the case where $p = 3$ and $q = 2$.

\begin{proposition}\label{twist-pqq} (i) $\mathcal A_0 \simeq (kG)_{\phi}$ as Hopf algebras;

 (ii) Let $1 \leq i \leq q-1$. Then the Hopf algebras $\mathcal A_i$ cannot be obtained  from a group algebra by twisting the comultiplication by means of a $2$-cocycle.
\end{proposition}

\begin{proof} (i) We claim that the Hopf algebra $(kG)_{\phi}$ is nontrivial.
Note that the twisted group algebra $k_{\omega}\Gamma$ is simple; so that $\phi$ is a minimal twist for $k\Gamma$ in the sense of \cite{eg-triangular}.
Therefore we have $$k\Gamma = \{ (f \otimes \id) (b); \ f \in k^{\Gamma} \},$$ where $b = \phi_{21}^{-1}\phi  \in k\Gamma \otimes k\Gamma$.

Suppose that $(kG)_{\phi}$ is cocommutative. This is equivalent to the condition  $b \Delta (g) = \Delta (g)  b$, for all $g \in G$, or
$$gb^1g^{-1} \otimes gb^2g^{-1} = b^1 \otimes b^2, \quad \forall g \in G.$$
Regard $\widehat \Gamma$ as a subgroup of $\widehat G$ by transposing the natural projection $G \to \Gamma$.
For all $f \in \widehat \Gamma$, after applying $f \otimes \id$ to both sides of the above equation, we have $g (f(b^1)b^2) g^{-1} = f(b^1)b^2$.
Since $\widehat \Gamma$ spans $k^{\Gamma}$, this implies that $\Gamma$ is central in $G$, which is a contradiction. This establishes the claim.

Since $\Gamma$ is abelian and $\phi \in k\Gamma \otimes k\Gamma$,
it follows  $G((kG)_{\phi}) = \Gamma$. We know that $(kG)_{\phi}$
and $(kG)_{\phi}^*$ are of Frobenius type \cite{eg-rtq}.
Therefore, by the classification results for semisimple Hopf
algebras of dimension $pq^2$   \cite[Theorem 5.4.2]{clspqq},
$(kG)_{\phi} \simeq \mathcal A_0$ as claimed.

(ii) In this case we have $G(\mathcal A_i) \simeq \mathbb Z_4$. Suppose on the contrary that $\mathcal A_i = (kN)_J$, for some group $N$ and $2$-cocycle $J \in kN \otimes kN$.
By \cite{eg-triangular} there exists a subgroup $H$ of $N$ such that $J \in kH \otimes kH$ is a minimal twist for $kH$. In particular, $|H|$ is a square and therefore $|H| = q^2$. This forces $H \simeq \mathbb Z_q \times \mathbb Z_q$ in view of the minimality of $J$. But then $kH \subseteq (kN)_J$ is a Hopf subalgebra of dimension $q^2$, necessarily isomorphic to $kG(\mathcal A_i)$. This is a contradiction.
Thus part (ii) follows. \end{proof}

\section{Examples in dimension $pq^2$; $q = 1 \mod p$} Suppose that $q = 1 \mod p$. There is a family of Hopf algebras
of dimension $pq^2$, constructed in \cite{pqq} as a generalization
of the examples in \cite{masuoka-further}.  This family is
parametrized by the Hopf algebras $\mathcal B_{\lambda}$, $0 \leq
\lambda < p-1$, where $\lambda \lambda' \neq 1 \mod p$.

We have $G(\mathcal B_{\lambda}) \simeq \mathbb Z_q \times \mathbb Z_q$ for all $\lambda$, $|G(\mathcal B_{\lambda}^*)| = p$ for $\lambda > 0$, and $G(\mathcal B_0^*) \simeq \mathbb Z_{pq}$.

Let $F$ and $\Gamma$ be as before. Let $G = G_{\lambda}$ be the semidirect product $\Gamma \rtimes F$ with respect to the action of $F$ on $\Gamma$ given by
$a. s = s^n$, $a. t = t^{n^{\lambda}}$, where $0 < n \leq q-1$ is such that $n^p = 1 \mod q$, and $0 \leq \lambda < p-1$.

Let $[\omega ] \in H^2 (\widehat \Gamma, k^{\times})$ be the $2$-cocycle corresponding to the skew-symmetric non-degenerate bicharacter $\Omega: \widehat \Gamma \times \widehat \Gamma \to k^{\times}$, $\Omega (u, v) = \det (u v)$.

We have $\Omega (a. u, a. v) = n^{\lambda + 1} \Omega (u, v)$.
Therefore, the bicharacter $\Omega$, and hence also the
class $[\omega]$ is  \emph{not} $\ad G$-invariant.

Consider the $2$-cocycle $\phi \in kG \otimes kG$ lifted from $\Gamma$ as in Remark \ref{rmk-lifted}, corresponding to $\omega$.

\begin{proposition} (i) $\mathcal B_{\lambda} \simeq (kG_{\lambda})_{\phi}$ as Hopf algebras;

 (ii) The Hopf algebras $\mathcal B_{\lambda}^*$ are not a cocycle deformation of any finite group.  \end{proposition}

In the case where $p = 2$, part (i) of the proposition is contained in \cite[Proposition 3.2]{ma-newdir}.

\begin{proof} Part (i) follows from the above discussion, in view of the classification results for semisimple Hopf algebras of dimension $pq^2$ and the description in \cite[1.4]{pqq}.

Part (ii) follows from an argument similar to the one used to prove part (ii) of Proposition \ref{twist-pqq}, since $\mathcal B_{\lambda}^*$ contains no Hopf subalgebra of square dimension. \end{proof}

\section{Normal Hopf subalgebras in cocycle twists} Let $H$ be a semisimple Hopf algebra and let $\phi \in (H \otimes H)^{\times}$ be a $2$-cocycle.

\begin{lemma}\label{hfsb-tw} Let $B \subseteq H$ be a Hopf subalgebra of $H$.
Then $B$ is a Hopf subalgebra of $H_{\phi}$ if and only if $\phi (B \otimes B) \phi^{-1} = B \otimes B$.
In this case, $B$ is normal in $H$ if and only if it is normal in $H_{\phi}$.
\end{lemma}

Note that the Hopf algebra structure  on $B$ as a Hopf subalgebra of $H_{\phi}$ is not {\it a priori} that of a $2$-cocycle twisting of $B$, since we may not have $\phi \in B \otimes B$.

\begin{proof} The first part of the lemma follows easily. As for the second part, recall that $B$ is normal in $H$ if and only if $HB^+ = B^+H$, which depends only on the multiplication and the counit of $H$. \end{proof}

The following corollary is an easy consequence of the lemma.

\begin{corollary}\label{andrus} Suppose that $B \subseteq H$ is a central Hopf subalgebra. Then $B$ is a central Hopf subalgebra of $H_{\phi}$. \qed \end{corollary}

In particular, let $G$ be a finite group and let $\phi  \in kG
\otimes kG$ be a two cocycle. Then if $Z(G) \neq 1$,
$(kG)_{\phi}$ contains a nontrivial central group-like element.

\bigskip The remaining chapters will be devoted, respectively, to the proof  of our main results on
semisimple Hopf algebras of dimension  24, 30, 36, 40, 42, 48, 54 and 56.

\chapter{Dimension $24$}\label{24}

\section{Possible (co)-algebra structures} Let $H$ be a nontrivial semisimple Hopf algebra of dimension
$24$.

\begin{lemma}\label{conteo24} The group $G(H^*)$ is of order $2$, $3$, $4$, $6$, $8$ or $12$.  As an algebra, $H$ is of one of the
following types:
\begin{itemize} \item $(1, 2; 2, 1; 3, 2)$,
\item $(1, 3; 2, 3; 3, 1)$, \item $(1, 4; 2, 5)$, \item $(1, 4; 2, 1; 4, 1)$, \item $(1, 6;  3, 2)$, \item $(1, 8; 2, 4)$, \item $(1, 8; 4, 1)$, \item $(1, 12; 2, 3)$. \end{itemize} \end{lemma}

We shall prove later that the type $(1, 4; 2, 1; 4, 1)$ is
impossible, that is, there exists no semisimple Hopf algebra with
this (co)algebra structure; see Lemma \ref{b-conm} below. The
analogous result for Kac algebras follows from \cite[Proposition
XIV.37]{IK}, since there is no group $G$ for which $kG$ has this
algebra structure.

\begin{proof} The proof follows from counting arguments using \ref{alg-struct}.
\end{proof}

\begin{remark}\label{rmk-24} (i) If $H$ is of type $(1, 12; 2, 3)$ as an algebra,  then $H$ is not simple by Corollary \ref{kob-mas}.

(ii) Suppose that $H$ is of type $(1, 2; 2, 1; 3, 2)$ as an algebra. Then $H$ has a unique irreducible character $\chi$ of degree 2, which necessarily satisfies \eqref{dim6} with $G(H^*) = \{ \epsilon, g \}$.
Therefore there is a quotient Hopf algebra $H \to k \mathbb S_3$, where $\mathbb S_3$ is the only non-abelian group of order $6$.

(iii) Suppose that $H$ is of type $(1, 4; 2, 5)$ or $(1, 4; 2, 1; 4, 1)$ as a coalgebra. Then, by  Proposition \ref{cociente8}, $H$ contains a non-cocommutative Hopf subalgebra of dimension $8$. \end{remark}

\begin{lemma}\label{24-3} Suppose that $H$ is of type $(1, 3; 2, 3; 3, 1)$ as a coalgebra. Then $G(H^*) \cap Z(H^*) \neq 1$. \end{lemma}

\begin{proof} It follows from Remark \ref{el-nr} (i) that $H$ has a Hopf subalgebra $K$ of dimension $12$.
Since the index of
$K$ in $H$ is $2$, then $G(H^*) \cap Z(H^*) \neq 1$.  \end{proof}

\begin{lemma}\label{24-6} Assume that $H$ is of type $(1, 6;  3, 2)$ as a coalgebra. Then $H$ is not simple. \end{lemma}

\begin{proof} Suppose on the contrary that $H$ is simple.
Then $H$ must be of type $(1, 2; 2, 1; 3, 2)$ or $(1, 6;  3, 2)$ as an algebra. Otherwise, by Lemma \ref{conteo24}, Remark \ref{rmk-24} and Lemma \ref{24-3}, there is a Hopf algebra quotient $H \to \overline H$, where $\dim \overline H = 8$.
Therefore $\dim H^{\co \overline H} = 3$. Decomposing $H^{\co \overline H}$ as a direct sum of irreducible left coideals of $H$ (of dimension $1$ or $3$), we find that $H^{\co \overline H} \subseteq kG(H)$ and thus $H^{\co \overline H} = kT \simeq k^T$ is a Hopf subalgebra, where $T$ is the unique subgroup of $G(H)$ of order $3$; this implies that $H$ is not simple.

Therefore, there is a quotient Hopf algebra $\pi: H \to A$, where $A$
is of dimension $6$. We have $\dim H^{\co \pi} = 4$, implying
that $H^{\co \pi} = k1 \oplus V$, where $V$ is an irreducible left
coideal of dimension $3$. It follows that the restriction
$\pi\vert_{kG(H)}: kG(H) \to A$ is an isomorphism. In particular, $A \simeq kG$ and $H \simeq
R \# kG$ is a biproduct, where $G = G(H)$.  The lemma follows from Proposition
\ref{dimR-3-4}.  \end{proof}

\begin{lemma}\label{hopfsb-8} Suppose that $H$ is simple.
If  $H$ has a Hopf subalgebra of dimension $8$, then $H$ is of
type $(1, 2; 2, 1; 3, 2)$ as an algebra.
\end{lemma}

\begin{proof}
By assumption,  there is a Hopf algebra quotient $H^* \to B$,
where $\dim B = 8$; so that $\dim (H^*)^{\co B} = 3$ and thus $(H^*)^{\co
B} = k1 \oplus V$ as a left coideal of $H^*$, where $V$ is an
irreducible left coideal of dimension $2$.

Suppose that $H$ is not of the claimed type as an algebra.  In view of Lemma \ref{conteo24} and the
previous lemmas, there is a Hopf subalgebra $A \subseteq H^*$, with
$\dim A = 8$. We have $A \cap (H^*)^{\co B} = k1$, unless $V \subseteq
A$. However, the last possibility implies that $(H^*)^{\co B}
\subseteq A$, which contradicts Lemma \ref{restriction}.
Therefore, $H^* = R \# A$ is a biproduct, where $R$ is a
$3$-dimensional Yetter-Drinfeld Hopf algebra over $A$.

By Corollary  \ref{bip-8} and Proposition \ref{dimR-3-4}, $H$ is
not simple in this case. Thus $H$ is of type $(1, 2; 2, 1; 3, 2)$, as claimed.
\end{proof}

\begin{corollary}\label{24-8} Suppose that $H$ is of coalgebra type $(1, 8; 2, 4)$ or $(1, 8; 4, 1)$. Then  $H$ is not simple. \end{corollary}

\begin{proof} In both cases we have $|G(H)| = 8$.
By  Lemma \ref{hopfsb-8} and Remark \ref{rmk-24} (ii) we may
assume that there is a Hopf algebra quotient $q: H \to B$, where
$\dim B = 6$. Since $\dim H^{\co B} = 4$, then $\dim H^{\co B}
\cap kG(H) = 2$ or $4$. The last possibility implies that $H$ is
not simple, and the first one implies that $\dim q(kG(H)) = 4$
which contradicts \cite{NZ}. This proves the corollary.
\end{proof}

\begin{lemma}\label{b-conm} There exists no semisimple Hopf algebra with coalgebra type  $(1, 4; 2, 1; 4, 1)$. \end{lemma}

\begin{proof} Suppose on the contrary that $H$ is of type $(1, 4; 2, 1; 4, 1)$ as a coalgebra.
The set of irreducible characters of
$H^*$ consists of $G(H)$, one irreducible character $\lambda$ of
degree $2$ and one irreducible character $\psi$ of degree $4$.
Then $G(H)$ and $\lambda$ span a standard subalgebra of $R(H^*)$,
which corresponds to the Hopf subalgebra $B \subseteq H$, with
$\dim B = 8$. By \cite{ma-6-8}, $B$ is isomorphic to $H_8$ or to
$k^G$, where $G$ is either the dihedral or the quaternionic group
of order $8$.

It is not hard to see that the fusion rules in $R(H^*)$ are
determined by $\lambda^2 = \sum_{g \in G(H)} g$, $\psi^2 = \sum_{g
\in G(H)} g + 2 \lambda + 2 \psi$; in particular, we have $g \psi
= \psi = \psi g$, for all $g \in G(H)$, and $\lambda \psi = 2 \psi
= \psi \lambda$.

\begin{claim} $B$ is commutative. \end{claim}

\begin{proof} Let $C \subseteq H$ be the
$16$-dimensional simple subcoalgebra corresponding to $\psi$. By Remark \ref{fusion-rules} the
fusion rules for $\widehat{H^*}$ imply that $BC = C = CB$. Let
$\overline C = C/CB^+$; so that $\overline C$ is a cocommutative
coalgebra of dimension $2$.

By Corollary \ref{comm-pair-cor}  there is a bijective
correspondence between simple $\overline C$-comodules and simple
$B_{\alpha}$-modules, for some invertible normalized $2$-cocycle
$\alpha: B \otimes B \to k$.

Suppose that $B$ is not commutative; so that $B \simeq H_8$ is the
unique nontrivial semisimple Hopf algebra of dimension 8. It
follows from \cite[Theorem 4.8 (1)]{masuoka-cont} that  any
Galois object of $B$ is trivial; so that  $B_{\alpha} \simeq B$ as
$B$-comodule algebras. This is impossible since $B$ has more than
$2$ non-isomorphic simple modules. Hence $B$ must be commutative.
\end{proof}

We have $B \simeq k^G$, where $G$ is not abelian of order 8. If
$B$ is normal in $H$, then $H$ fits into an abelian extension $1
\to k^G \to H \to kF \to 1$, where $F$ is the cyclic group of
order 3.

Consider the associated matched pair $\triangleright: G \times F
\to F$, $\triangleleft: G \times F \to G$; see \cite{Maext}.
There exists a subgroup
$G_0$ of $G$, with $|G_0| = 4$, which acts trivially on $F$.
The compatibility conditions between $\triangleright$ and $\triangleleft$ imply that $G_0$ is stable under the action of $F$.
This implies, in view of formulas (4.2) and (4.5) in \cite{Maext}, that
the subspace $A_0 = k^F \otimes kG_0$ is a Hopf subalgebra of
$H^*$ of dimension 12, necessarily normal.

We have $H \simeq k^G \#_{\sigma}kF$ is a crossed product, and the
irreducible $H$-modules are of dimension 1 or 3. We may assume $H$
is not commutative, since otherwise $H^* = k(F \rtimes G)$ would
have no irreducible modules of dimension 4. By Lemma
\ref{conteo24}, $H^*$ is of type $(1, 6; 3, 2)$ as a coalgebra.
Hence $A_0$ is commutative, and $H^*$ fits into an abelian
extension $1 \to A_0 \to H^* \to k\mathbb Z_2 \to 1$. Then the
irreducible $H$-comodules are of dimension 1 or 2. In particular,
$H^*$ cannot have irreducible modules of dimension 4. This shows
that $B$ is not normal in $H$.

\begin{claim} $H$ is of type $(1, 2; 2, 1; 3, 2)$ as an algebra. \end{claim}

\begin{proof} Since $B$ is not normal  in $H$, then $(H^*)^{\co B^*} = k1 \oplus V$, where $V$ is an irreducible left coideal of dimension 2. Suppose $H$ is  not of type $(1, 2; 2, 1; 3, 2)$. Combining Theorem \ref{coinvariantes} with  Lemma \ref{conteo24}, we see that the possible coalgebra types for $H^*$ are  $(1, 4; 2, 5)$ and $(1, 8; 2, 4)$.

By Remark \ref{rmk-24}, $H^*$ contains a Hopf subalgebra $K$
of dimension $8$. By Lemma \ref{restriction}, $K \cap (H^*)^{\co
B^*} = k1$. Hence $K \simeq B^*$ is cocommutative and $H^*$ is a
biproduct $H^* \simeq R \# K$, where $\dim R = 3$.

By Lemma \ref{nucleo-accion}, there exists a normal subgroup $1
\neq N \subseteq G(K)$ such that  $N$ commutes with $R$. Since
$|G(K)| = 8$, $N \subseteq Z(G(K))$; so that $N$ commutes also
with $K$. Therefore $N$ is central in $H^*$, and $H^*$ has a
central group-like element of order 2.  This is a contradiction,
since $H$ has no Hopf subalgebra of dimension 12. Hence $H$ is of
type $(1, 2; 2, 1; 3, 2)$ as an algebra, as claimed. \end{proof}

  The irreducible characters of degrees 1 and 2 belong to
a commutative Hopf subalgebra $A \subseteq H^*$.

\begin{claim} $A$ is normal in $H^*$. \end{claim}

\begin{proof} Consider the projection $q: H^* \to k^{G(H)}$. Then we have $\dim (H^*)^{\co q} = 6$.
Suppose that $A$ is not normal. Then $q$ is not normal, and  $G(H^*) \cap (H^*)^{\co q} = 1$.
Let $1 \neq g \in G(H^*)$. Thus there exists $s \in G(H)$ such that $\langle g, s \rangle \neq 1$.
Therefore, we must have $H^{\co A^*} = k1 \oplus kt \oplus V$,
where $V$ is an irreducible left coideal of dimension 2. In
particular, $H^{\co A^*} = B^{\co A^*} \subseteq B$.

  Since $[H: B] = 3$, it follows from Lemma \ref{index-3}
and Proposition \ref{doble} that  $G(D(H)^*) \neq 1$; hence $H$
has a one dimensional Yetter-Drinfeld module.

Note that $G(H^*) \cap Z(H^*) = 1$ and $G(H) \cap Z(H) = 1$, since
neither $H$ nor $H^*$ contain  Hopf subalgebras of dimension 12.

By Lemma \ref{yd-1}, a one dimensional Yetter-Drinfeld module of
$H$ is of the form $V_{a, g}$, for some $1 \neq a \in G(H)$.
Consider the projection $q: H \to k^{\langle g \rangle}$.
 Since $B$ is commutative, $a^{-1}ha = h$, for all $h \in B^{\co q}$. By
Theorem \ref{A-comm-g},  $B^{\co q}$ is a Hopf subalgebra of $B$,
which is a contradiction.  \end{proof}

The above claim implies that $H^*$ fits into the abelian extension
$$1 \to k^G \to H^* \to kF \to 1,$$ where $G \simeq \mathbb S_3$
and $F \simeq \mathbb Z_2 \times \mathbb Z_2$.

Then $H^*$ is isomorphic to a crossed product $k^G \#_{\sigma} kF$
as an algebra. Using Clifford theory, we know that the isomorphism
classes of irreducible modules of $H^*$ are exactly the classes of
the modules $$V_{x, \rho} : = \Ind_{k^G\#k_{\alpha}F_x}^{H^*} x
\otimes \rho,$$ where $x$ runs over a system of representatives of
the action of $F$ in $G$, $F_x \subseteq F$ is the stabilizer of
$x$, and $\rho$ is an irreducible representation of a twisted
group algebra $k_{\alpha}F_x$.

We have $\dim V_{x, \rho} = \dim \rho [F: F_x]$. The assumption on
the coalgebra structure of $H$ implies that $\dim V_{x, \rho} =
4$, for some $x \in G$. This implies that $[F:F_x] = 4$ (because
we cannot have $[F:F_x] = \dim \rho = 2$). Hence, the orbit of $x$
has 4 elements, and there must exist $1 \neq y \in G$ such that
$F_y = F$. This implies that $|G(H)| = 8$ and gives a
contradiction. This contradiction shows that $H$ cannot have this
coalgebra structure and finishes the proof of the lemma.
\end{proof}

\begin{lemma}\label{24-2} Assume that $H$ is of type $(1, 2; 2, 1; 3, 2)$ as a coalgebra.
Then $H$ is not simple.  \end{lemma}

\begin{proof} Suppose on the contrary that $H$ is simple. By previous lemmas, $H^*$ is of type  $(1, 4; 2, 5)$ or $(1, 2; 2, 1; 3, 2)$ as a coalgebra.

\begin{claim} $H^*$ is of type $(1, 4; 2, 5)$ as a coalgebra.
 \end{claim}

\begin{proof} Suppose not. Then $H^*$ is of type $(1, 2; 2, 1; 3, 2)$.

Hence,  there is a Hopf algebra quotient $q: H \to B$, where $\dim B
= 6$ and $B$ is cocommutative; in particular,  $\dim H^{\co B} =
4$. On the other hand, $H$ contains a unique Hopf subalgebra $A$
of dimension $6$ which is not cocommutative and such that $G(H)
\subseteq A$.

As a left coideal of $H$,  $H^{\co B} = kG(H) \oplus V$, where $V$
is an irreducible left coideal of dimension $2$, or else $H^{\co
B} = k1 \oplus W$, where $W$ is an irreducible left coideal of
dimension $3$.

The first possibility implies  $H^{\co B} \subseteq A$, which
contradicts Lemma \ref{restriction}. The second possibility
implies that $A \cap H^{\co B} = k1$ and therefore the restriction
$q\vert_A: A \to B$ is an isomorphism. This is impossible, since
$A$ is not cocommutative. Hence the claim follows. \end{proof}

By Remark \ref{rmk-24} (iii), there is a Hopf subalgebra  $B \subseteq H^*$ of coalgebra type $(1, 4; 2, 1)$.

Since $[H^*: B] = 3$, it follows from Lemma \ref{index-3} and Proposition \ref{doble}, $G(D(H)^*) \neq 1$. By
Lemma \ref{yd-1}, this is the form $V_{g, \eta}$,
for some $1 \neq g \in G(H)$ and $1 \neq \eta \in G(H^*)$.

Consider the projection $q: H \to k^{\langle \eta \rangle}$.
 Since $A$ is
commutative, $g^{-1}ag = a$, for all $a \in A^{\co q}$. By
Theorem \ref{A-comm-g},  $A^{\co q}$ is a Hopf
subalgebra of $A$. This implies that $A^{\co q} = A \subseteq H^{\co q}$. In particular, $\eta\vert_A = \epsilon$, and this  $\eta$ is the only group-like element of $H^*$ with this property: otherwise, $A = H^{\co G(H^*)}$ and $H$ is not simple.

Necessarily, $(H^*)^{\co A^*} = k1 \oplus ks \oplus U$, where $U$
is an irreducible left coideal of $B$ of dimension $2$. It follows from Corollary  \ref{cor-mk}, that $s \in G(B) \cap Z(B)$. Also, $s\vert_A = \epsilon$; therefore $s = \eta$.

On the other hand, $(H^*)^{\co A^*} \subseteq (H^*)^{\co G(H)}$. Since $B$ cannot be contained in $(H^*)^{\co G(H)}$ by \cite{NZ}, we see that $B^{\co G(H)} = B^{\co A^*} = k1 \oplus ks \oplus U$.

By Remark \ref{yd-*}, $V_{s, g}$ is a Yetter-Drinfeld module of $H^*$. Since
$s \in Z(B)$, applying again Theorem \ref{A-comm-g},  we see that $B^{\co G(H)}$ is a Hopf subalgebra of $B$. This is a contradiction, because of the decomposition of $B^{\co G(H)}$ as a left coideal of $H^*$. This contradiction finishes the proof of the Lemma.  \end{proof}

\section{Upper and lower semisolvability} Our aim in this section is to prove that semisimple Hopf algebras of
dimension $24$ are  \emph{both} upper and lower semisolvable. We first need the following
lemma:

\begin{lemma}\label{ss-24} Suppose that $H$ is not simple.
Then $H$ is upper and lower semisolvable. \end{lemma}

\begin{proof} We shall show that $H$ is lower semisolvable, and the other statement will follow by duality.
Note that every semisimple Hopf algebra of dimension less than $24$
is upper and lower semisolvable; see Table 1.  Therefore if $H$ has a
normal quotient $H \to \overline H$, where $\overline H$ is
commutative or cocommutative, then $H$ is lower semisolvable. In
particular, we may assume that $G(H^*) \cap Z(H^*) = 1$. This,
combined with Lemma \ref{24-3} and Remark \ref{rmk-24} (i), allows
us to suppose that $H$ is neither of type $(1, 3; 2, 3; 3, 1)$ nor $(1, 12; 2, 3)$ as a coalgebra.

By assumption, $H$ fits into an extension $$\begin{CD}1 @>>> A @>>> H
@>{\pi}>> \overline H @>>> 1,\end{CD}$$ and we may assume that
$\overline H$ is not trivial. Thus $\dim \overline H = 8$ or $12$,
and $\dim A = 3$ or $2$, respectively.

In the first case,  we may suppose that $H$ is
of type $(1, 6;  3, 2)$ as a coalgebra, by Lemma \ref{conteo24}. In
particular, $G(H)$ has a unique (normal) subgroup $G$ of order
$3$, which must coincide with the stabilizer of all simple
subcoalgebras of dimension $9$, such that $A = kG$; Remark
\ref{particular} now implies that the quotient $\overline H =
H/H(kG)^+$ is cocommutative, and $H$ is semisolvable.

  Suppose finally that  $\dim \overline H = 12$ and $\dim
A = 2$. We observe that if $|G(H^*)| = 12$, then any subgroup of
order $3$ must be contained in $(H^*)^{\co A^*} \simeq
\overline{H}^*$; thus $3 / |G(\overline{H}^*)|$ and
$\overline{H}^*$ is commutative or cocommutative. Hence, we may
assume that $|G(H^*)| \neq 12$.

If $\overline H$ is not trivial,  then its simple subcoalgebras
are of dimension $1$ or $4$ \cite{fukuda}. Therefore, by Corollary
\ref{cor-cociente}, we may assume that $H$ contains no simple
subcoalgebra of dimension $9$; this implies that $|G(H)| = 4$ or
$8$, by Lemma \ref{conteo24}. Also, $4 = |G(\overline{H}^*)|$ divides
$|G(H^*)|$, hence we must only consider the cases $|G(H^*)| = 4$
or $8$.

It turns out in this case, that both  $H$ and $H^*$ contain Hopf
subalgebras of dimension $8$. If $B \subseteq H$, $K \subseteq H^*$, with $\dim B = \dim K = 8$, then $B\cap H^{\co K^*} = k1$, by Lemma \ref{restriction}. Hence $H = R \# B$, where $B$ is
a Hopf subalgebra of dimension $8$.

If $B$ is not cocommutative, then by  Corollary \ref{bip-8}, $H$
has a normal Hopf subalgebra of dimension $12$, implying that
$G(H^*) \cap Z(H^*) \neq 1$; hence we are done in this case. If $B
= kG(H)$ is cocommutative, it follows from the proof of
Lemma  \ref{nucleo-accion} (ii) that $H$ has a normal Hopf subalgebra
of dimension $4$. Thus $H$ is semisolvable in this case as well.
This completes the proof of the lemma.
\end{proof}

\begin{theorem} Let $H$ be a semisimple Hopf algebra of dimension $24$.
Then $H$ is upper and lower semisolvable. \end{theorem}

\begin{proof} It will be enough to show that $H$ is not simple.
This follows from Remark \ref{rmk-24} (i), Lemmas \ref{24-3},
\ref{24-6} and Corollary \ref{24-8}, if $H$ is of type $(1, 12; 2, 3)$, $(1, 3; 2, 3; 3, 1)$, $(1, 6;  3, 2)$, $(1, 8; 2, 4)$ or $(1, 8; 4, 1)$,  respectively. Also, from Lemma \ref{24-2}, $H$ is not simple
if $H$ is of type $(1, 2; 2, 1; 3, 2)$ as a coalgebra. Hence, we may assume that $H$ is of type $(1, 4; 2, 5)$  as a coalgebra, and by
Remark \ref{rmk-24} (iii), $H$ contains a Hopf subalgebra of
dimension $8$. Applying Lemma \ref{hopfsb-8}, we find that
$H^*$ is of type $(1, 2; 2, 1; 3, 2)$ as a coalgebra, and $H$ is not simple by Lemma \ref{24-2}.
\end{proof}

\chapter{Dimension $30$}\label{30}

In this chapter $H$ will denote a semisimple Hopf algebra of dimension $30$
over $k$.
By \cite[Theorem 4.6]{pqq}, if $H$ fits into an abelian extension
$$1 \to k^{\Gamma} \to H \to kF \to 1,$$
where $\Gamma$ and $F$ are finite groups, then $H$ is trivial.

\section{Possible (co)-algebra structures} We shall assume in this section that $H$ is nontrivial, and reduce the possibilities for its algebra and coalgebra structures.

\begin{lemma}\label{conteo} The group $G(H^*)$ of group-like elements in $H^*$
is of order $2$, $3$,
$5$, $6$ or $10$.  As an algebra, $H$ is of one of the following types:
\begin{itemize}  \item $(1, 2; 2, 7)$, \item $(1, 3; 3, 3)$, \item $(1, 5; 5, 1)$, \item $(1, 6; 2, 6)$, \item $(1, 10; 2, 5)$. \end{itemize} \end{lemma}

\begin{proof} It follows from  \ref{alg-struct} that $n: = |G(H^*)| \neq 15$, and if $n
= 3, 5, 6$ or $10$, these are the only possible algebra
types for $H$. Suppose that $n = 1$; again by \ref{alg-struct}, it follows
 that the set of degrees of non-linear irreducible representations of $H$ has
at least three elements and necessarily $H$ has an irreducible
module of degree $2$. This is impossible by Corollary \ref{cor-g-1}.
Suppose finally
that $n = 2$. Then, as an algebra, $H$ must be either of type $(1, 2; 2, 7)$  or $(1, 2; 2, 3; 4, 1)$. By Theorem \ref{thm-nr},  $G[\chi] = G(H^*)$ is of order $2$, for all irreducible
character $\chi$ such that $\chi (1) = 2$. The second possibility implies, by Theorem \ref{corolario}, that $H$ has a quotient Hopf algebra of dimension $2 +
4.3 = 14$, which contradicts \cite{NZ}. This completes the proof of the lemma. \end{proof}

\begin{remark}\label{modulo-gl} (i) Observe that if $H$ is of type $(1, 5; 5, 1)$ as an algebra, the irreducible character of degree $5$ is stable under multiplication by $G(H^*)$. Also, if $H$ is of type $(1, 3; 3, 3)$ as an algebra, every irreducible character $\psi$ of degree 3 is stable under multiplication by $G(H^*)$: this can be seen decomposing the product $\psi\psi^*$ into irreducibles and using the relation \eqref{des-ss}.

In the other cases, after Theorem \ref{thm-nr}, all irreducible
characters $\chi$ of degree $2$ have nontrivial isotropy, {\it
i.e.}, $G[\chi]$ is of order $2$ for all such characters.
Combining  this observation with Corollary \ref{cor-cociente} and
Remark  \ref{particular} we find that, in any case, the quotient
coalgebra $H^*/H^*(kG(H^*))^+$ is cocommutative.

(ii) Suppose that $H$ is of type $(1, 10; 2, 5)$ or $(1, 6; 2, 6)$ as an algebra. Then the group $G(H^*)$ is abelian.

\begin{proof} By part (i), for all irreducible character $\chi$ of degree $2$ we have $G[\chi] \neq 1$. The claim follows from Proposition \ref{nab-pq}. \end{proof} \end{remark}

\begin{lemma}\label{30-2-comm} Suppose that $H$ is of type $(1, 2; 2, 7)$ as a coalgebra. Then $H$ is commutative. \end{lemma}

\begin{proof} Note first that if $G(H) \subseteq Z(H)$, in view of Remark \ref{modulo-gl}, $H$ fits
into an abelian extension $1 \to kG(H)  \to H \to H/H(kG(H))^+ \to
1$, thus implying the lemma, in view of \cite{pqq}.

By Lemma \ref{conteo}, $|G(H^*)| \neq 1$.  If $|G(H^*)|$ is even, the lemma follows from Theorem \ref{1-2-2-n} and \cite{pqq}.  We may thus assume that $|G(H^*)|$ is odd.

Consider the projection $q: H \to B$, where $B = k^{G(H^*)}$.
Since the dimension of $B$ is odd and $G(H)$ is of order $2$, we
must have $G(H) \subseteq H^{\co B}$. Therefore, as a left coideal
of $H$, $H^{\co B}$ decomposes in the form $H^{\co B} = kG(H)
\oplus V_1 \oplus \dots \oplus V_m$, where $V_i$ is an irreducible
left coideal of $H$ of dimension $2$, for all $i = 1, \dots, m$.
For $1 \leq i \leq m$, let  $C_i$ be the simple subcoalgebra of
$H$ containing $V_i$. By Remark \ref{modulo-gl}, we have $g C_i =
C_i = C_i g$, for all $g \in G(H)$; in particular, $kG(H)
\subseteq k[C_i]$.

\begin{claim} $G(H) \subseteq Z(k[C_i])$, for all $i = 1, \dots, m$.
\end{claim}

\begin{proof} Observe that $V_i$ appears in $H^{\co B}$ with multiplicity at most 2,
and moreover that $V_i$ appears with multiplicity exactly 2 if and
only if $C_i \subseteq H^{\co B}$. Consider first the case where
$m(V_i, H^{\co B}) = 2$. Then $k[C_i] \subseteq H^{\co B}$,
whence $\dim k[C_i] < \dim H$. Then $k[C_i]$ is commutative and
the claim follows.

Suppose now that $m(V_i, H^{\co B}) = 1$. Let $V = V_i$. Note that $gV$ and $Vg$ are irreducible left coideals of $H$ isomorphic to $V$, and  $gV, Vg \subseteq H^{\co B}$. The multiplicity condition on $V$ implies that $gV = V = Vg$, and the claim follows in this case from Corollary \ref{cor-mk}. \end{proof}

Let $A = k[C_i: \quad 1 \leq i \leq m]$. This is a Hopf subalgebra
of $H$, and the claim implies that $kG(H) \subseteq Z(A)$. Hence,
we may assume that $A \subsetneq H$. Since also $H^{\co B}
\subseteq A$, a dimension argument implies that $H^{\co B} = A$
and $A$ is commutative. Thus the map $q$ is normal, and $H$ fits
into an abelian extension $1 \to A \to H \to k^{G(H^*)} \to 1$.
This implies the lemma in view of \cite{pqq}. \end{proof}

\begin{lemma}\label{n-neq-3} Suppose that $H$ is of type $(1, 3; 3, 3)$ as a coalgebra. Then $H$ is commutative. \end{lemma}

\begin{proof} By Lemma \ref{30-2-comm},  $H$ is not of type $(1, 2; 2, 7)$ as an algebra.
Suppose first that $|G(H^*)|$ is divisible by $3$. Then, by Lemma
\ref{spl-ins}, $H$ is isomorphic to a biproduct $H \simeq R \#
kG(H)$, where $R$ is cocommutative by Remark \ref{modulo-gl} (i).
Then the lemma follows in this case by Proposition \ref{som} and
\cite{pqq}. Hence we may assume that $H$ is not of type $(1, 3; 3,
3)$ nor $(1, 6; 2, 6)$ as an algebra. Similarly, if $H$ is of type
$(1, 10; 2, 5)$ as an algebra, then $H$ fits into an  exact
sequence $1 \to kG(H) \to H \to k^{G(H^*)} \to 1$, and $k^{G(H^*)}
\simeq H/H(kG(H))^+$ is cocommutative. Hence the result follows by
\cite{pqq}.

Suppose finally that $H$ is of algebra type $(1, 5; 5, 1)$. Let $C$ be the unique simple subcoalgebra  of dimension $25$ of $H^*$. It is clear that $H^* = k[C]$ is generated by $C$ as an algebra.
Consider the projection $q: H^* \to k^{G(H)}$; we have $\dim (H^*)^{\co q} = 10$.
Since $|G(H)|$ and $|G(H^*)|$ are relatively prime, we have $G(H^*) \subseteq  (H^*)^{\co q}$.
By dimension restrictions, $(H^*)^{\co q} = kG(H^*) \oplus V$, where $V \subseteq C$ is an irreducible left coideal of dimension 5. Let $g$  be a generator of $G(H^*)$; since $V$ is the only 5-dimensional irreducible coideal contained in $(H^*)^{\co q}$, then we must have $gV = V = Vg$. By Corollary \ref{cor-mk} $kG(H^*)$ is normal in $H^*$. This implies the claim in this case, in view of \cite{pqq}, since the quotient Hopf algebra $H^* / H^*(kG(H^*))^+$ is cocommutative, by Remark \ref{modulo-gl} (i).
 This completes the proof of the lemma.
\end{proof}

\begin{lemma}\label{6-arrow-6} Suppose that $H$ is of type $(1, 6; 2, 6)$ as an algebra. Then $H$ is of type $(1, 6; 2, 6)$ as a coalgebra.
\end{lemma}

\begin{proof} Suppose not.
It follows from Lemmas \ref{conteo}, \ref{30-2-comm} and
\ref{n-neq-3},  that $G(H)$ has a subgroup $T$ of order $5$. Since
$|G(H^*)| = 6$, we have a sequence of Hopf algebra maps
\begin{equation}\label{sucesion}\begin{CD}kT @>{\iota}>> H @>{\pi}>>
k^{G(H^*)}.\end{CD} \end{equation}
 In particular, $\pi (a) = 1$, for all $a \in T$, implying that $kT \subseteq  H^{\co \pi}$, and since $\dim H^{\co \pi} = [H^*:kG(H^*)] = |T|$, we have $kT =  H^{\co \pi}$. Therefore
\eqref{sucesion} is an exact sequence of Hopf algebras. By Remark \ref{modulo-gl} (ii),
$kT$ and $k^{G(H^*)}$ are both commutative and cocommutative; hence the extension is abelian. This implies that $H$ is trivial, against the assumption. This finishes the proof of the lemma. \end{proof}

\begin{lemma}\label{z5} Suppose that $H$ is of type $(1, 5; 5, 1)$ as a coalgebra. Then $H$ is commutative. \end{lemma}

\begin{proof} By Lemmas \ref{conteo}, \ref{30-2-comm},  \ref{n-neq-3} and \ref{6-arrow-6}, $G(H^*)$ is
of order $5$ or $10$. Let $\Gamma \subseteq G(H^*)$ be a subgroup
of order $5$, and consider the sequence of Hopf algebra maps
$\begin{CD} kG(H) @>{\iota}>> H @>{\pi}>> k^{\Gamma}\end{CD}$.
Then $\dim H^{\co \pi} = [H^*:k\Gamma] = 6$ and $kG(H) \cap H^{\co
\pi} = k1$ by \cite{NZ}. Then the composition $\pi \iota : kG(H)
\to k^{\Gamma}$ is an isomorphism, and therefore $H$ is isomorphic
to a biproduct $R \# k G$, where $G = G(H)$ is a group of order
$5$ and $R$ is a semisimple Hopf algebra of dimension $6$ in the
category of  Yetter-Drinfeld modules over $G$.

In particular, $R \simeq H/H(kG(H))^+$ as a coalgebra, and thus
$R$ is cocommutative by Remark \ref{modulo-gl} (i). By Proposition
\ref{som},  $H$ fits into an abelian extension and the result
follows from \cite{pqq}.  \end{proof}

\begin{lemma}\label{z10} Suppose that $H$ is of type $(1, 10; 2, 5)$ as a coalgebra. Then $H$ is commutative. \end{lemma}

\begin{proof} By  previous lemmas $H$ is, as an algebra,
necessarily of type  $(1, 10; 2, 5)$. In view of Lemma \ref{spl-ins}, this implies that $G(H) \simeq G(H^*)$ are
abelian and $H$ is isomorphic to a biproduct $R  \# k G(H)$, where
$\dim R = 3$. Thus $R$ is cocommutative and commutative, and since the action
of $G(H)$ on $R$ is by coalgebra automorphisms, we have that the
unique subgroup of order $5$, $\Gamma$, of $ G(H)$ acts
trivially on $R$.

Then $\Gamma \subseteq Z(H)$, and there is an exact sequence $1
\to k\Gamma \to H \to \overline H \to 1$, where $\dim \overline H
= 6$. Note that $\overline H$ is cocommutative; otherwise
$\overline H \simeq k^{F}$, where $F$ is the only non-abelian
group of order $6$. But this implies that $F = G((k^{F})^*)$ is a
subgroup of $G(H^*)$, which is a contradiction. The lemma follows
from \cite{pqq}. \end{proof}

\section{Classification}
In this section we aim to give a proof of Theorem \ref{cls30}.  For this, we shall consider the possible
structures given by Lemma \ref{conteo24}.

{\it Proof of Theorem \ref{cls30}.} By Lemmas \ref{conteo}, \ref{modulo-gl}, \ref{z5} and \ref{z10}, we may assume that $H$ and $H^*$ are both of type $(1, 6; 2, 6)$ as coalgebras.  In view of Remark \ref{modulo-gl} (ii), the group $G(H)$ is cyclic.

Let $F$ be the unique subgroup of order $2$ of $G(H)$. By Remark
\ref{modulo-gl} (i) all irreducible characters of degree $2$ are
stable under multiplication by $F$, and therefore the quotient
coalgebra $H/H(kF)^+$ is cocommutative. Clearly, $H$ is a
biproduct $R \# kF$, and $R \simeq H/H(kF)^+$ is a cocommutative
coalgebra. This implies that $H$ fits into an abelian extension in
view of Proposition \ref{som}, and  $H$ is trivial by \cite{pqq}.
\qed

\chapter{Dimension $36$}\label{36}

\section{Reduction of the problem} Let $H$ be a nontrivial semisimple Hopf algebra of dimension $36$ over $k$.

\begin{lemma}\label{conteo36} The order of $G(H^*)$ is either $2$, $3$, $4$, $6$, $9$,
$12$ or $18$ and as an algebra $H$ is of one of the following types:
\begin{itemize}\item $(1, 2; 2, 4; 3, 2)$, \item $(1, 3; 2, 6; 3, 1)$, \item $(1, 4; 2, 8)$, \item $(1, 4; 2, 4; 4, 1)$, \item $(1, 4; 4, 2)$, \item $(1, 6; 2, 3; 3, 2)$, \item $(1, 9; 3, 3)$, \item $(1, 12; 2, 6)$, \item $(1, 18; 3, 2)$.  \end{itemize} \end{lemma}

By \cite{twist-simple} there is a simple Hopf algebra $\mathcal H$
with algebra and coalgebra type $(1, 4; 2, 4; 4, 1)$. By
construction, $\mathcal H = (kG)_{\phi}$ is a twisting of the
group algebra of $G = D_3 \times D_3$ with respect to a
non-degenerate 2-cocycle $\phi \in \mathbb Z_2 \times \mathbb
Z_2$. Moreover, $\mathcal H$ is the only twisting of a group of
order $36$ that is simple, and we have $\mathcal H \simeq \mathcal
H^{*\cop}$ \cite[4.3]{twist-simple}.

\begin{proof} It follows from \ref{alg-struct} that $G(H^*) = 1$ is impossible.
It follows as well that the only possibilities for the algebra
types are the prescribed  ones if  $|G(H^*)|$ equals either $3$,
$4$, $6$, $9$, $12$ or $18$.

  Finally, if $|G(H^*)| = 2$, a counting argument gives
that $H$ must be of type $(1, 2; 2, 4; 3, 2)$ or $(1, 2; 3, 2; 4,
1)$. In the last case, $H$ has two characters of degree $1$,
$\epsilon$ and $g$, two irreducible characters of degree $3$,
$\psi_1$ and $\psi_2$, and one irreducible character of degree
$4$, $\zeta$.  Since $G[\psi_i] = 1$, because $G(H^*)$ is of order
$2$, then $g \psi_1 = \psi_2$; on the other hand,  $g \zeta =
\zeta$ because $\zeta$ is the only irreducible character of degree
$4$. Counting dimensions, we obtain $\psi_i \psi_i^* = \epsilon +
2 \zeta, \quad i = 1, 2$. Thus $m(\psi_i, \zeta \psi_i) = m(\zeta,
\psi_i \psi_i^*) = 2$, and therefore
$$\zeta \psi_1 = 2 \psi_1 + r \psi_2 + s \zeta = g \zeta \psi_1 = 2 \psi_2 + r
\psi_1 + s \zeta, $$
implying that $r = 2$ and $\zeta \psi_1 = 2 \psi_1 + 2 \psi_2 = \zeta \psi_2$.

Since $*$ fixes $\zeta$ and permutes the set $\{ \psi_1, \psi_2
\}$, and also since $\psi_i \zeta = (\zeta \psi_i^*)^*$, we find
that $m (\zeta, \psi_i \zeta) = 0$, $i = 1, 2$. But $m(\zeta,
\psi_i \zeta) = m(\psi_i,  \zeta \zeta^*)$, so that $\zeta
\zeta^*$ decomposes in the form $\zeta \zeta^* = \epsilon + g + r
\zeta, \quad r \in \mathbb Z^+$. Taking degrees we see that this
is impossible. This discards the possibility $(1, 2; 3, 2; 4, 1)$
for the algebra type of $H$ and finishes the proof of the lemma.
\end{proof}

\begin{remark}\label{rmkk-36} (i) Suppose that  $H$ is of type $(1, 18; 3, 2)$.  Then,  by Corollary \ref{kob-mas}, $G(H^*) \cap Z(H^*) \neq 1$.

 (ii) Suppose that $H$ is of type $(1, 2; 2, 4; 3, 2)$ as a coalgebra. Since $|G(H)| = 2$, $H$ contains no Hopf subalgebra of dimension $12$. By Remark \ref{nr} (i), the irreducible characters of degrees $1$ and $2$ give rise to a  Hopf subalgebra of dimension $18$. Therefore, in view of Corollary \ref{kob-mas}, $G(H^*) \cap Z(H^*) \neq 1$.

 (iii) Suppose that $H$ is of type $(1, 6; 2, 3; 3, 2)$ as a coalgebra. Then for every irreducible character $\chi$ of degree $2$, we must have $|G[\chi]| = 2$ (because there are only $3$ such characters) and $H$ has no simple subcoalgebra of dimension $16$.
Hence, also here, the irreducible characters of degrees $1$ and $2$ give rise to a  Hopf subalgebra of dimension $18$. Therefore $G(H^*) \cap Z(H^*) \neq 1$.

(iv) If $H$ is of type $(1, 4; 4, 2)$ as a coalgebra, then  $H$ does
not contain any  left coideal subalgebra of dimension $3$.
Therefore $H^*$ contains no Hopf subalgebra of dimension $12$.
\end{remark}

\begin{lemma}\label{36-3} Suppose that $H$ is of type $(1, 3; 2, 6; 3, 1)$ as a coalgebra. Then $H$ is not simple. \end{lemma}

\begin{proof} By Theorem \ref{thm-nr}, $H$ has a unique Hopf subalgebra $A$ of dimension $12$, which contains $G(H)$ and the unique simple subcoalgebra of dimension $9$ of $H$.

The classification of semisimple Hopf algebras of dimension $12$ \cite{fukuda} implies that $A$ is commutative. Therefore $\dim V \leq [H:A] = 3$ for all irreducible $H$-module $V$; see for instance Corollary 3.9  in \cite{harmonic}. Then, by the previous remark, we may assume that $H^*$ is, as a coalgebra, of one of the types $(1, 3; 2, 6; 3, 1)$, $(1, 4; 2, 8)$, $(1, 9; 3, 3)$ or $(1, 12; 2, 6)$.

Suppose that $H^*$ is of type $(1, 9; 3, 3)$ as a coalgebra. We have a Hopf algebra projection $q: H^* \to A^*$, such that $\dim {H^*}^{\co A^*} = 3$. This projection is necessarily normal, thus $H$ is not simple in this case.

If $H^*$ is of type $(1, 3; 2, 6; 3, 1)$ as a coalgebra, then also $H^*$ contains a commutative Hopf subalgebra $B \simeq A$ of dimension $12$.
We have a Hopf algebra projection $q: H \to B^*$, for which $\dim H^{\co B^*} = 3$. Then we may assume $A \cap H^{\co B^*} = k1$ and the restriction $q\vert_A : A \to B^*$ is an isomorphism. Thus $A$ is cocommutative, which is absurd. This contradiction discards this type.

Similarly, if $H^*$ is of type $(1, 12; 2, 6)$ as a coalgebra,
then $H$ is a biproduct  $H = R \# A$, where $R$ is a braided Hopf
algebra over $A$ of dimension $3$. Since the irreducible
$A$-comodules are of dimension $1$ or $3$, we must have $\rho (R)
\subseteq kG(A) \otimes R$. This implies that $R$ is a braided
Hopf algebra over $kG(A)$ and the biproduct $K = R \# kG(A)$ is a
Hopf subalgebra of $H$. But $\dim K = 9$ and $K$ is thus
cocommutative, which contradicts $|G(H)| = 3$. Thus this type is
not possible.

Finally suppose that $H^*$ is of type $(1, 4; 2, 8)$ as a
coalgebra. Consider  the surjective Hopf algebra map $H^* \to
A^*$, so that $\dim (H^*)^{\co A^*} = 3$ and we may assume that
$(H^*)^{\co A^*} = k1 \oplus V$ as a left coideal of $H^*$, where
$V$ is an irreducible coideal of dimension $2$. By Theorem
\ref{coinvariantes} $H^*$ contains a Hopf subalgebra $B$ of
dimension $12$ such that $B$ is not cocommutative. As in the
previous paragraphs, it turns out that $A \simeq B^*$, which is a
contradiction. This finishes the proof of the lemma.
\end{proof}

\begin{lemma}\label{l-36g} Suppose that $H$ is of type $(1, 9; 3, 3)$ as a coalgebra. Then $H$ is not simple. \end{lemma}

\begin{proof} If $|G(H^*)| = 9$, then $H = R \# kG(H)$ is a biproduct and the lemma follows from Proposition \ref{dimR-3-4}. Thus, by Remark \ref{rmkk-36} and Lemma \ref{36-3}, we may assume that $|G(H^*)| = 12$ or $4$. In any case, there is a quotient $H \to k\Gamma$, where $\Gamma$ is a group of order $4$, for which necessarily $kG(H) = H^{\co k\Gamma}$. Hence $H$ is not simple. \end{proof}

\begin{lemma}\label{red-36} Suppose that $H$ is simple. Then the order of $G(H^*)$ equals either $4$ or $12$.
Moreover, $H$ is isomorphic to a biproduct $H \simeq
R \# kG$, where $G$ is a group of order $4$ and $R$ is a Yetter-Drinfeld Hopf
algebra over $G$ of dimension $9$. \end{lemma}

\begin{proof} By Remark \ref{rmkk-36} and Lemmas \ref{36-3} and \ref{l-36g}, $|G(H^*)|$ is not $18$, $9$, $6$, $3$ or $2$.
Therefore, by Lemma \ref{conteo36}, $|G(H^*)| = 4$ or $12$.

Let  $\Gamma \subseteq G(H^*)$ be a subgroup of order 4.
Consider the quotient $p: H \to k^{\Gamma}$. Then $\dim H^{\co p} = 9$. Let $G$
be a subgroup of $G(H)$ of order $4$. Then, by \cite{NZ}, $kG \cap H^{\co p} =
k1$ and thus the restriction of $p$ gives an isomorphism $kG \to k^{\Gamma}$. This implies the lemma. \end{proof}

\begin{lemma}\label{subalg-36} Suppose that $H$ is of type $(1, 4; 2, 4; 4, 1)$ as a coalgebra. Then every irreducible left coideal of dimension $2$ in $H$ is contained in a Hopf subalgebra of dimension $12$.

Moreover, $H$ contains two Hopf subalgebras $B$ and $B'$, which are of type $(1, 4; 2, 2)$ as coalgebras and such that

 (i) $B \cap B' = kG(H)$ and $G(H) \simeq \mathbb Z_2 \times \mathbb Z_2$;

 (ii) $B \cup B'$ generates $H$ as an algebra. \end{lemma}

\begin{proof} Since $H$ contains no irreducible coideals of dimension $3$, then for all  irreducible
 $\chi$ such that $\deg \chi = 2$, we have $G[\chi ] \neq
1$.

\begin{claim} The action of $G(H)$ by left multiplication on the set $X_2$ of irreducible characters of degree $2$ has two disjoint orbits $\{ \chi, \chi' \}$ and $\{ \psi, \psi' \}$.
In addition, we have $G[\chi] = G[\chi']$ and $G[\psi] = G[\psi']$ are distinct subgroups of order $2$ of $G(H)$. \end{claim}

\begin{proof}  Note that the tensor product $\chi \chi'$ is irreducible for some irreducible characters $\chi$ and $\chi'$ of degree $2$, since otherwise, there would be a Hopf subalgebra of dimension $4 + 4 . 4 = 20$, which is impossible.
Therefore there exist $\chi, \psi \in X_2$ such that $G[\chi]$ is not contained in $G[\psi]$; c.f. \ref{productodesimples}. Thus $G[\chi]$ and  $G[\psi]$ are distinct subgroups of order $2$ of $G(H)$.  In particular, the group $G(H)$ is isomorphic to $\mathbb Z_2 \times \mathbb Z_2$.

Since $|G[\chi]| = |G[\psi]| = 2$, then the orbits $G(H) \chi$ and $G(H) \psi$ are of order $2$. Moreover, since $G[\chi'] = G[\chi]$, for all $\chi' \in G(H)\chi$ (because $G[\chi]$ is normal in $G(H)$), then $G(H) \chi$ and $G(H) \psi$ are disjoint orbits. This proves the claim. \end{proof}

\begin{claim} We have $\chi^* = \chi$, for all $\chi \in R(H^*)$. \end{claim}

\begin{proof}
We have already shown that $a^2 = 1$, for all $a \in G(H)$.
If $\chi$ is irreducible of degree $2$, then we have $\chi \chi^* = \sum_{a \in G[\chi]} a + \lambda$, for some $\lambda \in X_2$ such that $\lambda^* = \lambda$. Also $G[\chi] \subseteq G[\lambda]$, which implies that $\lambda \in G(H) \chi$.
Hence every orbit has a self-dual element.
Moreover, if $g \in G(H)$ is such that $g \lambda = \chi$, then $\chi^* = \lambda g$ and therefore $G[\chi] = G[\lambda] \subseteq G[\chi^*]$.
Hence $\chi^* \in G(H) \chi$. This implies, since each orbit has two elements, that $\chi^* = \chi$. Since $\chi$ was arbitrary, the claim follows. \end{proof}

The claim implies the character algebra $R(H^*)$ is commutative.

\smallbreak As before, we have decompositions $\chi^2 = \chi
\chi^* = \sum_{a \in G[\chi]} a + \lambda$, where $\lambda \in
G(H) \chi$. On the other hand, if $\chi' \in G(H) \chi$, say $b
\chi = \chi'$, then
 $\chi' \chi =  b \chi^2 = b + ba + b \lambda$, and $b \lambda  \in G(H) \chi$.
Thus, the set $G(H) \cup G(H) \chi$ spans a standard subalgebra of
$R(H^*)$ which corresponds to a Hopf  subalgebra $A$ of dimension
$12$  such that $R(A^*) = k(G(H) \cup G(H) \chi$.

The same argument applies to the other orbit and concludes the proof of the lemma. \end{proof}

\begin{lemma}\label{subalg-36-bis} Suppose that $H$ is of type $(1, 4; 2, 8)$ as a coalgebra. Assume in addition that $H$ is simple. Then  $H$ contains two Hopf subalgebras $B$ and $B'$, which are of type $(1, 4; 2, 2)$ as coalgebras and such that $B$ and $B'$ satisfy the conditions (i) and (ii) of Lemma \ref{subalg-36}. \end{lemma}

\begin{proof} First note that it is enough to show that $H$ contains Hopf subalgebras $B, B'$ of dimension $12$ such that $B \cap B' = kG$. Indeed, in this case, a dimension argument implies that necessarily the subalgebra generated by $B$ and $B'$ coincides with $H$.
On the other hand, this also implies that the group $G = G(H)$ is
isomorphic to $\mathbb Z_2 \times \mathbb Z_2$: otherwise $G$
would contain a unique subgroup of order $2$, which would be
central both in $B$ and in $B'$ (\cite{fukuda}) and a fortiori in
$H$.

  By Lemma \ref{red-36}  $H = R \# kG$, where $G = G(H)$
and $\dim R = 9$. Therefore $R$ decomposes in the form $R = k1
\oplus V_1 \oplus V_2 \oplus V_3 \oplus V_4$ as a left coideal of
$H$, where $V_i$ is an irreducible left coideal of dimension $2$,
for all $i = 1, \dots, 4$.

\smallbreak It follows from Corollary \ref{formula} and
Proposition \ref{R-coideal} that $V_i$ is a subcoalgebra of $R$,
for all $i$,  and $V_i$ is not isomorphic to $V_j$, as left
coideals of $H$, if $i \neq j$. Also $|G[\chi_i]| = 2$, where
$\chi_i = \chi_{V_i}$. In particular, $\# G \chi_i = 2$, for all
$i$. We have $V_i \# kG = C_i \oplus C_ih$, where $C_i$  is the
simple subcoalgebra of $H$ containing $V_i$ and $C_ih \neq C_i$.
Thus, if $V_j \simeq V_i h$ then $V_j \# kG = V_i \# kG$ and $V_j
= (\id \otimes \epsilon) (V_j \# kG) = V_i$. Therefore $\chi_i$,
$1 \leq i \leq 4$,  form a system of representatives of the orbits
of the action of $G$ on $X_2$ by right multiplication.

  The multiplicity-one condition for the  decomposition of
$R$ implies that $V_1, \dots, V_4$ are uniquely determined, and
these are the only irreducible left coideals of $H$ contained in
$R$. Thus the action of $G$ on $R$ must permute the $V_i$'s.

\begin{claim}\label{yd-36} Let $1 \leq i \leq 4$.
There  exists a subgroup $1 \neq S$ of $G$ such that $s . V_i =
V_i$. \end{claim}

\begin{proof} The possible decompositions of $R$ as an $H$-Yetter-Drinfeld submodule of $H$ are the following
$$k1 \oplus X_1 \oplus X_2, \quad k1 \oplus Y_1 \oplus Y_2 \oplus Y_3 \oplus Y_4, \quad k1 \oplus X_1 \oplus Y_1 \oplus Y_2, \quad k1 \oplus Y_1 \oplus Z, $$
where $\dim X_i = 4$, $\dim Y_i = 2$, $\dim Z = 6$.
Decomposing the $X_i$'s and $Z$ as left coideals of $H$, and using that the action of $G$ permutes the $V_i$'s, the claim follows. \end{proof}

Let $1 \neq S$ be a subgroup of $G$ such that $s.V_1 = V_1$, for all $s \in S$, as in Claim \ref{yd-36}. The action of $S$ must  permute either transitively or trivially the set $G(V_1)$.

Suppose first that the action of $S$ on $V_1$ is trivial. Let $V_1' = G.V_1$. Then the action of $S$ on $V_1'$ is also trivial and so is the action of $S$ on the subalgebra $k[V_1']$ generated by $V_1'$.

Since $k[V_1']$ is a Yetter-Drinfeld Hopf subalgebra of $R$, $\dim
k[V_1']$ is either $3$ or $9$. The last possibility implies that
the action of $S$ on $R$ is trivial. Therefore, we may assume that
$\dim k[V_1'] = 3$, and then $V_1$ is stable under the action of
$G$. Thus $B = k[V_1] \# kG$ is a Hopf subalgebra of $H$ of
dimension $12$. Moreover, since the action of $G$ permutes the set
$\{  V_1, \dots, V_4 \}$ and by the above fixes $V_1$, we may
assume that also $gV_2g^{-1} = V_2$, for all $g \in G$. Then the
Hopf subalgebras $B$ and $B' = k[V_2] \# kG$ do the job.

Suppose finally that the action of $S$ on $V_1$ is transitive. By Lemma \ref{tr-act} $|G_{V_1}| = 2$.
Let $V_1' :  = G_{V_1} . V_1$ as above; note that $G_{V_1'} = G_{V_1}$. Then, by Lemma \ref{util},
 $\widetilde B = k[V_1'] \# kG_{V_1} = k[V_1'] \# kG_{V_1'}$ is a Hopf subalgebra of $H$.
If $k[V_1'] = R$, then $\widetilde B$ is an index-$2$ (and
therefore normal) Hopf subalgebra of $H$. So we may assume that
$\dim k[V_1'] = 3$, implying that $g.V_1 = V_1$ for all $g \in G$.
We may also assume that  $g.V_2 = V_2$ for all $g \in G$. Let  $B
=  k[V_1] \# kG$ and $B' = k[V_2] \# kG$. Then we have $\dim B =
\dim B' = 12$, and we are done. \end{proof}

\begin{lemma}\label{36-4} Suppose that $H$ is of coalgebra type $(1, 4; 2, 4; 4, 1)$ or $(1, 4; 2, 8)$.
Assume in addition that $H$ is simple. Then $H \simeq \mathcal H$.
\end{lemma}

Here $\mathcal H$ denotes the simple (self-dual) example in
\cite{twist-simple}.

\begin{proof}
Let  $B$ and $B'$ be the Hopf subalgebras of $H$ as in  Lemmas
\ref{subalg-36} and \ref{subalg-36-bis}, respectively. Then
neither $B$ nor $B'$ is cocommutative.

\begin{claim} We have $B \simeq B' \simeq \mathcal A_0$; see Chapter \ref{twist}. \end{claim}

\begin{proof} Suppose not.
Since  $G(H) \simeq \mathbb Z_2 \times \mathbb Z_2$, then $B$ and
$B'$ are not isomorphic to $\mathcal A_1$. In virtue of the
classification in dimension $12$ \cite{fukuda}, at least one of
them, say $B$, is commutative. Then there exists $1 \neq g \in
Z(B') \cap G(B')$. Since $G(B') = G(B) \subseteq B$, $g$ also
commutes with $B$. Since $B$ and $B'$ generate $H$ as an algebra,
this implies that $g$ is a central group-like element in $H$,
contradicting the assumption. \end{proof}

  We know from Chapter \ref{twist} that there exists a
normalized $2$-cocycle $\phi \in kG(H) \otimes kG(H)$ such that
$B_{\phi} \simeq {B'}_{\phi} \simeq kG$, where $G$ is the group in
Proposition \ref{twist-pqq}. Moreover, $B_{\phi}$  and
${B'}_{\phi}$ are Hopf subalgebras of $H_{\phi}$ and $B_{\phi}
\cup {B'}_{\phi}$ generates $H_{\phi}$ as an algebra, because the
multiplication was unchanged. It turns out that $H_{\phi}$ is
cocommutative.

More precisely,  $H_{\phi} = kG'$, where $G' = F' \rtimes \Gamma$
is a semidirect product with $|F'| = 9$ and $\Gamma = G(H)$. Since
$H$ is simple by assumption, then by the results in
\cite[4.3]{twist-simple}, $G \simeq D_3 \times D_3$ and $H \simeq
\mathcal H$.
\end{proof}

\begin{lemma}\label{36-12} Suppose that $H$ is of type $(1, 12; 2, 6)$ as a coalgebra. Then $H$ is not simple. \end{lemma}

\begin{proof} By Lemma \ref{red-36}, we may assume that $|G(H^*)| = 4$ or $12$. By Lemma \ref{36-4}, if $|G(H^*)| = 4$ we may also assume that $H^*$ is of type $(1, 4; 4, 2)$ as a coalgebra; this is impossible in view of Remark \ref{rmkk-36} (iv). Therefore, we must consider the case $|G(H^*)| = 12$, {\it i.e.}, $H$ and $H^*$ are both of type $(1, 12; 2, 6)$.

Consider the quotient Hopf algebra $p: H \to k^{G(H^*)}$. We have
$\dim H^{\co p} = 3$. If $g \in H^{\co p}$ for some $g \in G(H)$,
then by \cite{NZ} $g$ is of order $3$ and $H^{\co p} = k \langle g
\rangle$ is a normal Hopf subalgebra of $H$, implying that $H$ is
not simple. Thus, we  may assume that  $p(g) \neq 1$, for all $g
\in G(H)$. This implies that the restriction $p\vert_{kG(H)}:
kG(H) \to k^{G(H^*)}$ is an isomorphism. Therefore the groups
$G(H)$ and $G(H^*)$ are isomorphic and abelian and $H$ is a
biproduct $H = R \# kG(H)$, where $\dim R = 3$. The lemma follows
from Proposition \ref{dimR-3-4}. \end{proof}

\section{Main result}
In view of the previous results, the only case it remains to
consider is that where $H$ and $H^*$ are both of type $(1, 4; 4,
2)$ as coalgebras.

\begin{theorem} Let $H$ be a semisimple Hopf algebra of dimension
$36$. Suppose that $H$ is simple. Then $H \simeq \mathcal H$.
\end{theorem}

\begin{proof} By Lemma \ref{red-36}, $H$ is a biproduct $H = R \# kG(H)$,
where $R$ is a left coideal subalgebra of dimension $9$.
Therefore, as a left coideal of $H$, $R = k1 \oplus V_1 \oplus
V_2$, where $V_i$  are irreducible left coideals of dimension $4$,
$i = 1, 2$.

\begin{claim}\label{af-cc-tw} There exists a normalized invertible $2$-cocycle $\phi \in kG(H) \otimes kG(H)$ such that $H_{\phi}$ is not simple. \end{claim}

\begin{proof} We have $g V_i \simeq V_i \simeq V_i g$ as a left coideal of $H$, for all $g \in G(H)$; so that by Corollary \ref{formula}, $V_i$ is a subcoalgebra of $R$ in the category of Yetter-Drinfeld modules over $G(H)$, and we have $V_i \# kG(H)$ is a simple subcoalgebra of dimension 16. It is also clear that $V_i$ generates $R$ as an algebra, because $k[V_1]$ is a braided Hopf subalgebra of $R$ containing $V_i$; then $\dim V_i^G = 1$ by \cite[Lemma 4.4.1]{clspqq}.

Write $G = G(H)$. By \cite[1.3]{pqq2} there exist normalized 2-cocycles $\alpha_i: G\times G \to G$ such that $V_i \simeq (k_{\alpha_i}G)^*$, $i = 1, 2$, and moreover, in view of the results in \cite[Section 4]{clspqq}
 there exists a normalized invertible $2$-cocycle $\phi \in kG \otimes kG$ such that $H_{\phi} = R_{\phi} \# kG$, and $V : = V_1$ is a cocommutative subcoalgebra of $R_{\phi}$. See \cite[Remark 4.4.3]{clspqq}.

Suppose that $H_{\phi}$ is simple. Then $H_{\phi} \simeq H$ as
algebras and as coalgebras, whence $V \# kG(H)$ is a simple
subcoalgebra of $H_{\phi}$. The claim can be established using the
argument in the proof of  \cite[Theorem 1.0.1]{clspqq}. We sketch
 this proof below for the sake of completeness.

Write $G(V) = \{ x_1, \dots, x_4 \}$. By \cite[Lemma
2.2.1]{clspqq}  $(\mathcal S_{R_{\phi}}V) V$ is a Yetter-Drinfeld
subcoalgebra of  $R_{\phi}$,  and   $(\mathcal S_{R_{\phi}} V)x_i$
is a left coideal of $R$ containing 1, for all $i  = 1, \dots, 4$,
where $\mathcal S_{R_{\phi}}$ is the antipode of $R_{\phi}$. Since
the dimensions of the irreducible comodules of a twisted dual
group algebra of $G$ divide the order of $G$, we have that the
dimensions of the irreducible $R_{\phi}$-comodules are either $1$
or  $2$.

The group $G$ acts on $R_{\phi}$ by $\leftharpoonup : R_{\phi}
\otimes kG \to R_{\phi}$, $r \leftharpoonup g : = (g, r_{-1})
r_0$, $r \in R_{\phi}$, where $(\ , \ )$ is a non-degenerate
bicharacter on $G$. Moreover, this action restricts to an action
of $G$ on $(\mathcal S_{R_{\phi}} V) V$ by coalgebra
automorphisms. Let $\omega_1, \dots, \omega_l$ be  the group-like
elements of $R_{\phi}$ which belong to $(\mathcal S_{R_{\phi}} V)
V$. Then the action  $\leftharpoonup$ of $G$ permutes the set $\{
\omega_1, \dots, \omega_l \}$. As in \cite[Claim 5.3.1]{clspqq},
one can see that every group-like element $\omega$  in $(\mathcal
S_{R_{\phi}} V) V$ belongs to $(\mathcal S_{R_{\phi}} V) x_i$ for
all $1 \leq i \leq n$.

Therefore $\mathcal S_{R_{\phi}} V x_i  = k\omega_1 \oplus \dots
\oplus k\omega_l \oplus \bigoplus_j W_j$,  for all $i$,  where
$W_j$ are irreducible right coideals of $R_{\phi}$ of dimension
bigger than $1$.  Thus $l$ is even, and since $1 \in \{ \omega_1,
\dots \}$, there exists a nontrivial  group-like element in
$R_{\phi}$ invariant under the action $\leftharpoonup$ of $G$;
this is exactly the same as a coinvariant group-like element in
$R_{\phi}$. Then $G(H_{\phi}) \cap R_{\phi}$ has an element of
order $3$. This is a  contradiction, because $G(H_{\phi})$ is of
order 4. The contradiction comes from the assumption that
$H_{\phi}$ is simple, hence the claim follows. \end{proof}

Let $A \subseteq H_{\phi}$ be a proper normal Hopf subalgebra.
Since $H_{\phi} \simeq H$ as algebras, $H_{\phi}$  contains no
Hopf subalgebras of dimension $12$, by Remark \ref{rmkk-36} (iv).
Then $\dim A = 2, 3, 4, 6$ or $18$. Note that $kG(H)_{\phi} =
kG(H)$ is a Hopf subalgebra of $H_{\phi}$;  so that 4 divides the
order of $G(H_{\phi})$. If $\dim A = 3$, then either $\vert
G(H_{\phi})\vert = 12$, which contradicts Remark \ref{rmkk-36}
(iv),  or else $H_{\phi}$ is cocommutative. The last possibility
implies $H \simeq \mathcal H$ by \cite[4.3]{twist-simple}; hence
we may assume that $\dim A \neq 3$.

If $\dim A = 2$, then $A \subseteq kG(H_{\phi}) \cap Z(H_{\phi})$
by Corollary \ref{kob-mas}. Thus $A \subseteq kG(H) \cap Z(H)$;
see Corollary \ref{andrus}. Also, if $\dim A = 4$, then there is a
quotient Hopf algebra $H_{\phi} \to H_{\phi} / H_{\phi} A^+ = :
\overline{H_{\phi}}$, such that $\dim \overline{H_{\phi}} = 9$. In
particular, $9 / |G(H_{\phi}^*)|$, which is a contradiction.

Suppose that $\dim A = 6$. Then, with the above notation, $\dim \overline{H_{\phi}} = 6$. Hence $\overline{H_{\phi}}^*$ contains a group-like Hopf subalgebra of dimension $3$ or else a simple subcoalgebra of dimension $4$, which again contradicts the fact that $H_{\phi}$ is of type $(1, 4; 4, 2)$ as an algebra.

Finally, if $\dim A = 18$, then $|G(A)| = 2$ and by \cite{masuoka-further} $A$ is commutative. It follows from \cite[Corollary 3.9]{harmonic}, $\dim V \leq [H_{\phi}: A] = 2$, for all irreducible $H_{\phi}$-module $V$. This is impossible since $H_{\phi}$ is of type $(1, 4; 4, 2)$ as an algebra.
This finishes the proof of the theorem. \end{proof}

\chapter{Dimension $40$}\label{40}

\section{Reduction of the problem} Let $H$ be a nontrivial semisimple Hopf algebra of dimension $40$ over $k$.

\begin{lemma}\label{conteo40} The order of $G(H^*)$ is either
$4$, $8$ or $20$ and as an algebra $H$ is of one of the following types:
\begin{itemize}\item $(1, 4; 2, 9)$, \item $(1, 4; 2, 5; 4, 1)$, \item $(1, 4; 2, 1; 4, 2)$, \item $(1, 8; 2, 8)$, \item $(1, 8; 2, 4; 4, 1)$, \item $(1, 8; 4, 2)$, \item $(1, 20; 2, 5)$. \end{itemize} \end{lemma}

We shall prove that the type $(1, 4; 2, 5; 4, 1)$ is impossible.
See Lemma \ref{vs-4}.

\begin{proof} Let $n = |G(H^*)|$. It follows from \ref{alg-struct} that $n \neq 1, 5, 10$. The possibility $n = 2$ is discarded using Theorem \ref{thm-nr}. In the cases $n = 4, 8, 20$ the only possibilities excluded are the types $(1, 4; 3, 4)$ and $(1, 4; 6, 1)$; however, in this cases $H$ has an irreducible character $\psi$ of degree $3$ (respectively $6$), and decomposing the product $\psi \psi^*$ as  a sum of irreducible characters gives a contradiction. This finishes the proof of the lemma. \end{proof}

\begin{remark} (i) If $H$ is of type $(1, 20; 2, 5)$, then $H$ is not simple, by Corollary \ref{kob-mas}.

(ii) Let $H$ be as in Lemma \ref{conteo40}. It follows from Proposition \ref{cociente8} that there is a Hopf subalgebra $A \subseteq H$, with $\dim A = 8$.  \end{remark}

\begin{lemma}\label{red-40} We have $H = R \# A$ is a biproduct where $A$ is a semisimple Hopf algebra of dimension $8$ and $R$ is a $5$-dimensional Yetter-Drinfeld Hopf algebra over $A$. \end{lemma}

\begin{proof} There are a Hopf subalgebra $A \subseteq H$ and a quotient Hopf algebra $q: H \to B$, such that $\dim A = \dim B = 8$. In particular, $\dim H^{\co B} = 5$, and therefore $G(H) \cap H^{\co B} = 1$ by \cite{NZ}.

Consider the restriction $q\vert_A: A \to B$. It will be enough to show that $q$ is an isomorphism, or equivalently, that $A \cap H^{\co B} = k1$.
Suppose on the contrary that $A \cap H^{\co B} = A^{\co B}$ is of dimension $2$ or $4$; see Lemma \ref{restriction}. Then there must exist $1 \neq g \in G(A) \cap A^{\co B}$, since the irreducible left coideals of $A$ are of dimension $1$ or $2$. Then $2 / \dim H^{\co B} = 5$, which is absurd. This contradiction finishes the proof of the lemma. \end{proof}

\begin{lemma}\label{r=k1+u} Suppose that $H$ is simple and $|G(H)| = 4$. Then $R = k1 \oplus U$, where $U$ is an irreducible left coideal of $H$ of dimension $4$. In particular, $H$ is not of type $(1, 4; 2, 9)$ as a coalgebra. \end{lemma}

\begin{proof} Since $|G(H)| = 4$, $A$ is not cocommutative.
Suppose on the contrary that $R = k1 \oplus V_1 \oplus V_2$, where
$V_i$ is an irreducible left coideal of  dimension $2$. By  Lemma
\ref{bip-2-dim} (i), $\rho (V_i) \subseteq kG(H) \otimes V_i$, $i
= 1, 2$. Thus $\rho (R) \subseteq kG(H) \otimes R$ and $R \#
kG(H)$ is a  normal Hopf subalgebra of $H$ (of index 2).
\end{proof}

\begin{lemma}\label{vs-4} There exists no semisimple Hopf algebra $H$ of type  $(1, 4; 2, 5; 4, 1)$ as a coalgebra. \end{lemma}

\begin{proof} Suppose on the contrary that $H$ is of type $(1, 4; 2, 5; 4, 1)$.
Let $A \subseteq H$ be a Hopf subalgebra of dimension $8$.  For
all irreducible $H^*$-characters $\chi$ such that $\deg \chi = 2$,
we have that $G[\chi ] \neq 1$, because $H$ does not contain
irreducible coideals of dimension $3$. Also, the tensor product
$\chi \chi'$ is irreducible for some irreducible characters $\chi$
and $\chi'$ of degree $2$, since otherwise, there would be a Hopf
subalgebra of dimension $4 + 4 . 5 = 24$, which is impossible.
Therefore there exist $\chi, \psi \in X_2$ such that $G[\chi]$ is
not contained in $G[\psi]$; c.f. Lemma \ref{productodesimples}.
Thus $G[\chi]$ and  $G[\psi]$ are distinct subgroups of order $2$
of $G(H)$.

Since $|G[\chi]| = |G[\psi]| = 2$, then the orbits $G(H) \chi$ and
$G(H) \psi$ are of order $2$. Moreover, since $G[\chi'] =
G[\chi]$, for all $\chi' \in G(H)\chi$ (because $G[\chi]$ is
normal in $G(H)$), then $G(H) \chi$ and $G(H) \psi$ are disjoint
orbits. Let $\tau \notin A$ be an irreducible character of degree
$2$. Then we have
\begin{equation}\label{relacion}\tau \tau^* = \sum_{a \in G[\tau]}
a + \lambda,\end{equation} for some $\lambda \in X_2$. By
\eqref{relacion}, we have $\lambda^* = \lambda$ and $G[\tau]
\subseteq G[\lambda]$. Therefore, either $\lambda \in G(H) \tau$
or $\lambda \in A$.

\begin{claim}\label{af-la} We have $\lambda \in A$. \end{claim}

\begin{proof} Suppose on the contrary that $\lambda \in G(H) \tau$; hence, every orbit has a self-dual element.
Moreover, if $g \in G(H)$ is such that $g \lambda = \tau$, then $\tau^* = \lambda g$ and therefore $G[\tau] = G[\lambda] \subseteq G[\tau^*]$.
Hence $\tau^* \in G(H) \tau$. This implies, since each orbit has two elements, that $\tau^* = \tau$.
As before, we have decompositions
$\tau^2 = \tau \tau^* = \sum_{a \in G[\tau]} a + \lambda$, where $\lambda \in G(H) \tau$ and if $\tau' \in G(H) \tau$, say $b \tau = \tau'$, then
 $\tau' \tau =  b \tau^2 = b + ba + b \lambda$, and $b \lambda  \in G(H) \tau$.
Then the set $G(H) \cup G(H)\tau$ spans a standard subalgebra of
$R(H^*)$ which corresponds to a Hopf subalgebra of dimension $12$.
This contradicts \cite{NZ}; so we may conclude that $\lambda$ is
not in the orbit $G(H)\tau$. Therefore, we must have $\lambda \in
A$ as claimed. \end{proof}

Since $\chi \psi$ is irreducible, it follows from Lemma
\ref{productodesimples} that $m(\lambda, \psi^*\psi) = 0$. This
contradicts Claim \ref{af-la}. Then the lemma follows.
\end{proof}

\begin{remark}\label{rmk-40} Keep the notation in Lemma \ref{red-40}.

(i) Suppose that $R = k1 \oplus U$, where $U$ is an irreducible left coideal of dimension $4$ and let $\chi \in H$ be the corresponding character.
It follows from Corollary \ref{formula}, that $|G[\chi^*]| \leq 4$.

(ii) If $|G(H^*)| = |G(H)| = 8$, then $H$ is not simple by \cite[Lemma 2.2.5]{pqq2}. So we may assume that $A$ is not a group algebra of an abelian group.
In particular, after dualizing if necessary, we may assume that $A$ is not cocommutative.  \end{remark}

\begin{lemma} Suppose that $H$ is of type $(1, 8; 2, 4; 4, 1)$ as a coalgebra. Then $H$ is not simple. \end{lemma}

\begin{proof} By Remark \ref{rmk-40} (ii), we may assume that $G(H)$ is not abelian of order $8$. Consider the action by right multiplication of $G(H)$  on the set $X_2$ of irreducible characters of degree $2$. Then either this action is transitive or else $|G[\chi^*]| = 4$, for all $\chi \in X_2$.

The first possibility implies, in view of Remark \ref{nr} (iii)
that $H$ has a Hopf subalgebra of dimension $24$, which is
impossible. Consider the second possibility. Since $G(H)$ is not
abelian, $1 \neq [G; G] \subseteq G[\chi'] \cap G[\chi^*]$, for
all $\chi', \chi \in X_2$. Therefore, by Theorem \ref{corolario},
also in this case there is a Hopf subalgebra of dimension $24$,
which is impossible. This proves the lemma. \end{proof}

\begin{lemma} Suppose that $H$ is of type $(1, 8; 4, 2)$ as a coalgebra. Then $H$ is not simple. \end{lemma}

\begin{proof} Suppose that $H$ is simple. We have $\dim R(H^*) = 10$ and $G(H)$ is not abelian, by Remark \ref{rmk-40} (ii).
Let $e_0, e_1, e_2, e_3, f_1, f_2$ be orthogonal primitive
idempotents in $kG(H)$ such that $\dim kG(H) e_j = 1$ and $\dim
kG(H) f_j = 2$. Then $\dim He_j = 5$ and $\dim Hf_j = 10$, and
moreover $Hf_1 \simeq Hf_2$.

\smallbreak We may assume that $|G(H^*)| \neq 8$. Then, by Lemma
\ref{r=k1+u}, $R^* \simeq k1
\oplus U$, where $U$ is an irreducible left coideal of $H^*$,  and
{\it a fortiori} an irreducible Yetter-Drinfeld submodule of
$H^*$.

By Proposition \ref{ref-corr-z2}, there is a bijective correspondence between the irreducible Yetter-Drinfeld submodules of $H$ appearing in $R^* = (H^*)^{\co kG(H)^*}$ and primitive idempotents $e \in R(H^*)$ such that $ee_0 \neq 0$.

Therefore we have a decomposition $e_0 = \Lambda + E_0$, where
$\Lambda \in H$ is the normalized integral and $E_0$ is a
primitive idempotent of $R(H^*)$ such that $\dim HE_0 = 4$. This
gives at least $7$ orthogonal idempotents in $R(H^*)$, implying
that $R(H^*) \simeq k^{(6)} \times M_2(k)$ as an algebra. In
particular, since $f_1$ and $f_2$ are not central in $kG(H)$, they
correspond necessarily to the $4$-dimensional simple component  of
$R(H^*)$.

We find that there must exist some $1 \leq j \leq 3$ such that $e_j$ is not primitive in $R(H^*)$. By the Kac-Zhu Theorem, this implies the existence of a primitive idempotent $\Lambda \neq E \in R(H^*)$ such that $\dim HE = 1$. Then $G(H^*) \cap Z(H^*) \neq 1$ and the lemma follows. \end{proof}

\begin{lemma} Suppose that $H$ is of type $(1, 8; 2, 8)$ as a coalgebra. Then $H$ is not simple. \end{lemma}

\begin{proof} We have $H = R \# kG$, where $G = G(H)$ and
$R = k1 \oplus V_1 \oplus V_2$, where $V_i$ are irreducible left
coideals of $H$ such that $\dim V_i = 2$. By Theorem \ref{thm-nr}
we must have $G[\chi] \neq 1$, for all irreducible characters
$\chi$ of degree $2$. Then $V_i$ is a subcoalgebra of $R$, by
Proposition \ref{R-coideal}, and $V_1$ is not isomorphic to $V_2$
as left coideals of $H$. The group $G$ permutes the set $G(V_1)
\cup G(V_2)$ and we may also assume that the group homomorphism $G
\to \mathbb S_4$ is injective. In particular, $G \simeq D_4$ and
acts transitively on $G(V_1) \cup G(V_2)$.

Let $G_i = \{ g: g.V_i \subseteq V_i \}$. Then $|G_i| \geq 4$ and
in particular, $Z(G) \subseteq G_i$, $i = 1, 2$. By Lemma
\ref{tr-act}, if $Z(G)$ acts transitively on $G(V_i)$, then
$|G_{V_i}| \leq 2$.  Moreover, the action of $Z(G)$ is either
transitive or trivial on $G(V_1)$. Observe in addition that $V_1$
generates $R$ as an algebra. Hence, the last possibility implies
that $Z(G) \subseteq Z(H)$ and $H$ is not simple.

In the other case, we may conclude that $\rho (R)  \subseteq kG_{V_i} \otimes R$ and since $|G_{V_i}| = 2$, $B = R \# kG_{V_i}$ is a Hopf subalgebra of dimension $10$.
By Remark \ref{rmk-40} (ii), we may assume that $|G(H^*)| = 4$.
Then $(H^*)^{\co B^*} = k1 \oplus kt \oplus V$, where $1 \neq t \in G(H)$ and $V$ is an irreducible left coideal of dimension $2$. Consider the possible decompositions of $(H^*)^{\co B^*}$ as a Yetter-Drinfeld submodule of $H^*$, as in Lemma \ref{YD-coinv}. Since $3$ does not divide $\dim H$, the dimensions of the irreducible summands are either $1$ or $2$. Decomposing each irreducible Yetter-Drinfeld summand as a sum of irreducible left coideals, implies that $kt$ is necessarily a Yetter-Drinfeld submodule of $H$. Hence $t \in G(H) \cap Z(H)$ and $H$ is not simple.
This concludes the proof of the lemma. \end{proof}

\section{Main result} In view of the lemmas in the previous section, we may assume that $H$ is of type $(1, 4; 2, 1; 4, 2)$ as an algebra and as a coalgebra. Moreover we have $H = R\# A$, where $A \simeq H_8$  and $R = k1 \oplus U$, for some irreducible left coideal  $U$ of dimension $4$. We keep this notation along this section.

\begin{lemma}\label{cocom-40} If $R$ is cocommutative then $H$ is not simple. \end{lemma}

\begin{proof} We have $U \subseteq R$ is a subcoalgebra in ${}_A^A\mathcal{YD}$ and $\mathcal S_RU = U$.
On the other hand, $Ut \simeq U$ for some $1 \neq t \in G(H) =
G(A)$. Let $G(U) = \{ x_1, \dots, x_4 \}$. By \cite[Lemmas 2.2.1
and 2.2.3]{clspqq} $Ux_i$ is a left coideal of $R$ and $Ux_i
\simeq U(t.x_i)$; so if $t.x_i \neq x_i$, by counting dimensions,
$x_i \in Ux_i$ implying that $1 \in U$, which is a contradiction.
Therefore  $t$ acts trivially on $R$. The lemma will follow from
the following claim:

\begin{claim} Let $A \to \End V$ be an action of $A \simeq H_8$ on a vector space $V$ such that there exists $1 \neq t \in G(A)$ with $t\vert_V = \id_V$. Then $t'\vert_V = \id_V$ for all $t' \in G(A) \cap Z(A)$.
\end{claim}

\begin{proof}
By \cite{ma-6-8} we know that $A$ is generated as an algebra by
$G(A) = \{ 1, x, y, xy \}$ and an invertible element $z$ such that
$zx = yz$. Moreover, $G(A) \cap Z(A) = \{ 1, xy \}$. If $t = xy$,
then there is nothing to prove. Thus, without loss of generality
we may assume that $t = x$. In this case the relation $zx = yz$
plus the fact that $z$ is invertible imply that also $y\vert_V =
\id_V$. Then we have $xy\vert_V = \id_V$ as claimed. \end{proof}
\end{proof}

We now prove the main result of this chapter.

\begin{theorem} Let $H$ be a semisimple Hopf algebra of dimension $40$. Then $H$ is not simple.
\end{theorem}

\begin{proof} We may assume that $H = R \# A$ is a biproduct, where $A \simeq H_8$ and $R = k1 \oplus U$.
By Lemma \ref{cocom-40} it will be enough to show that $U$ is
cocommutative. Suppose on the contrary that $U$ is a simple
coalgebra. Observe that $U \# A$ is a subcoalgebra of $H$;
moreover,  there exist $t, s \in G(H)$ such that $Ut \simeq U$ and
$Us$ is an irreducible coideal not isomorphic to $U$. In
particular, $(U \# A) \cap C \neq 0$ and $(U \# A) \cap C' \neq
0$, where $C$ and $C'$ are the simple subcoalgebras of $H$
containing $U$ and $Us$, respectively. Thus, $U \# A = C \oplus C'
\simeq M_4(k) \oplus M_4(k)$.

On the other hand, the assumption that $U$ is simple implies that
$U \# A \simeq U \otimes {}_{\phi} A$ for some invertible
normalized $2$-cocycle $\phi \in A \otimes A$; see
\cite[7.3.1]{Mo}. By \cite[Theorem 4.8 (1)]{masuoka-cont}
${}_{\phi} A \simeq A$ as coalgebras, which gives $U \# A \simeq
M_2(k)^{(4)} \oplus M_4(k)$. This is a contradiction, which
implies that $U$ is cocommutative and finishes the proof of the
theorem. \end{proof}

\chapter{Dimension $42$}\label{42}

Let $H$ be a semisimple Hopf algebra of dimension $42$ over $k$. We
shall assume that $H$ is nontrivial.
Suppose that  $H$ fits into an abelian extension
\begin{equation*} 1 \to k^{\Gamma} \to H \to kF \to 1,\end{equation*}
where $\Gamma$ and $F$ are finite groups. It follows from the
results in  \cite[Section 4]{pqq} that in this case $H$ is
isomorphic to one of the Hopf algebras $\mathcal A_7(2, 3)$ or
$\mathcal A_7(3, 2) = \mathcal A_7(2, 3)^*$ constructed in
\cite[(2.3.1)]{examples}; indeed $\mathcal A_7(2, 3)$
(respectively, $\mathcal A_7(3, 2)$) fits into an extension as
above, where $\Gamma$ is the non-abelian group of order $14$ and
$F$ is  the cyclic group of order $3$, (respectively, $\Gamma$ is
the non-abelian group of order $21$ and $F$ is  the cyclic group
of order $2$).

For these Hopf algebras, $G(\mathcal A_7(3, 2)) \simeq G(\mathcal
A_7(2, 3)) \simeq \mathbb Z_6$. As coalgebras, $\mathcal A_7(3,
2)$ is of type $(1, 6; 2, 9)$ and $\mathcal A_7(2, 3)$ is of type
$(1, 6; 3, 4)$.

\section{Possible (co)-algebra structures} Let $H$ be a nontrivial semisimple Hopf algebra of dimension 42.

\begin{lemma}\label{conteo42} The order of $G(H^*)$ is either $2$, $6$ or $14$.
As an algebra $H$ is of one of the following types:
\begin{itemize} \item $(1, 2; 2, 1; 6, 1)$, \item $(1, 2; 2, 1; 3, 4)$, \item  $(1, 2; 2, 10)$, \item $(1, 6; 6, 1)$, \item $(1, 6; 2, 9)$, \item $(1, 6; 3, 4)$, \item $(1, 14; 2, 7)$. \end{itemize} \end{lemma}

\begin{proof} Let $n = |G(H^*)|$. A counting argument, using \ref{alg-struct},
shows that $n \neq 7, 21$, and the only possibilities are the
prescribed ones if $n = 6$ or $14$. If $n = 1$  we find that $H$
must have an irreducible module of degree $2$, which  contradicts
Corollary \ref{cor-g-1}. If  $n = 3$, then  $H$ is necessarily of
type $(1, 3; 2, 3; 3, 3)$, which is discarded by Theorem
\ref{thm-nr}.

Suppose finally that $n = 2$. Then, either $H$ is of one of the
prescribed algebra types, or else $H$ is of one of the following
types: $$(1, 2; 2, 6; 4, 1), \quad  \quad (1, 2; 2, 2; 4, 2).$$
However, by Remark \ref{el-nr} (i), in these cases we have $G[\chi] = G(H^*)$, for all irreducible character $\chi$ of degree 2. By Theorem \ref{corolario}, $H$ has a Hopf algebra
quotient $H \to \overline H$, where $\dim \overline H = 2 + 6.4$
or $\dim \overline H = 2 + 2.4$, respectively. By \cite{NZ}, this is not
possible. Then the lemma follows. \end{proof}

\begin{remark}\label{obs} (i) Since $12$ does not divide the dimension of $H$, by Theorem \ref{thm-nr}, every irreducible character of degree $2$ is stable under left multiplication by some group-like element $g \in G(H^*)$ of order $2$.

(ii) It follows from Lemma \ref{conteo42} that $G[\chi]$ is nontrivial for all
irreducible character $\chi$ of degree  $6$. Also, if $|G(H^*)| = 6$, then for
all irreducible character $\chi$ of degree $3$, $G[\chi]$ is of order $3$.

(iii) It follows from Lemma \ref{conteo42} that $G(H^*)$ always contains a subgroup of order $2$. \end{remark}

\begin{lemma}\label{2} (i) Suppose $H$ is of type $(1, 2; 2, 10)$ as a coalgebra. Then $H$ is commutative.

(ii) Suppose  $H$ is of type $(1, 2; 2, 1; 6, 1)$ or $(1, 2; 2, 1;
3, 4)$ as a coalgebra. Then  $H$ contains a  Hopf subalgebra
isomorphic to $k^{\mathbb S_3}$, where $\mathbb S_3$ is the symmetric group on
$3$ symbols. \end{lemma}

\begin{proof} (i) This follows from Theorem \ref{1-2-2-n}, since by Remark \ref{obs} (iii), $|G(H^*)|$ is even.

(ii) In this case $H$ has a unique simple subcoalgebra $C$ of dimension $4$, necessarily stable under left and  right multiplication by $G(H)$. Therefore $A = kG(H) \oplus C$ is a Hopf subalgebra of $H$ of dimension $6$, which is not
cocommutative;
hence, in view of  the classification of semisimple Hopf algebras of dimension $6$,  $A \simeq k^{\mathbb S_3}$. \end{proof}

\begin{lemma}\label{14} Suppose that $H$ is of type $(1, 14; 2, 7)$ as a coalgebra. Then $H$ is commutative.  \end{lemma}

\begin{proof} Note that the group $G(H^*)$ must be abelian by Proposition \ref{nab-pq}.
Suppose that $G(H^*)$ is of order $2$ or $6$. By Lemmas
\ref{conteo42} and \ref{2}, either there is a quotient Hopf
algebra $p: H \to kG$, where $G$ is a group of order $6$, or else
$G(H^*) \simeq \mathbb S_3$ is nonabelian of order 6. In the first
case, a dimension argument implies that $H^{\co p} = k\Gamma$,
where $\Gamma \subseteq G(H)$ is the only subgroup of order $7$.
Therefore, $H$ fits into the (abelian) extension $1 \to k\Gamma
\to H \to kG \to 1$, and the lemma follows from \cite[Section
4]{pqq}. In the second case, there is a projection $q: H \to kT
\simeq k^T$, where $T \subseteq G(H^*)$ is the only subgroup of
order 3; hence $\dim H^{\co q} = 14$, and $\Gamma \subseteq H^{\co
q}$. A dimension argument shows that $G(H) \subseteq H^{\co q}$
and thus $kG(H) = H^{\co q}$. Therefore $H$ fits into an abelian
exact sequence $1 \to kG(H) \to H \to kT \to 1$, and the lemma
follows also from \cite[Section 4]{pqq}.

We claim  that $H^*$ is not of type $(1, 14; 2, 7)$ as a
coalgebra. Indeed, in this case, we have $G(H) \simeq G(H^*)$ is
abelian  and $H \simeq R \# kG(H)$ is a biproduct, where $R$ is of
dimension $3$. In particular, $R$ is commutative and
cocommutative, and it follows easily that the only subgroup $F$ of
order $7$ of $H$ acts trivially on $R$. Then $F$ is central in
$H$, and there is a central extension $0 \to kF \to H \to K \to
1$. Now, since every $6$-dimensional semisimple Hopf algebra is
trivial and since $|G(H^*)|$ is not divisible by $3$, we find that
$K$ is necessarily isomorphic to a  group algebra. In particular,
this extension is abelian. This is a contradiction since $|G(H)|
\neq 6$. This completes the proof of the lemma. \end{proof}

\begin{lemma}\label{61} $H$ is not of type $(1, 6; 6, 1)$ as a coalgebra.
\end{lemma}

\begin{proof} Suppose on the contrary that $H$ is of type $(1, 6; 6, 1)$ as a
coalgebra. It follows from Lemmas \ref{conteo42}, \ref{2} and
\ref{14} that there is a quotient Hopf algebra $H \to A$, where
$A$ is of dimension $6$. Hence, $A \simeq kG(H)$ and $H$ is a
biproduct $H \simeq R \# kG$, where $R \simeq H/H(kG)^+$ is a
braided Hopf algebra over $G$ of dimension $7$, and $G = G(H)$.

As a left coideal of $H$, $R$ decomposes as a direct sum $R = k1 \oplus
V$, where $V$ is an irreducible left coideal of dimension $6$. In particular,
$V$ is a subcoalgebra of $R$ in the category of Yetter-Drinfeld modules over
$G$, and the smash coproduct coalgebra $V \# kG$ coincides with the unique
simple subcoalgebra of $H$ of dimension $36$.
Since $H^2(G, k^{\times}) = 1$, it follows from \cite[1.3.1]{pqq2} that $V$ is
cocommutative and the action of $G$ permutes the $6$ distinct group-like
elements. Thus $V$ has a basis $x_i$, $1 \leq i \leq 6$,
consisting of group-like elements of $R$.

By \cite[2.2]{clspqq}, $Vx_i$ is a left coideal of $R$ containing
$1$, and   $V x_i \simeq Vx_j$, for all $1 \leq i, j \leq 6$ such
that $gx_i = x_j$, for some $g \in G$. Note that $x_i$ does not
belong to $Vx_i$: indeed, since $x_i$ is invertible (with inverse
$x_i^{-1} = \mathcal S_R(x_i)$), if $x_i = \sum_jx_jx_i$, then $1
\in \sum_jx_j \in V$, which is absurd.

Decomposing $Vx_i$ into a direct sum of irreducible left coideals
of $R$, we get $Vx_i = k1 \oplus \bigoplus_{l \neq i} kx_l$, where
$g_i \neq g$, for all $i$. If $g \in G$ is such that $gx_i = x_j
\neq x_i$, we have $gx_i \in Vx_i$. But then $x_j \in Vx_j \simeq
Vx_i$. This is a contradiction. Then the lemma follows.
\end{proof}

\begin{lemma}\label{gl-z2} Suppose  $H$ is of coalgebra type $(1, 2; 2, 1; 6, 1)$ or $(1, 2; 2, 1; 3, 4)$. Then $H$ is commutative. \end{lemma}

\begin{proof} We shall show that $H$ fits into an abelian extension, and hence is trivial by \cite[Section 4]{pqq}.

By Lemma \ref{2}, $H$ has  a  Hopf subalgebra $A$ isomorphic to
$k^{\mathbb S_3}$. We claim that there is no Hopf algebra
surjection $p: H \to k\mathbb S_3$. If this were the case, then $A
\cap H^{\co p} = k1$. Indeed, as a left coideal of $H$, $H^{\co
p}$ decomposes as a direct sum or irreducible left coideals. Note
that $G(A) = G(H)$ intersects trivially $H^{\co p}$ by \cite{NZ}.
If $A \cap H^{\co p} \neq k1$, then there is a $2$-dimensional
irreducible left coideal $U$ of $A$ such that $U \subseteq H^{\co
p}$. Then $H^{\co p} = k1 \oplus U \oplus W$, where $W$ is a left
coideal of $H$ of dimension $4$, such that $W$ contains no
one-dimensional left coideal of $H$. This is not possible and
therefore $A \cap H^{\co p} = k1$ as claimed. But then the
restriction $p : k^{\mathbb S_3} \to k\mathbb S_3$ is an
isomorphism, which is not possible. This establishes the claim. In
particular, $|G(H^*)| \neq 2$.

We may then assume that $|G(H^*)| = 6$. Hence, by Lemma \ref{61},
$H^*$ is of type $(1, 6; 3, 4)$ or  $(1, 6; 2, 9)$ as a coalgebra.
Dualizing the inclusion $A \subseteq H$, we get a Hopf algebra
quotient $p: H^* \to k{\mathbb S_3}$. By \cite{NZ}, $G(H^*) \cap
H^{\co p} = k1$; so that the restriction $p : kG(H^*) \to
k{\mathbb S_3}$ is an isomorphism. In particular, $G(H^*)$ is
non-abelian and $H^*$ is isomorphic to a biproduct $H^* \simeq R
\# kG(H^*)$, where  $\dim R = 7$. Let $\Gamma \subseteq G(H^*)$ be
the unique subgroup of order $3$.

\smallbreak \noindent {\bf Case I.} $H^*$ is of type $(1, 6; 3,
4)$.

In this case, $R = k1 \oplus W_1 \oplus W_2$, where $W_i$ are
irreducible left coideals of $H^*$ of dimension $3$. Since $\Gamma
= G[\chi]$ for all irreducible $\chi$ of degree $3$, then $gW_i
\simeq W_i \simeq W_ig$, for all $g \in \Gamma$, $i = 1,2$.

\begin{claim} $\rho(W_i)$ is not contained in $k\Gamma \otimes W_i$. \end{claim}

\begin{proof} Suppose on the contrary that $\rho(W_1) \subseteq k\Gamma \otimes W_1$. Let $\widetilde R : = k[W_1]$; then $\widetilde R$ is a subalgebra and subcoalgebra of $R$ which is stable under the action of $\Gamma$ and $\rho (\widetilde R) \subseteq k\Gamma \otimes \widetilde R$.
Therefore $\widetilde A = \widetilde R \# k\Gamma$ is a Hopf
subalgebra of $H^*$, and $\dim \widetilde A = 3 \dim \widetilde R
> 12$. Then $\dim \widetilde A = 21$, thus $\widetilde A$ is
commutative and normal, and moreover the quotient $H^*/H^*\widetilde A^+$ is
cocommutative, since $k\Gamma \subseteq \widetilde A$. Hence, $G(H^*)$ is abelian, which is a contradiction. This
proves the claim.  \end{proof}

Fix $i = 1, 2$. We may assume that there exists $0 \neq w \in
W_i$,  such that $\rho(w) = a \otimes w$, where $a \in G(H^*)$ is
an element of order $2$. Let $1 \neq g \in \Gamma$, so that $0
\neq g . w \in W_i$ is homogeneous of degree $gag^{-1} \neq a$.
Hence $W_i^{\co G(H^*)} = 0$, and $R^{\co G(H^*)} = k1$.

Dualizing, $H = R^* \# A$ and the subalgebra $(R^*)^A$ of
invariants in $R^*$ is one-dimensional,  since $(R^*)^A \simeq
\Hom^{G(H^*)}(R, k)$. By \cite[Section 8]{harmonic}, $(R^*)^A$ is
isomorphic to the {\it Hecke algebra} of the pair $A, H$;  that
is, $(R^*)^A \simeq e_AHe_A$,  where $e_A \in A$ is the normalized
integral. By \cite[Theorem 4.4]{harmonic} there is a bijection
between  irreducible representations of $(R^*)^A$ and irreducible
representations of $H$ appearing in $\Ind_A^H\epsilon_A$ with
positive multiplicity. This is a contradiction, hence the lemma
follows.

\smallbreak \noindent {\bf Case II.} $H^*$ is of type $(1, 6; 2,
9)$.

We shall use the action defined in \eqref{accion}.

\begin{claim}\label{claim42-II} The action $\Gamma \times \Gamma \times X_2 \to X_2$, given by $(g, h) . \chi = g \chi h^{-1}$ is freely transitive. \end{claim}

\begin{proof} Let $\chi \in X_2$.
It  will be enough to show that the stabilizer $(\Gamma \times
\Gamma)_{\chi}$ is trivial. Let $g, h \in \Gamma$ such that $g
\chi h^{-1} = \chi$. Then $G[\chi] = G[g \chi h^{-1}] = g G[\chi]
g^{-1}$. This implies that $g = 1$, because $G[\chi]$ is of order
2 and $G(H^*)$ is not abelian. Therefore $\chi h^{-1} = \chi$, and
$h = 1$. This proves the claim.
\end{proof}

Observe that there exists $\chi \in X_2$ such that $\chi^* =
\chi$, because $|X_2|$ is odd. Hence $\chi^2 = 1 + a + \chi'$,
where $(\chi')^* = \chi'$, and $G[\chi] = \{ 1, a \} = G[\chi']$.
Note also that $\chi \neq \chi'$, since $H$ has no quotient of
dimension $6$ which is not commutative. Write $\chi' = g \chi
h^{-1}$, $g, h \in \Gamma$. Hence $G[\chi] = G[\chi'] = g G[\chi]
g^{-1}$, implying $g = 1$. So $\chi' = \chi h^{-1}$ and $h \neq
1$. Thus $\chi  h^{-1} = \chi' = (\chi')^* = h \chi^* = h \chi$;
which implies that $h \chi h = \chi$. This contradicts Claim
\ref{claim42-II}. The proof of the lemma is now complete.
\end{proof}

\section{Classification} We shall now prove  Theorem \ref{cls42}. In view of the previous results, we may assume that $H$ is  of type $(1, 6; 2, 9)$
or $(1, 6; 3, 4)$ as a coalgebra. The proof of the theorem of will follow from the
next two lemmas.

\begin{lemma} Suppose that $H$ is of type $(1, 6; 2, 9)$ as a coalgebra. Then
$H$ is isomorphic to $\mathcal A_7(3, 2)$. \end{lemma}

\begin{proof} By the above, we also have $|G(H^*)| = 6$. Therefore the groups
$G(H)$ and $G(H^*)$ are both abelian. This  implies that $G(H)$
contains a unique subgroup $F$ of order $2$. In particular, $H
\simeq R \# kF$ is a biproduct over $F$. The subgroup $F$
necessarily stabilizes all $4$-dimensional simple subcoalgebras.
Therefore the quotient coalgebra $H/H(kF)^+ \simeq R$ is
cocommutative. It follows from  Proposition \ref{som} that $H$
fits into an abelian extension. This implies the lemma.
\end{proof}

\begin{lemma} Suppose that $H$ is of type $(1, 6; 3, 4)$ as a coalgebra. Then
$H$ is isomorphic to $\mathcal A_7(2, 3)$. \end{lemma}

\begin{proof} Also here we have $H \simeq R \# kG$ is a biproduct, where
$G$ is the unique subgroup of order $3$ of $G(H)$. By Remark \ref{obs}, $G$ stabilizes all simple subcoalgebras of $H$ of dimension $9$. Then $R$ is a cocommutative coalgebra and the lemma follows from Proposition \ref{som}. \end{proof}

\chapter{Dimension $48$}\label{48}

\section{First reduction} Let $H$ be a nontrivial semisimple Hopf algebra of dimension
$48$.

\begin{lemma}\label{conteo48} The order of $G(H^*)$ is either $2$, $3$, $4$, $6$, $8$, $12$, $16$ or $24$.
As an algebra, $H$ is of one of the following types:
\begin{itemize} \item $(1, 2; 2, 3; 3, 2; 4, 1)$, \item $(1, 3; 3, 1; 6, 1)$,
\item $(1, 3; 2, 9; 3, 1)$, \item $(1, 3; 3, 5)$,
\item $(1, 4; 2, 2; 3, 4)$, \item $(1, 4; 2, 3; 4, 2)$, \item $(1, 4; 2, 7; 4, 1)$, \item $(1, 4; 2, 11)$, \item $(1, 4; 2, 2; 6, 1)$, \item $(1, 6; 2, 6; 3, 2)$, \item  $(1, 8; 2, 2; 4, 2)$, \item $(1, 8; 2, 10)$, \item $(1, 8; 2, 6; 4, 1)$, \item $(1, 12; 2, 9)$, \item $(1, 12; 3, 4)$, \item $(1, 12; 6, 1)$, \item $(1, 16; 2, 8)$, \item $(1, 16; 4, 2)$, \item $(1, 24; 2, 6)$. \end{itemize} \end{lemma}

\begin{proof} The possibility $|G(H^*)| = 1$ is discarded by \ref{alg-struct}
and  Theorem \ref{thm-nr}; see Corollary \ref{cor-g-1}.

The  only possibilities with $|G(H^*)| = 2$ are the types
$(1, 2; 2, 7; 3, 2)$ and $(1, 2; 2, 3; 3, 2; 4, 1)$. In the first
case, $H$ cannot have Hopf algebra quotients of type $(1, 3; 3,
1)$. Therefore, by Remark \ref{el-nr} (i) every irreducible
character of degree $2$ is stable under left multiplication by
$G(H^*)$; hence by Theorem \ref{corolario} there is a quotient Hopf
algebra of dimension $30$, which is impossible.

Suppose that $|G(H^*)| = 16$. The only  possibility excluded in
the list is the type $(1, 16; 2, 4; 4, 1)$. Suppose on the
contrary that $H$ is of this type. Then $H$ has four irreducible
characters of degree $2$ and one irreducible character of degree
$4$. In particular, the action of $G(H^*)$ on $X_2$ is transitive
and $|G[\chi]| = 4$, for all irreducible character of degree $2$.
By Remark \ref{nr} (iii), there is a quotient Hopf algebra of
dimension $32$, which is impossible. The rest of the lemma follows
from \ref{alg-struct}. \end{proof}

\begin{remark}\label{rmk-48} (i) If $H$ is of type $(1, 24; 2, 6)$, then $H$ is not simple, by Corollary \ref{kob-mas}.

  (ii) Suppose that $H$ is  of type $(1, 4; 2, 3; 4, 2)$
as a coalgebra.  By Theorem \ref{corolario} the irreducible
characters of degrees $1$ and $2$ give rise to a Hopf subalgebra
of dimension $16$.

(iii) If $H$ is of type  $(1, 8; 2, 2; 4, 2)$  as a
coalgebra then the irreducible characters of degrees $1$ and $2$
give a Hopf subalgebra of dimension $16$; see Remark \ref{nr} (iii).

(iv) Suppose that $H$ is of type  $(1, 2; 2, 3; 3, 2;
4, 1)$ as a coalgebra. Then $H$ is not simple.

\begin{proof} We have that $H$ does not contain Hopf subalgebras
of dimension $12$ \cite{fukuda}. On the other hand, there must
exist an irreducible character $\chi$ of degree $2$, such that
$G[\chi] = 1$ (since otherwise there would be a Hopf subalgebra of
dimension $14$). By Theorem \ref{thm-nr}, $H$ has a Hopf
subalgebra of dimension $24$ and then $H$ is not simple.
\end{proof}

(v) Suppose that $H$ is of type $(1, 4; 2, 2; 3, 4)$
as a coalgebra. Then $H$ contains unique Hopf subalgebras $A_1
\subseteq A$ of dimension $6$ and $12$, respectively.

\begin{proof} The irreducible characters of degrees $1$ and $2$
give rise to a (unique) Hopf subalgebra of dimension $12$.

Let  $\psi_1, \dots, \psi_4$ be the irreducible characters of
degree $3$, and let $\chi_1, \chi_2$ be the irreducible characters
of degree $2$. We have $\psi_1^* \psi_1 = 1 + \chi_1 + \psi +
\psi'$, where $\psi, \psi'$ are irreducible of degree $3$; this
implies that $\chi_i^* = \chi_i$, $i = 1, 2$. Thus $\psi_1 \chi_1
= \psi_1 +  \psi_1 a$, where $1 \neq a \in G[\chi_1]$. In
particular, $|G[\chi_1]| = |G[\chi_2]| = 2$.

Let  $\psi_i$ be any irreducible character of degree $3$; thus
$\psi_i = g\psi_1$, for some $g \in G(H)$. Hence $\psi_i^* \psi_i
= \psi_1^* \psi_1 = 1 + \chi_1 + \psi + \psi'$.

Then  $\chi_1^2 = 1 + a + \chi$, where $\chi$ is irreducible of
degree $2$. Comparing decompositions of $\psi_1\chi_1^2$, we find
that $\psi_1 \chi = \psi_1 + \psi_1 a$, implying that $m(\chi,
\psi_1^* \psi_1) = 1$ and in turn that $\chi = \chi_1$. Then
$G[\chi_1] \cup \{ \chi_1 \}$ give rise to a Hopf subalgebra $A_1$
of dimension $6$. We have  also $m(\chi_1, \chi_2\chi_1) =
m(\chi_2, \chi_1^2) = 0$. Therefore $m(\chi_2, \chi_2\chi_1) =
m(\chi_2, \chi_1\chi_2) = m(\chi_1, \chi_2^2) > 0$, implying that
$m(\chi_2, \chi_2\chi_1) = m(\chi_1, \chi_2^2) = 1$. Thus $\chi_2
= b \chi_1$, $b \neq a \in G(H)$ and $\chi_2^2 = 1 + a + \chi_1$.
In particular, $A_1$ is the only six-dimensional Hopf subalgebra
of $H$ as claimed. \end{proof}

Note that the set $\{ \psi, \psi' \}$ is stable under the adjoint
action of $G(H)$ and also under $*$. Also, if $\psi = \psi'$, then
$H$ is not simple; indeed, in this case, $G[\chi_1]$, $\chi_1$,
$\psi$ and $a\psi = \psi a$ span a standard subalgebra
corresponding to a Hopf subalgebra of dimension $24$.

(vi) Suppose that $H$ is of type $(1, 4; 2, 2; 6, 1)$ as a
coalgebra. Then $H$ contains a unique Hopf subalgebra $A$ of
dimension  $12$, which coincides with the sum of simple
subcoalgebras  of dimension 1 and 4.\end{remark}

\begin{lemma} Suppose that $H$ is of type  $(1, 4; 2, 7; 4, 1)$ as a coalgebra.
Then $H$ is not simple. \end{lemma}

\begin{proof} By Proposition \ref{cociente8}, $H$ contains
a Hopf subalgebra $A$ of dimension $8$. In particular, the group
$G(H) = G(A)$ is not cyclic.

Let $\zeta \in H$ be the unique irreducible character of degree 4.
Hence we have $g \zeta = \zeta = \zeta g$, for all $g \in G(H)$.
Let $\lambda \in A$ be the unique irreducible character of degree
2. We claim that \begin{equation*}\lambda \zeta = 2 \zeta = \zeta
\lambda. \end{equation*} Indeed, otherwise, there would exist
$\chi \in X_2$ such that $m(\chi, \lambda \zeta) > 0$; hence
$\lambda \zeta = \chi + \dots$. Multiplying on the left by
$\lambda$, and using that $\lambda^2 =  \sum_{g \in G(H)} g$, we
find that $\lambda \chi =  \zeta$ is irreducible. This contradicts
Lemma \ref{productodesimples}. Thus $\lambda \zeta = 2 \zeta$, and
then $\zeta \lambda = (\lambda \zeta)^* = 2\zeta^* = 2\zeta$.
Hence the claim is established.

This implies that $AC = C = CA$. It follows from Corollary
\ref{comm-pair-cor} and \cite[Theorem 4.8 (1)]{masuoka-cont} that
$A$ is commutative. Suppose that $B \subseteq H$ is another Hopf
subalgebra of dimension 8. Then the same argument applies, and
thus $k[A, B]C = C = Ck[A, B]$. Hence $\dim k[A, B] = 16$ and the
Hopf subalgebra $k[A, B]$ is normal in $k[C] = H$. Therefore $H$
is not simple in this case.

Suppose next that $B \subseteq H$ is a Hopf subalgebra of
dimension $6$; we may assume that $k[A, B] = H$. We have $G(B)\cap
Z(B) \neq 1$ and, since $A$ is commutative, this group is central
in $k[A, B] = H$. Hence $H$ is not simple in this case.

  Since $H$ cannot contain Hopf subalgebras of dimension
32, there exist $\chi, \psi \in X_2$ such that $\chi^* \psi =
\zeta$ is irreducible of degree 4. It follows from Lemma
\ref{productodesimples} that $G[\chi]$ and $G[\psi]$ are distinct
subgroups of order 2 of $G(H)$. Write $\psi \psi^* = 1 + a +
\tau$, where $\tau \in X_2$ and $G[\psi] = \{ 1, a \}$; thus $\tau
= \tau^*$ and $G[\psi] \subseteq G[\tau]$. Similarly, $\chi \chi^*
= 1 + b + \mu$, where $\mu \in X_2$ and $G[\chi] = \{ 1, b \}$;
thus $\mu = \mu^*$ and $G[\chi] \subseteq G[\mu]$. In view of
Lemma \ref{productodesimples}, we have $\mu \neq \tau$. Hence,
after changing if necessary $\psi$ and $\chi$, we may assume $\tau
\neq \lambda$.

Suppose $G[\tau] = G(H)$. Then $\tau$ and $G(H)$  span a standard
subalgebra corresponding to a Hopf subalgebra of dimension 8, and
we know $H$ is not be simple in this case. Hence we may assume
$|G[\tau]| = 2$. If $\tau = \psi g$, $g \in G(H)$, then $\tau^2 =
\psi g (\psi g)^* = 1 + a + \tau$. Thus $\tau$ and $G[\tau]$ span
a standard subalgebra corresponding to a Hopf subalgebra of
dimension 6, and $H$ is not simple.

Therefore we may assume that the orbits $\tau G(H)$ and $\psi
G(H)$ are disjoint. Then the orbits of the right action of $G(H)$
on $X_2$ are $$\lambda G(H), \quad \chi G(H), \quad \psi G(H),
\quad \tau G(H),$$ and the only orbits with (left) stabilizer
$G[\psi]$ are $\psi G(H)$ and $\tau G(H)$.

Note that if $\pi \in X_2$ is such that $m(\pi, \psi\zeta) > 0$,
then $m(\zeta, \pi^*\psi) = m(\pi, \psi\zeta) > 0$. Therefore
$G[\pi] \cap G[\psi] = 1$ and $\pi \in \chi G(H)$. Thus $\psi\zeta
= \chi + \chi t + \zeta$, for some $t \in G(H)$. In particular,
$m(\zeta, \zeta\psi^*) = m(\zeta, \psi\zeta) = 1$.

On the other hand, $\zeta\psi^* = \chi^*\psi\psi^* = \chi^* +
\chi^*a +\chi^*\tau$. Hence $\chi^*\tau = \zeta$ is irreducible.
Letting $X = \{ \rho \in X_2: \, \chi^*  \rho  = \zeta \}$, we
have  $X = \psi G(H) \cup \tau G(H)$. Moreover, we have $G(H) X =
X = X G(H)$: the left hand side equality, because $G[g\rho] =
G[\rho]$, for all $g \in G(H)$, and the only elements in $X_2$
with stabilizer $G[\psi]$ are in $\psi G(H) \cup \tau G(H) = X$.

Let $\sigma, \rho \in X' = X \cup \{ \lambda \}$. Then the product
$\sigma^* \rho \neq \zeta$. Note also that, for all $\rho \in X$,
we have $m(\chi, \rho\zeta) = m(\chi g, \rho\zeta) = m(\zeta,
\chi^* \rho) = 1$. Hence we have $\rho\zeta = \chi + \chi t +
\zeta$, and thus $m(\rho, \zeta^2) = 1$, for all $\rho \in X$.

  This allows us to write $\zeta^2 = \sum_{g \in G(H)} g +
2 \lambda + \sum_{\rho \in X} \rho$. In particular, it follows
that $X'$ is closed under $*$. If $\sigma \in X$, then
\begin{align*}\sigma \zeta^2  & = \sum_{g \in G(H)} \sigma g  + 2 \sigma \lambda  + \sum_{\rho \in X} \sigma \rho  \\
& = (\sigma \zeta)\zeta  = (\chi + \chi t + \zeta)\zeta  = 2 \chi
\zeta  + \zeta^2 \\ & = 2 \sum_{\nu \in X} \nu + \sum_{g \in
G(H)}g + 2 \lambda + \sum_{\rho \in X} \rho. \end{align*} This
implies that for all $\sigma, \rho \in X$, the product $\rho
\sigma$ decomposes as a sum of elements of $X' \cup G(H)$. Also,
since $\lambda \zeta^2 = 2 \zeta^2 = \zeta^2 \lambda$, then all
products $\lambda \sigma, \sigma \lambda$ decompose as  sums of
elements of $X$, for all $\sigma \in X$.

Hence $X' \cup G(H)$ spans a standard subalgebra of $R(H^*)$ which
corresponds to a Hopf subalgebra of dimension 24. This proves the
lemma. \end{proof}

\begin{lemma} Assume that $H$ is of type $(1, 6; 2, 6; 3, 2)$ as a coalgebra. Then $H$ is not simple. \end{lemma}

\begin{proof} Since the
dimension of $H$ is not divisible by $30$, there must exist $\chi
\in X_2$ such that $G[\chi] = 1$. By Theorem \ref{thm-nr}, $H$ has
a Hopf subalgebra $A$ of type $(1, 3; 3, 1)$ as a coalgebra.

Let $\psi \in \widehat{H^*}$ be the unique irreducible character
of  degree $3$ in $A$. Then we have  $\psi^2 = 1 + a + a^2 +
2\psi$, where $G[\psi] = \{ 1, a, a^2 \}$. Let $b \in G(H)$ of
order $2$, so that $b \psi = \psi' = \psi b$ is the remaining
degree $3$ irreducible character. We have $\psi' \psi$, $\psi
\psi'$ belong to the span of $G(H) \cup \{ \psi, \psi' \}$. Thus
$G(H) \cup \{ \psi, \psi' \}$ spans a standard subalgebra of $H$
which corresponds to a Hopf subalgebra of dimension $24$. This
implies that $H$ is not simple, as claimed. \end{proof}

\begin{lemma} Assume that $H$ is of type $(1, 3; 3, 5)$ as a coalgebra. Then $H$ is not simple. \end{lemma}

\begin{proof} Let $\psi \in H$ be an irreducible character of degree 3.
Decomposing the product $\psi \psi^*$, we see that $G[\psi] = G(H)$ is necessarily of order 3.

By Remark \ref{particular}, the quotient coalgebra $H/H(kG(H))^+$
is cocommutative.  In particular, we may assume that $|G(H^*)|$ is
not divisible by 3: otherwise $H$ would be a biproduct $H = R \#
kG(H)$, with $R \simeq H/H(kG(H))^+$ cocommutative, implying that
$H$ is not simple, by Proposition \ref{som}.

  By Lemma \ref{conteo48} and the previous lemmas,  we may
assume that there is a Hopf algebra quotient $H \to B$, where
$\dim B = 4$; hence $\dim H^{\co B} = 12$, and by \cite{NZ},
$kG(H) \subseteq H^{\co B}$.

Let $H^{\co B} = kG(H) \oplus V_1 \oplus V_2 \oplus V_3$ be a
decomposition of $H^{\co B}$ as a sum of irreducible left coideals
of $H$. If $V_1$, $V_2$ and $V_3$ are pairwise isomorphic, then
$H^{\co B}$ is a subcoalgebra of $H$ and the map $H \to B$ is
normal. Hence we may assume that $V_1 $ appears with multiplicity
1 in $H^{\co B}$.

  Suppose that $V_i$ appears with multiplicity 1 in
$H^{\co B}$.  Then necessarily $gV_i = V_i = V_ig$, for all $g \in
G(H)$, $i = 1, 2, 3$. Let $C_i \subseteq H$ be the simple
subcoalgebra containing $V_i$. By Corollary \ref{cor-mk} $kG(H)$
is normal in $k[C_i]$. Hence, we may assume that $\dim k[C_i] =
12$.

  We claim that $V_i$  appears with multiplicity 1, for
all $i$.  As above, implies that $G(H)$ is normal in $k[C_i]$ for
all $i$. Since $k[C_1, C_2, C_3] = H$, it follows that $G(H)$ is
normal in $H$ and we are done.

  To prove the claim we argue as follows.  By Lemma
\ref{idf-coinv}, $(H^{\co B})^* \simeq \Ind_{B^*}^{H^*}1$ as left
$H^*$-modules. Hence, $\Ind_{B^*}^{H^*}1 = \sum_{g \in G(H)}g +
\psi_1 + 2\psi_2$, where $\psi_1 \neq \psi_2 \in H$ are
irreducible characters of degree 3. In particular, $\psi_2^* =
\psi_2$.

We may also assume that the Hopf subalgebra $k[C]$, where $C$ is
the subcoalgebra  corresponding to $\psi_2$, is all of $H$;
otherwise, $\dim k[C] = 12$ and $G(H)$ would be normal in $k[C]$,
hence $G(H)$ would be normal in $H = k[C, C_1]$ (by dimension).

By Frobenius reciprocity, $\psi_2\vert_{B^*} = 2. 1 + t$, where $1 \neq t \in G(B)$. Since $\psi_2^* = \psi_2$, then $t^2 = 1$.
Let $A = k\langle t \rangle \subseteq B$. Consider the composition
$$\pi: H \to B \to \overline B : = B/BA^+.$$
We have $\psi_2\vert_{\overline B} = \pi(\psi_2) = 3. 1$.
Applying again the Frobenius reciprocity, this gives $m(\psi_2, \Ind_{{\overline B}^*}^{H^*}1) = 3 = \deg \psi_2$.

Therefore, by Lemma \ref{idf-coinv}, $H^{\co \overline B}$ contains the simple subcoalgebra $C$ corresponding to $\psi_2$. This is absurd, since it implies that $k[C] = H \subseteq H^{\co \overline B}$.
This finishes the proof of the lemma. \end{proof}

\begin{lemma} Assume that $H$ is of type $(1, 3; 2, 9; 3, 1)$ as a coalgebra. Then $H$ is not simple. \end{lemma}

\begin{proof} Let $\psi$ be the unique irreducible character of degree 3.
For all irreducible character $\mu$ of degree 2, we have $\mu
\mu^* = 1 + \psi$. By Theorem \ref{thm-nr}, $G(H)$ and $\psi$ span
a standard subalgebra of $R(H)$,  which corresponds to a
commutative Hopf subalgebra of coalgebra type $(1, 3; 3, 1)$. In
particular, the dimension  of an irreducible $H$-module is at most
$[H: A] = 4$. Thus $H^*$ cannot be of type $(1, 3; 3, 1; 6, 1)$.
Moreover, by previous results in this section, we may assume
that either $H^*$ is of type $(1, 3; 2, 9; 3, 1)$ as a coalgebra,
or else there is a Hopf algebra quotient $H \to B$, where $\dim B
= 4$.

  In the last case, we have $\dim H^{\co B} = 12$ and
$G(H) \subseteq H^{\co B}$, by \cite{NZ}. Also, unless $A = H^{\co
B}$ is normal in $H$, we may assume that $H^{\co B}$ contains an
irreducible left coideal $V$ of dimension 2. Therefore, $H^{\co B}
= kG(H) \oplus \oplus_{g \in G(H)}gV \oplus W$ as a left coideal of $H$, where $W$  is an irreducible left coideal of dimension 3.
This implies that $G(H) \chi = \chi G(H)$, where $\chi$ is the
character of $V$. Moreover, by Lemma \ref{idf-coinv}, there exists $g_0 \in G(H)$ such that
$\chi^* = g_0\chi$.

The relation $\chi^*\chi = 1 + \psi$ implies that $\psi \chi =
\sum_{g \in G(H)} g \chi$. Also, for all $g, h \in G(H)$, we have
$\chi h = h_0\chi$, for some $h_0 \in G(H)$ and therefore, $$(g
\chi) (h \chi) = gh_0 \chi^2 = gh_0(g_0)^{-1} \chi^* \chi,$$ so that $(g
\chi) (h \chi)$ belongs to the span of $G(H)$ and $\psi$.  It
follows that $G(H), \psi$ and $G(H) \chi$ span a standard
subalgebra, corresponding to a Hopf subalgebra of dimension 24.
Hence $H$ is not simple in this case.

\bigbreak Suppose finally that $H^*$ is of type $(1, 3; 2, 9; 3,
1)$ as a coalgebra. Then there is a projection $q: H \to B$, where
$B = kG(B)$ is of dimension 12. The coalgebra structure of $B^*$
implies that $G(B)$ has a unique normal subgroup of order 4. We
have necessarily $H^{\co B} = k1 \oplus V$, where $V$ is an
irreducible left coideal of dimension 3; whence $\dim q(A) = 3$
and  $q(\psi) = \sum_{x \in G(q(A))}x$.

Since $|X_2|$ is odd, there exists an irreducible character $\chi$
of degree 2, such that $\chi^* = \chi$. Then $q(\chi) = a + b$,
where $a, b \in G(B)$. The relation $\chi^2 = 1 + \psi$ implies
that $G(q(A)) \subseteq \langle a, b \rangle$. We cannot have $a^2
= b^2 = 1$, since otherwise $\langle a, b \rangle$ would be
contained in the unique subgroup of order 4 of $G(B)$. Hence $b =
a^{-1}$, and $G(q(A)) = \langle a \rangle$ of order 3.

Let $C$ be the simple subcoalgebra containing $\chi$. Then
$q(k[C]) = \langle a \rangle$; in particular, $k[C] \neq H$. In
addition, $A \subseteq k[C]$, so that $12 / \dim k[C]$. Hence,
$\dim k[C] = 24$ and thus $H$ is not simple. This finishes the
proof of the lemma. \end{proof}

\begin{lemma}\label{48-8} Let $H$ be of
type $(1, 8; 2, 2; 4, 2)$ or $(1, 8; 2, 10)$ as a coalgebra.
Assume $H$ is simple. Then $H$ contains a Hopf subalgebra of
dimension $16$.
\end{lemma}

\begin{proof} If $H$ is of
type $(1, 8; 2, 2; 4, 2)$, the claim follows from  Remark
\ref{rmk-48} (iii). So assume $H$ is of type  $(1, 8; 2, 10)$. For
all $\lambda \in X_2$ we have $|G[\lambda]| = 2$ or $4$.

Suppose $\chi \in X_2$ is such that $G[\chi] = G[\chi^*]$ of order
4. Let  $A = kG[\chi]$ and $C$ the simple subcoalgebra containing
$\chi$. We have $AC = C = CA$, whence $A$ is normal in $k[C]$ by
Proposition \ref{dima=dimc}. Since $|G[\chi]| = 4$, $A$ is normal
in $kG(H)$ and then it is also normal in $K = k[C, G(H)]$. But
$\dim K > 8$ is divisible by 8, and we are assuming that $H$ is
simple. Hence $\dim K = 16$ and we are done.

 Assume first that $|G[\lambda]| = 4$ for all $\lambda
\in X_2$. We claim that there exists $\chi \in X_2$ such that
$G[\chi] = G[\chi^*]$, implying the lemma. To prove the claim,
consider the action of $G(H) \times G(H)$ on $X_2$ given by $(g,
h).\chi = g \chi h^{-1}$. The orbit of an element $\chi$ is
$G(H)\chi G(H)$, so it has order 2 or 4, and clearly $G(H)\chi
G(H)$ and $G(H)\chi^* G(H)$ are of the same order, for all $\chi$.
Also, because $G[\chi]$ is normal in $G(H)$, $G[\lambda] =
G[\chi]$ for all $\lambda \in G(H)\chi G(H)$. Suppose on the
contrary that $G[\chi] \neq G[\chi^*]$ for all $\chi \in X_2$.
Then $G(H)\chi G(H)$ and $G(H)\chi^* G(H)$ are disjoint. Let
$\alpha \in X_2$ such that $\alpha$ is not conjugate to $\chi$ nor
to $\chi^*$ (such an $\alpha$ exists because $|X_2| = 10$). Then
the orbits $G(H)\chi G(H)$, $G(H)\chi^* G(H)$, $G(H)\alpha G(H)$
and $G(H)\alpha^* G(H)$ are pairwise disjoint.

Then $|X_2| = 10 \geq  2|G(H)\chi G(H)| + 2 |G(H)\alpha G(H)| \geq
8$, and there must exist $\beta$ not in any of these orbits.
Hence, by dimension, $G(H)\beta G(H) = G(H)\beta^* G(H)$. Thus
$G[\beta] = G[\beta^*]$ and the lemma holds in this case.

  Suppose next that $|G[\lambda]| = 2$ for some $\lambda
\in X_2$. Then $\lambda \lambda^* = \sum_{g \in G[\lambda]}g +
\chi$ for some $\chi \in X_2$ such that $\chi^* = \chi$. If
$|G[\chi]| = 4$, then we are done. Otherwise, we may assume $\chi
G(H)$ is the only orbit with 4 elements: if not, since $|X_2| =
10$, there would be a unique orbit $\alpha G(H)$ with stabilizer
of order 4, and thus $G[\alpha^*] = G[\alpha]$ implying the lemma.

In particular, for all $g \in G(H)$, $g\chi g^{-1} = \chi a$, $a
\in G(H)$, and therefore $G[g\chi g^{-1}] = g G[\chi]g^{-1} =
G[\chi]$. That is, $G[\chi] \trianglelefteq G(H)$.

Also $\lambda = \chi t$, $t \in G(H)$, and $\chi^2 = \sum_{g \in
G[\chi]}g + \chi$. Hence $G[\chi]$ and $\chi$ span a standard
subalgebra corresponding to a Hopf subalgebra $A$ of dimension 6.
Moreover, $G[\chi]$ is central in $A$. Then $kG[\chi]$ is normal
in $k[A, G(H)]$. But this implies $\dim k[A, G(H)] = 24$ or 48,
and  $H$ is not  simple in this case.  \end{proof}

\begin{lemma} Suppose  that $H$ is of type  $(1, 8; 2, 6; 4, 1)$ as a coalgebra. Then $H$ is not  simple.  \end{lemma}

\begin{proof} Let $\zeta \in C$ be the unique irreducible character of degree 4 in $H$.
Let $X = \{ \chi \in X_2: m(\chi, \zeta^2) > 0 \}$, so that
$\zeta^2 = \sum_{g \in G(H)}g + \sum_{\chi \in X}n_{\chi}\chi +
n\zeta$, where $n_{\chi} = m(\chi, \zeta^2) = m(\zeta, \chi\zeta)
= 1$ or 2,  $\chi \in X$. We have $X$ is stable under $*$ and left
and right multiplication by $G(H)$.

If $X = \emptyset$ then $G(H)$ and $\zeta$ span a standard
subalgebra corresponding to a Hopf subalgebra of dimension $24$.
Hence $H$ is not simple in this case. Thus we may assume $X \neq
\emptyset$.

Let $\chi \in X$. We have $\chi\zeta = n_{\chi} \zeta + \sum
\lambda$, where $\lambda$ runs over the set $\Lambda = \{ \lambda
\in X_2: \, m(\lambda, \chi\zeta) = m(\zeta, \lambda^*\chi) > 0
\}$. Then $\deg \chi \zeta = 4n_{\chi} + 2 |\Lambda|$. Note that
$\Lambda G(H) = \Lambda$, so that $|\Lambda| = 0$ (iff $n_{\chi} =
2$) or $2$ (iff $n_{\chi} = 1$).

If $|\Lambda| = 0$, we have $n_{\chi} = 2$. Since $g\chi$ appears
in $\zeta^2$ with the same multiplicity as $\chi$, for all $g \in
G(H)$. Then $|G[\chi]| = 4$ and  $\zeta^2 = \sum_{g \in G(H)}g +
2\chi + 2t\chi$, for some $t \in G(H)$. In particular $G(H)\chi =
\chi G(H)$. Also $\chi \zeta^2 = 2\zeta^2$. Then $G(H)$,
$G(H)\chi$ and $\zeta$ span a standard subalgebra corresponding to
a Hopf subalgebra of dimension 32. This is impossible by
\cite{NZ}. Hence $|\Lambda| = 2$; say $\Lambda = \{ \lambda,
\lambda a \}$, $a\in G(H)$. In particular, $|G[\lambda^*]| =
|G[\lambda]| = 4$. Then $n_{\chi} = 1$, for all $\chi \in X$. So
that $\zeta^2 = \sum_{g \in G(H)}g + \chi_1 + \chi_2 + \zeta$,
with $X = \{ \chi_1, \chi_2 \}$, if $|X| = 2$, or $\zeta^2 =
\sum_{g \in G(H)}g + \chi_1 + \chi_2 + \chi_3 + \chi_4$, with $X =
\{ \chi_1, \chi_2, \chi_3, \chi_4 \}$, if $|X| = 4$.

Then, for each $\chi \in X$, $\chi\zeta = \zeta + \lambda +
\lambda a$, and thus $\chi\zeta^2 = \zeta^2 + \lambda (1 + a)
\zeta = \zeta^2 + 2\lambda \zeta$. Now since $\lambda^*\chi =
\zeta$, then $\lambda \zeta = \lambda\lambda^*\chi \in G(H)\chi$.
Hence $\chi\zeta^2$ belongs to the span of $G(H)$, $\zeta$ and
$X$. Since this happens for all $\chi \in X$, then the irreducible
summands of $\zeta^2$ span a standard subalgebra corresponding to
a Hopf subalgebra of dimension $32$ or $24$. By \cite{NZ} the
first case is impossible, and the second one implies that $H$ is
not simple.  \end{proof}

\begin{lemma} Suppose $H$ is of type $(1, 3; 3, 1; 6, 1)$ as a coalgebra. Then $H$ is not simple. \end{lemma}

\begin{proof} The irreducible characters of degrees 1 and 3 span a standard subalgebra of $R(H)$, which corresponds to a Hopf subalgebra $A$ of dimension  $12$. By \cite{fukuda} $A$ is commutative. In particular, $\dim V \leq [H:A] = 4$, for all irreducible $H$-module $V$.

Suppose on the contrary that $H$ is simple. Then there is no Hopf
algebra quotient $H \to B$ with $\dim B = 16$: otherwise, by
\cite{NZ}, $kG(H) = H^{\co B}$ is a normal Hopf subalgebra of $H$.
In view of Lemma \ref{conteo48} and the previous results, it
follows that $\dim V = 1, 2$ or $3$ for all irreducible $H$-module
$V$. It follows as well that $|G(H^*)| = 4$ or $12$, and thus
$G(H^*)$ contains a subgroup $\Gamma$ of order 4.

  Consider the Hopf algebra projection $H^* \to A^*$. Then
$(H^*)^{\co A^*} \subseteq (H^*)^{\co q}$, where $q: H^* \to
k^{G(H)}$ is the natural projection. By \cite{NZ}, we have
$k\Gamma \subseteq (H^*)^{\co q}$. In particular, since
$(H^*)^{\co q} \neq kG(H^*)$, $\Gamma = G(H^*) \cap (H^*)^{\co q}$
is a normal subgroup of $G(H^*)$, and it is therefore its only
subgroup of order 4.

\begin{claim} $(H^*)^{\co q}$ contains no irreducible left coideal of dimension 3. \end{claim}

\begin{proof} Suppose otherwise that $U \subseteq (H^*)^{\co q}$ is an irreducible left coideal of dimension 3.
Hence $H^*$ is of type $(1, 4; 2, 2; 3, 4)$ or $(1, 12; 3, 4)$ as a coalgebra,
and there is a Hopf subalgebra $B \subseteq H^*$ with $\dim B = 12$, such that $B$ is cocommutative or of type $(1,  4; 2, 2)$.

Since $g U$ is not isomorphic to $U$, for all $g \in \Gamma$, the
sum $\sum_{g \in \Gamma}gU$ is direct and by dimension,
$(H^*)^{\co q} = k\Gamma \oplus \oplus_{g \in \Gamma}gU$. Hence,
if $H$ is simple, $(H^*)^{\co A^*} = k1 \oplus gU$, for some $g
\in \Gamma$. This implies that $B^{\co A^*} = B \cap (H^*)^{\co
A^*} = k1$, and therefore that $B \simeq A^*$ is cocommutative and
$H^* \simeq R \# B$ is a biproduct, where $R$ is a  braided Hopf
algebra over $B$ of dimension 4. Proposition \ref{dimR-3-4}
implies that $H$ is not simple in this case. Hence the claim
follows. \end{proof}

In view of the claim, we may assume that $(H^*)^{\co A^*} = k1 \oplus kt \oplus V$, where $1 \neq t \in G(H^*)$, and  $V$ is an irreducible left coideal of dimension 2. In particular, $tV = V = Vt$.  So that, letting $C$ be the simple subcoalgebra of $H^*$ containing $V$, we find that $t$ is central in the Hopf subalgebra $k[C]$,  by Corollary \ref{cor-mk}. In particular, $k[C] \neq H^*$ and by \cite{NZ},  $\dim k[C] = 12$ or 8.
In the first case, the restriction $q\vert: k[C] \to k^{G(H)}$ is surjective; then by \cite{fukuda} $k[C]$ is trivial, whence necessarily commutative. But this implies that $\dim U \leq [H^*: k[C]] = 4$, for all irreducible $H$-comodules $U$, against the assumption on the coalgebra type of $H$.
Hence $\dim k[C] = 8$ and  $k[C]  \subseteq (H^*)^{\co q}$.

We have therefore a decomposition $(H^*)^{\co q} = k[C] \oplus V_1
\oplus V_2 \oplus V_3 \oplus V_4$, where $V_i$ is an irreducible left coideal of
dimension 2 of $H^*$, for all $i  = 1, \dots, 4$.

Note that $V_i$ appears in $(H^*)^{\co q}$ with multiplicity
1, for all $i$. Otherwise, say $C_1 \subseteq (H^*)^{\co q}$,  where $C_1$ is the simple subcoalgebra containing $V_1$.  Then $k[C, C_1] \subseteq (H^*)^{\co q}$; but the inclusion $k[C] \subseteq k[C, C_1]$ is strict, so $k[C, C_1] = (H^*)^{\co q}$ by dimension restrictions,  implying that $H$ is not simple.

Let $\tau_i \in X_2$  be the irreducible character  corresponding to $V_i$, $i  = 1, \dots, 4$, and let $\lambda$ be the character of $V$. By Frobenius reciprocity $m(1, q(\tau_i\tau_i^*)) = 2$, thus $G[\tau_i] \neq \Gamma$, because $q(g) = 1$, for all $g \in \Gamma$. Hence $G[\tau_i] $ is a subgroup of order 2, for all $i  = 1, \dots, 4$. Similarly, we see that $m(\tau_i, \lambda \tau_i) = m(\lambda, \tau_i\tau_i^*) = 0$

  Without loss of generality we may write $\tau_2 = h\tau_1$ and $\tau_4 = g\tau_3$, for some $h, g \in \Gamma \backslash \{ 1\}$. By the above, $\lambda  \tau_1 = \tau_3 + \tau_4$. Write $q(\tau_1) = 1 + x$, where $1  \neq x$ is of order 3. Since $q(\lambda) = 2.1$, we find that $q(\tau_i) = 1 + x$, for all $i  = 1, \dots, 4$.
This is impossible, since by Lemma \ref{idf-coinv}, we must have $\tau_1^* = \tau_i$ for some $i$, and since $q(\tau_1^*) = 1 + x^{-1}$. This contradiction finishes the proof of the lemma. \end{proof}

\begin{lemma}\label{12-no} Suppose that $H$ is simple. Then $|G(H)| \neq 12$. \end{lemma}

This discards the possibility that the types $(1, 12; 2, 9)$, $(1, 12; 3, 4)$ and $(1, 12; 6, 1)$ in Lemma \ref{conteo48} correspond to a simple Hopf algebra.

\begin{proof} Suppose on the contrary that  $|G(H)| = 12$.

\begin{claim}\label{g-h^*}  $|G(H^*)| = 12$. \end{claim}

\begin{proof} There is a Hopf algebra surjection $H^* \to B$, where
$B \simeq k\mathbb Z_3$.  Therefore, $H^*$ does not contain any Hopf subalgebra $A$ of dimension $16$. Indeed, if this were the case,  necessarily $(H^*)^{\co B} = A$, by \cite{NZ}; hence  $H$ would not be simple.

Suppose that $|G(H^*)| \neq 12$. By previous
results, $H^*$ must be of type $(1, 4; 2, 2; 3, 4)$, $(1, 4; 2, 11)$ or $(1, 4; 2, 2; 6, 1)$ as a coalgebra.

Consider first the case where $H^*$ is of one of the types $(1, 4; 2, 2; 3, 4)$ or $(1, 4; 2, 2; 6, 1)$; so that $H^*$ contains a Hopf subalgebra $A$ of dimension $12$, which is of type $(1, 4; 2, 2)$ as a coalgebra.
Since $kG(A) \subseteq
 (H^*)^{\co B}$, and $\dim (H^*)^{\co B} = 16$, then $(H^*)^{\co B} \cap A = kG(H^*)$. Let $\pi: H^* \to k^{G(H)}$ be the natural Hopf
algebra projection.
Then we have $(H^*)^{\co \pi} \subseteq (H^*)^{\co B}$ and $\dim (H^*)^{\co \pi} = 4$. Hence,
$(H^*)^{\co \pi} = k1 \oplus W$, for some
irreducible left coideal $W$ of dimension $3$.
(In particular, $H^*$ is not of type $(1, 4; 2, 2; 6, 1)$.)

Thus $(H^*)^{\co \pi} \cap A = k1$ and the restriction $\pi\vert_A: A \to k^{G(H)}$ is an isomorphism. This implies that $H$
is a biproduct $H = R \# kG(H)$. Then the claim follows in this case from
Proposition \ref{dimR-3-4}.

Suppose next that $H^*$ is of type $(1, 4; 2, 11)$.
By Proposition \ref{cociente8}, there is a Hopf algebra quotient $H \to B'$ where $B'$ is of algebra type $(1, 4; 2, 1)$; so that $\dim H^{\co B'} = 6$. This implies that $H$ is not of type $(1, 12; 6, 1)$.

By \cite{NZ}, any subgroup of order 3 of $G(H)$ is contained in $H^{\co B'}$. Hence, $G(H)$ contains a unique subgroup $F$ of order 3, and $F$ is the unique subgroup of $G(H)$ contained in $H^{\co B'}$.

Suppose $H$ is of type $(1, 12; 3, 4)$. Then $H^{\co B'} = kF \oplus V$, where $V$ is an irreducible left coideal of $H$ of dimension $3$. Then we have $gV = V = Vg$, for all $g \in F$. Let $C$ be the simple subcoalgebra of $H$ containing $V$. By Corollary \ref{cor-mk}, $kF$ is normal in $k[C]$.

Besides, since $H^{\co B'}$ is normal in $H$, we have $gCg^{-1} = C$, for all $g \in G(H)$. Since $k[C]$ and $G(H)$ necessarily generate $H$ as an algebra, it follows that $k[C]$ is normal in $H$. This discards this type as the coalgebra structure of $H$.

Finally, suppose $H$ is of type $(1, 12; 2, 9)$.
Then we must have $H^{\co B'} = kF \oplus \oplus_{g \in F} gV$, where $V$ is an irreducible left coideal of $H$ of dimension $2$.
This is not possible, since $\dim H^{\co B'} = 6$.
Thus  $H$ is not type $(1, 12; 2, 9)$ and the proof is complete.
 \end{proof}

Claim \ref{g-h^*} implies that $H$  is not of type $(1, 12; 6, 1)$ as a coalgebra: indeed, in this case, $H^{\co k^{G(H^*)}} \subseteq kG(H)$ and
therefore $kG(H^*)$ is a normal Hopf subalgebra of $H^*$.

Since $|G(H)| = |G(H^*)| = 12$,  then $H = R \# k\mathbb Z_3$ is a
biproduct. In particular, $G(H)$ contains a unique (normal) subgroup
of order $4$. There is in addition a Hopf algebra surjection $H
\to B$, where $\dim B = 4$. Thus if $F \subseteq G(H)$ is a
subgroup of order $3$, then $F \subseteq H^{\co  B}$. This shows
that $G(H)$ also contains a unique (normal) subgroup of order $3$.
Thus $G(H)$ is abelian. This implies that $H$ is not of type $(1,
12; 3, 4)$, in view of Proposition \ref{som}, since in this case $R$ will be a
cocommutative coalgebra, by Remark \ref{particular}.

  It remains to discard the case where  $H$ is of type
$(1, 12; 2, 9)$. Let $\Gamma \subseteq G(H)$ be the unique
subgroup of order $4$. Then there is an irreducible character
$\lambda$ of degree $2$, such that $\Gamma \cup \{ \lambda \}$
spans a standard subalgebra of $R(H^*)$, corresponding to  a  Hopf
subalgebra $A$ of dimension $8$.

By \cite{NZ}, we must have $A \subseteq R$. Also, since $R$ is not a Hopf subalgebra, $R$
can contain only one simple subcoalgebra of dimension $4$, which
implies that this subcoalgebra (and hence all of $A$) is stable
under the action of $\mathbb Z_3$ by conjugation.

Hence $B = k[A, \mathbb Z_3]$ contains $A$ as a normal Hopf subalgebra, and thus $B \neq H$. Since $\dim B$ is divisible by $3$ and  $8 = \dim A$, we get $\dim B = 24$. This implies that $H$ is not simple. The lemma follows. \end{proof}

\section{Further reductions} In this section we further reduce the possibilities for the (co)algebra structure of an eventual simple $H$.

\begin{lemma}\label{tw-2} Suppose that $H$ is of type $(1, 4; 2, 11)$ as a coalgebra.  Then we have

(i) There exists an 8-dimensional non-cocommutative Hopf subalgebra $A_0 \subseteq H$.

(ii) Assume in addition that $H$ is simple and contains a Hopf
subalgebra $A$ of dimension 12. Then neither $A$ nor $A_0$ is
commutative. In this case,  there exists a twist $\phi \in kG(H)
\otimes kG(H)$ such that $|G(H_{\phi})| = 12$, $24$ or $48$. In
particular, $H_{\phi}$ is not simple. \end{lemma}

\begin{proof} Part (i) follows from Proposition \ref{cociente8}. We show part (ii).  Note that a Hopf subalgebra of $H$ of dimension 12 is necessarily of coalgebra type $(1, 4; 2, 2)$. Let $B = k[A_0, A]$ be the subalgebra generated by $A_0$ and $A$, so that $B$ is a Hopf subalgebra of $H$. Since $8 = \dim A_0$ divides the dimension of $B$ and also $3 / \dim A / \dim B$, then $24$ divides the dimension of $B$.
$H$ being simple by assumption, we may assume that $B = H$. In other words, $A_0$ and $A$ generate $H$ as an algebra.

Also $G(A_0) =  G(A) = G(H)$ is not cyclic.

Suppose that $A$ is commutative. We know that there exists a central group-like element $1 \neq g \in G(A_0)$. Since $A$ is commutative, and $G(A_0) =  G(A)$, then $g$ commutes with $A$. Therefore $g$ is central in $k[A_0, A] = H$, contradicting the simplicity of $H$. A similar argument shows that $A_0$ cannot be commutative.  The last statement of the lemma follows from Proposition \ref{twist-pqq}.
\end{proof}

\begin{lemma}\label{bip-16} Suppose $H = R \# A$ is a biproduct, where $\dim A = 16$. Then $H$ is not simple. \end{lemma}

\begin{proof} We have $\dim R = 3$,  thus $R$ is commutative and cocommutative.
As a left coideal of $H$,  we may assume that $R = k1 \oplus V$,
where $V$ is an irreducible left coideal of dimension $2$.

If $A$ is cocommutative, then  the lemma follows from Proposition
\ref{dimR-3-4}.  So we may assume $A$ is not cocommutative. By
Lemma \ref{bip-2-dim}, $\rho (R) \subseteq kG(A) \otimes R$.

It is not difficult to see that  there must exist a normal Hopf
subalgebra $B$ of $A$ such that $kG(A) \subseteq B$ and $\dim B =
8$. Thus $\rho (R) \subseteq B \otimes R$ and $R \# B$ is a normal
Hopf subalgebra (of index $2$) of $H$.
\end{proof}

\begin{lemma}\label{ns-16} Suppose that $H$ contains a Hopf subalgebra
$A$ with $\dim A = 16$. If there is a quotient Hopf algebra $q: H \to B$, with $\dim B = 16$, then $H$ is not simple. \end{lemma}

\begin{proof}  We have $\dim H^{\co B} = 3$. Hence, by Lemma \ref{restriction}, $A \cap H^{\co B} = k1$,
thus $H = R \# A$ is a biproduct. The lemma follows from Lemma \ref{bip-16}.  \end{proof}

\begin{lemma}\label{nq-16} Suppose that $H$ is of coalgebra type $(1, 4; 2, 2; 3, 4)$ or $(1, 4; 2, 2; 6, 1)$. Assume in addition that there is a Hopf algebra quotient $q: H \to B$, with $\dim B = 16$. Then $H$ is not simple. \end{lemma}

\begin{proof} We consider first the case where $H$ is of type $(1, 4; 2, 2; 6, 1)$.  Let $\psi$ be the unique irreducible character of degree 6, and let $\chi_1, \chi_2$ be the irreducible characters of degree 2.
Then $\chi_1, \chi_2 \in A$, where $A$ is the unique Hopf subalgebra of dimension 12 of $H$; see Remark \ref{rmk-48} (vi). We shall show that $A$ is normal in $H$.

It is not hard to see that
\begin{equation}\label{psi2}\psi^2 = \sum_{g \in G(H)}g + 2 \chi_1 + 2 \chi_2 + 4 \psi.\end{equation}

Necessarily, we must have $H^{\co B} = k1 \oplus V$, where $V$ is an irreducible left coideal of $H$ of dimension 2.
In particular, $H^{\co B} = A^{\co B}$ and $q(A)$ is a four-dimensional Hopf subalgebra of $B$.

We first claim that $B$ is of type $(1, 4; 2, 3)$ as a coalgebra. To prove this, we consider the possible decompositions of $q(\psi) = \psi\vert_{B^*}$ into a sum of irreducible characters in $B$. It follows from Frobenius reciprocity, together with equation \eqref{psi2}, that $q(\psi) = \lambda_1 + \lambda_2 + \lambda_3$, where $\lambda_i$ are pairwise distinct irreducible characters of degree 2 in $B$. Hence $B$ must be of the prescribed coalgebra type.

Therefore, $q(A) = kG(B)$; since this is the unique Hopf subalgebra of dimension 4 in $B$.

On the other hand, it follows from the classification results in \cite{kashina}, that $kG(B)$ is normal in $B$. For the sake of completeness, we give a proof of  this fact in what follows. We may assume that $B$ is not commutative. By \cite[3.3]{kashina}, $B$ fits into a cocentral abelian extension
$$1 \to K \to B \to kF \to 0,$$
where $F = \langle t: t^2 = 1\rangle \simeq \mathbb Z_2$ and $K$ is a commutative Hopf algebra of dimension 8. In particular, $kG(B)$ is contained in $K$ and therefore $G(B) = G(K)$ is central $K$.

Let $e = \sum_{x \in G(B)}g$. As an algebra, $H$ is a smash product $H = K \# kF$, with respect to an action $\rightharpoonup : F \times K \to K$ by Hopf algebra automorphisms. Hence, the action of $t \in F$ permutes the elements of $G(B)$, and therefore $t \rightharpoonup e = e$.

Note that $H$ is generated as algebra by $K$ and $T : = 1 \# t$. It follows from the above discussion that
$$Te = (t \rightharpoonup e) \# t = e \# t = Te.$$
Hence $e \in Z(H)$, which proves the desired fact.

Consider the sequence of surjective Hopf algebra maps
$$\begin{CD}H @>q>> B @>q'>> B', \end{CD}$$
where $B' = B / B (kG(B))^+$.  Since $q(A) = kG(B)$, we have $A \subseteq H^{\co q'q}$. Thus, by dimension, $A = H^{\co q'q}$ and it is a normal Hopf subalgebra, as claimed.

\bigbreak
Now consider the case where $H$ is of type $(1, 4; 2, 2; 3, 4)$. We shall keep the notation in Remark \ref{rmk-48} (v).

For all irreducible character $\psi_i$ of degree 3, we have relations
\begin{equation}\label{relgral}\psi_i^* \psi_i = 1 + \chi_1 + \psi +
\psi', \end{equation} and $\chi_1 \psi_i = \psi_i + a \psi_i$, where $\psi, \psi'$ are fixed irreducible of degree $3$, $G[\chi_1] = \{ 1, a \}$, and  $\chi_1 \in A_1$, where $A_1$ is the
unique six-dimensional Hopf subalgebra of $H$.

Observe first that $\psi' =  \psi g$, for some $g \in G(H)$. We
shall show that $g \in G[\chi_1]$.  So that the irreducible
characters $1, a, \chi_1, \psi, a \psi$, span a standard
subalgebra of $R(H^*)$ corresponding to a Hopf subalgebra of
dimension  24, whence $H$ is not simple in this case.

We have $H^{\co B} = k1 \oplus V$, where $V$ is an irreducible left coideal of dimension 2, and necessarily $V \subseteq A_1$, since $\dim A_1$ does not divide $\dim B$. Hence, $q(\chi_1) = (\chi_1)\vert_{B^*} = 1 + \alpha$, where $1 \neq \alpha \in G(B)$. Since $a\chi_1 = \chi_1$, and $q(a) \neq 1$, then $q(a) = \alpha$.

On the other hand, $m(1, q(\psi_i)) = m(\psi_i,
\Ind_{B^*}^{H^*}1) = 0$, for all $i = 1, \dots, 4$. The relation \eqref{relgral} implies that $m(1, q(\psi_i^*\psi_i)) = 2$.

This implies that $q(\psi) = \beta + \lambda$ and $q(\psi') = \beta' + \lambda'$, for some $\beta, \beta' \in G(B)$, and $\lambda, \lambda' \in X_2(B)$. By Frobenius reciprocity, $\beta, \beta'$ and $\alpha$ are pairwise distinct elements of $G(B)$.

Using again equation  \eqref{relgral} for $\psi_i = \psi$, we find that $\lambda^* \lambda = 1 + \beta + \beta' + \alpha$; hence $\{ 1, \beta, \beta', \alpha \} = G[\lambda^*]$ is a subgroup of $G(B)$.
Therefore $\beta \alpha = \beta'$.

Using Frobenius reciprocity, we have
$$\psi a =  (\Ind_{B^*}^{H^*}\beta) a = \Ind_{B^*}^{H^*}(\beta a\vert_{B^*}) = \Ind_{B^*}^{H^*}(\beta \alpha) = \Ind_{B^*}^{H^*}(\beta') = \psi'.$$
This establishes the claim and finishes the proof of the lemma.
\end{proof}

\begin{lemma}\label{qt16} Suppose that $H$ contains a Hopf subalgebra
$A$ with $\dim A = 16$. Then we have:

(i) If $H$ is simple, then $H$ is of type  $(1, 4; 2, 11)$ as an algebra, and there exists a normalized 2-cocycle $\phi \in kG(H^*) \otimes kG(H^*)$ such that $(H^*)_{\phi}$ is not simple.

(ii) Assume in addition that $A$ is cocommutative. Then $H$ is not simple. \end{lemma}

By part (ii), if $H$ is of type $(1, 16; 2, 8)$ or $(1, 16; 4,
2)$, then $H$ is not simple. Therefore, if $H$ is simple and
contains a Hopf subalgebra of dimension 16, $H$ is of type  $(1,
8; 2, 2; 4, 2)$, $(1, 8; 2, 10)$, $(1, 4; 2, 3; 4, 2)$ or $(1, 4;
2, 11)$ as a coalgebra.

\begin{proof} (i) By previous lemmas, $|G(H^*)| \neq 2, 3, 6, 12$ and there is no Hopf subalgebra
$B \subseteq H^*$ with $\dim B = 16$.

Consider the projection  $\pi: H^* \to A^*$; we may assume
$(H^*)^{\co A^*} = k1 \oplus V$, where $V$ is a left coideal of
$H^*$ of dimension $2$.
 It follows from  Lemma \ref{conteo48} that $G[\chi] \neq 1$,
 where $\chi = \chi_V$ is the character corresponding to $V$. Then $|G[\chi]| = 4$ or $2$.

By Theorem \ref{coinvariantes}, $H^*$ contains a  Hopf subalgebra $B$ of
dimension $3|G(H^*)|$. We may assume that
$|G(H^*)| = 4$ or $8$, and therefore, $\dim B = 12$ or $24$.

Since $H$ is simple, $\dim B = 12$, and therefore $|G(H^*)| = 4$. It
follows from Lemma \ref{ns-16} and the results of the previous section, that $H^*$ is of type $(1, 4; 2, 2; 3, 4)$, $(1, 4; 2, 2; 6, 1)$ or $(1, 4; 2, 11)$ as a coalgebra.   In view of Lemma \ref{nq-16},  $H^*$ is of type $(1, 4; 2, 11)$ as a coalgebra.

Since  $H^*$ contains a Hopf subalgebra $A$ of dimension 12, by
Lemma \ref{tw-2},  there  exists a normalized 2-cocycle $\phi \in
kG(H^*) \otimes kG(H^*)$ such that  $|G(H^*_{\phi})| = 12$, $24$
or $48$. Then $H^*_{\phi}$ is not simple. This proves (i).

  (ii) We may assume that $|G(H)| = 16$ and $H$ is simple. As in the proof of part (i), it follows that $H^*$ contains a  Hopf subalgebra $B$ of
dimension $12$.  Consider the dual projection $q: H \to B^*$; we have $\dim H^{\co B^*} = 4$. Also, by a dimension argument using \cite{NZ}, the kernel of the restriction of $q$ to $G(H)$ must be of order $4$. Thus $H^{\co B^*} \subseteq kG(H)$ and $q$ is  normal. Hence $H$ is not simple in this case. \end{proof}

We summarize the results of this section in the following
corollary.

\begin{corollary}\label{red-48} If $H$ is simple, the possible coalgebra
types for $H$ and $H^*$ must be among the following:
\begin{align*} & (1, 4; 2, 2; 3, 4), \quad (1, 4; 2, 2; 6, 1), \quad (1, 4; 2, 11)
\\ & (1, 4; 2, 3; 4, 2), \quad (1, 8; 2, 2; 4, 2), \quad (1, 8; 2, 10).
\end{align*}
Moreover, either $H$ or $H^*$ must be of one of the types listed
in the first row, and if the type of $H$ is in the second row then
the type of $H^*$ is $(1, 4; 2, 11)$. \qed \end{corollary}

\section{Main result up to cocycle twists} We shall prove in this section that semisimple Hopf algebras of dimension 48 are not simple, \emph{up to a
cocycle twist}.

\begin{lemma}\label{qt12-cc} Suppose $H$ is simple.
Then $H$ admits no Hopf algebra quotient $H \to B$, where $B$ is a
cocommutative Hopf algebra of dimension 12.
\end{lemma}

\begin{proof} We shall consider separately the coalgebra types listed in Corollary \ref{red-48}.

\begin{claim} Suppose $H \to B$ is a Hopf algebra quotient where $B$ is
cocommutative and $\dim B = 12$. Then $H^*$ is of type $(1, 4; 2,
2; 3, 4)$ or $(1, 4; 2, 2; 6, 1)$ as a coalgebra. \end{claim}

\begin{proof} By assumption $B^* \subseteq H^*$ is a commutative
Hopf subalgebra of dimension 12. If $H^*$ is not of the prescribed
types, then there is a Hopf subalgebra $A\subseteq H^*$ of
dimension 8. Since $H$ is simple, we may assume $k[A, B^*] = H$.
On the other hand, there exists $1 \neq g \in Z(A)\cap G(A)$. If
$g \in B^*$, then $g$ is central in $k[A, B^*] = H$ and $H$ is not
simple. If $g \notin B^*$, then $|G(H^*)| = 8$ and $G(B^*)
\subseteq G(H^*)$ is a normal subgroup (of index 2). Then
$kG(B^*)$ is normal in $k[G(H^*), B^*] = H$. \end{proof}

Suppose first that $H$ is of type $(1, 4; 2, 2; 6, 1)$ as a
coalgebra. The lemma follows in this case from \ref{alg-struct},
since a cocommutative quotient Hopf algebra $B$ of dimension 12 of
$H$ is the same as a commutative Hopf subalgebra $B^*$ of index 4
of $H^*$, and since $H^*$ has an irreducible module of dimension
$6 > [H^*: B^*]$.

  Suppose that, as a coalgebra,  $H$ is of one of the
types $$(1, 4; 2, 3; 4, 2), \quad (1, 8; 2, 2; 4, 2), \quad (1, 8;
2, 10).$$ Then $H$ contains a Hopf subalgebra of dimension 16 and
$H^*$ is of type $(1, 4; 2, 11)$ as a coalgebra. Then the lemma
follows from Lemma \ref{tw-2}.

  Suppose next that  $H$ is of type $(1, 4; 2, 2; 3, 4)$.
Recall from Remark \ref{rmk-48} (v), that $G(H)$ and the
irreducible characters $\chi_1$ and $\chi_2$ of degree 2 give rise
to a Hopf subalgebra $A$ of dimension $12$, and for all
irreducible characters $\psi_1, \dots, \psi_4$ of degree 3 we have
$\psi_i^* \psi_i =  1 + \chi_1 + \psi + \psi'$, where $\psi$ and
$\psi'$ are fixed elements in $X_3$. We also have  $\chi_1^2 = 1 +
a + \chi_1$, where $\{ 1, a \} = G[\chi_1] \subseteq G(H)$.

Suppose that $H$ has a Hopf algebra quotient $q: H \to B = kG(B)$,
where $\dim B = 12$. We may assume that $m(1, q(\psi_i)) = 0$, for
all $i = 1, \dots, 4$. Otherwise, we would have a decomposition
$H^{\co B} = k1 \oplus W$, where $W$ is an irreducible left
coideal of dimension 3, implying that the restriction $q: A \to B$
is an isomorphism, and thus that $A$ is cocommutative, which is
absurd.

We claim that also $m(1, q(\chi_1)) = 0$. If not, since $q(\chi_1)
= \chi_1\vert_{B^*}$, then $m(1, q(\chi_1)) = 1$ by Frobenius
reciprocity. Then $q(\chi_1) = 1 + g$, where $1 \neq g \in G(B)$.
The relation $\chi_1^2 = 1 + a + \chi_1$, implies that $q(a) = g$.

Let $1 \neq t \in G(H) \cap H^{\co B}$. Then $q(t) = 1$, and
$(\Ind_{B^*}^{H^*}1)t = \Ind_{B^*}^{H^*}1$; so that $t = a \in
G[\chi_1]$. This is a contradiction. Therefore $m(1, q(\chi_1)) =
0$, as claimed.

  The relation $\psi_i^* \psi_i =  1 + \chi_1 + \psi +
\psi'$ implies that $m(1, q(\psi_i^*\psi_i)) = 1$. Write
$q(\psi_i) = \sum_{s \in G(B)}n_s s$. Then $\sum_{s}n_s = \deg
\psi_i = 3$, and $n_s = m(s, q(\psi_i)) = m(\psi_i,
\Ind_{B^*}^{H^*}s) \leq 1$. Hence, the set of all $s \in G(B)$
such that $n_s \neq 0$ has at least 3 elements. Now we have
$$q(\psi_i^*\psi_i) = q(\psi_i)^*q(\psi_i) = \sum_{s, u \in
G(B)}n_s n_u s^{-1}u, $$ whence $m(1, q(\psi_i^*\psi_i)) \geq 3$.
This is again a contradiction. Therefore the lemma is established
in this case.

It remains to consider the case where $H$ is of type $(1, 4; 2, 11)$.
Suppose that there is Hopf algebra quotient $H \to B$, where $B = kG(B)$ is a
cocommutative Hopf algebra of dimension 12. By Corollary \ref{red-48}, $B^*$ is of type $(1, 4; 2, 2)$ as a coalgebra. Hence, $G(B)$ has a normal subgroup of order 3, and there is a normal Hopf subalgebra $B_0 \subseteq B$, with $\dim B_0 = 3$.

We may write $H^{\co B} = k1 \oplus kt \oplus V_1$, where $V_1$ is
an irreducible left coideal of  dimension 2, and $1 \neq t \in
G(H)$. By Corollary \ref{cor-mk}, $t$ is central in $k[C_1]$,
where $C_1$ is the simple subcoalgebra containing $V_1$.

Consider the sequence of surjective Hopf algebra maps
$$\begin{CD}H @>q>> B @>q'>> B', \end{CD}$$
where $B' = B / B B_0^+$.  We have $\dim H^{\co B'} = 12$, and there exists an irreducible left coideal $U$ of $H$, not isomorphic to $V_1$ such that $m(U, H^{\co B'}) > 0$.

Let $\chi \in H$ and $C \subseteq H$ be, respectively, the
irreducible character and the simple subcoalgebra corresponding to
$U$. Decompose $q(\chi) = \chi\vert_{B^*}$ in the form $q(\chi) =
g + h$, where $g, h \in G(B)$. Then the restriction of $q$ induces
an epimorphism $q: k[C] \to k\langle g, h \rangle \subseteq
kG(B)$.

Similarly, the restriction of $q'q$ induces an epimorphism $q: k[C] \to k\langle g', h' \rangle \subseteq kG(B')$, where $q'q(\chi) = g' + h'$. In particular, $g'$ and $h'$ are the natural projections of $g$ and $h$ in $G(B')$.

By Frobenius reciprocity, we may assume that $g' = 1$. Hence
$$q'q(k[C]) = k\langle h' \rangle \neq B',$$ because, by Claim \ref{notcyclic}, $G(B')$ is not cyclic.
In particular $k[C] \neq H$.

On the other hand, since $g' = 1$, then $g \in B_0$ is of order 3
(because $m(1, q(\chi)) = 0$). Hence 3 divides the dimension of
$k[C]$. Since $H$ is simple, then $\dim k[C] = 12$ or $6$. Also,
$\dim k[C]^{\co B} = 2$, and it turns out that $1 \neq t$ is a
central group-like in $k[C]$.

Since this happens for every irreducible constituent $U$ of $H^{\co B'}$, not isomorphic to $V$, it follows that $t$ is central in the Hopf subalgebra generated by all simple subcoalgebras intersecting $H^{\co B'}$. This implies that $t$ is central in $H$. This finishes the proof of the lemma.  \end{proof}

\begin{lemma}\label{red-tw} Suppose that $H$ is simple. Then there exists an
invertible normalized 2-cocycle $\phi \in kG(\widetilde H) \otimes
kG(\widetilde H)$ such that $G(\widetilde H_{\phi}) \simeq G$,
where $G$ is the group defined in \ref{caso1}. Here, $\widetilde
H$ is one of the Hopf algebras $H$ or $H^*$.
\end{lemma}

Recall from Chapter \ref{twist} that $G$ is the semidirect product
$G = F \rtimes \Gamma$, where $F = \langle a: \ a^3 = 1 \rangle$,
and $\Gamma = \langle s, t: \ s^2 = t^2 = sts^{-1}t^{-1} = 1
\rangle$, corresponding to the action by group automorphisms of
$\Gamma$ on $F$ defined on generators by $s.a = a^{-1}$ and $t . a
= a^{-1}$.

Observe that, by Lemma \ref{conteo48}, $\widetilde H_{\phi}$ is of
one of the types $$(1, 12; 2, 9), \quad (1, 12; 3, 4), \quad (1,
12; 6, 1).$$ In particular, $\widetilde H_{\phi}$ is not simple.

\begin{proof} We shall show that there exist a Hopf subalgebra $A \subseteq \widetilde H$ such that $A \simeq \mathcal A_0$. Thus there is
an invertible normalized 2-cocycle $\phi \in kG(\widetilde H)
\otimes kG(\widetilde H)$ such that $A_{\phi} \simeq kG$, in view
of Proposition \ref{twist-pqq}. In particular, $|G(\widetilde
H_{\phi})|$ is divisible by 12, where $\widetilde H$ is either $H$
or $H^*$. Thus $|G(\widetilde H_{\phi})| = 12, 24$ or $48$. The
following claim implies that indeed $|G(\widetilde H_{\phi})| =
12$. Hence, $G(\widetilde H_{\phi}) \simeq G$.

\begin{claim} Suppose $|G(\widetilde H_{\phi})| = 24$ or $48$. Then $H$ is
not simple. \end{claim}

\begin{proof} If $|G(\widetilde H_{\phi})| = 24$, then $k G(\widetilde H_{\phi})$ is a normal Hopf subalgebra of
$\widetilde H_{\phi}$. On the other hand we must have
$G(\widetilde H) \subseteq G(\widetilde H_{\phi})$. Then
$\phi^{-1} \in k G(\widetilde H_{\phi}) \otimes k G(\widetilde
H_{\phi})$. Therefore, in view of Lemma \ref{hfsb-tw}, $B =
(kG(\widetilde H_{\phi}))_{\phi^{-1}}$ is a normal Hopf subalgebra
of $\widetilde H$ (because it has index 2).

Assume now that $|G(\widetilde H_{\phi})| = 48$; that is,
$\widetilde H_{\phi} = k \Gamma$, where $\Gamma = G(\widetilde
H_{\phi})$ is a group of order 48. For such a group $\Gamma$,
either $Z(\Gamma) \neq 1$ or $\Gamma$ contains a normal subgroup
of order 16, that necessarily contains $G(\widetilde H)$. In any
case, $\widetilde H$ (and thus $H$) is not simple, as claimed.
\end{proof}

Since $H$ is simple by assumption, we know that $H$ is of one of
the coalgebra types listed in Corollary \ref{red-48}. Eventually
taking $\widetilde H = H^*$, we may further assume that $H$ has
one of the following coalgebra types: $$(1, 4; 2, 2; 3, 4), \quad
(1, 4; 2, 2; 6, 1), \quad (1, 4; 2, 11).$$

\begin{claim}\label{notcyclic} The group $G(H)$ is not cyclic. \end{claim}

\begin{proof} If $H$ is of type $(1, 4; 2, 11)$,  the claim follows from Lemma \ref{tw-2} (i). If $H$ is of type $(1, 4; 2, 2; 6, 1)$, then $H$ contains a simple subcoalgebra $C$ of dimension 36 such that $gC = C = Cg$, for all $g \in G(H)$, and the quotient coalgebra $C/C(kG(H))^+$ is of dimension 9. If $G(H)$ is cyclic, then any twisted group algebra $k_{\alpha}G(H)$ is cocommutative, but this contradicts Corollary \ref{comm-pair-cor}, because  $C/C(kG(H))^+$ cannot have 4 isomorphism classes of simple comodules. Hence the claim follows also in this case.

It remains to consider the type $(1, 4; 2, 2; 3, 4)$. Suppose on the contrary that $G(H)$ is cyclic.  By Corollary \ref{red-48}, there is a quotient Hopf algebra $q: H \to Q$ of dimension 4, so that $\dim H^{\co q} = 12$. Keep the notation in Remark \ref{rmk-48} (v). Considering the eventual decomposition of $A^{\co q}$, $A_1^{\co q}$, we distinguish the following two possibilities for $H^{\co q}$:

(a) $H^{\co q} = k1 \oplus V \oplus W_1  \oplus W_2 \oplus W_3$,  where $V$ is an irreducible left coideal of dimension 2, and $W_i$'s are irreducible left coideals of dimension 3, $i = 1, 2, 3$.

(b) $H^{\co q} = k1 \oplus kt \oplus V  \oplus V' \oplus W_1 \oplus W_2$,  where $t \in G(H)$, $V, V'$ are  irreducible coideals of dimension 2, and $W_i$'s are irreducible left coideals of dimension 3, $i = 1, 2$.

Note that in this case $t$ is an element of order 2 in $G(H)$, and since this group is cyclic by assumption, necessarily $t = a \in G[\chi_1]$.

{\bf Case (a).} In this case the restriction of $q$ to $kG(H)$ induces an isomorphism $kG(H) \simeq Q$. Hence $G(H^*)$ contains a cyclic subgroup of order 4, and by Corollary \ref{red-48}, we may assume that $H^*$ is also of type $(1, 4; 2, 2; 3, 4)$ as a coalgebra. Hence there is a quotient Hopf algebra $H \to   B$, where $\dim B = 12$ and $H^{\co B} \subseteq H^{\co q}$ is a left coideal of dimension 4. Thus we must have a decomposition  $H^{\co B} = k1 \oplus  W$, where $W$ is an irreducible left coideal of dimension 3. Let $\psi_W \in H$ be the character of $W$. Then $g\psi_Wg^{-1} = \psi_W$, for all $g \in G(H)$, and  $(\psi_W)^* = \psi_W$ by Lemma \ref{idf-coinv}. As in Remark \ref{rmk-48} (v) we have a decomposition
\begin{equation}\label{main-rel}(\psi_W)^*\psi_W = 1 + \chi_1 + \psi +\psi'.\end{equation}
By Frobenius reciprocity, $m(1, (\psi_W)\vert_B) = 1$, hence $(\psi_W)\vert_B = 1 + \lambda$, where $\lambda \in B$ is a (not necessarily irreducible) character such that $m(1, \lambda) = 0$; then $m(1, (\psi_W)^*(\psi_W)\vert_B) \geq  2$.
On the other hand,
\begin{align*}m(1, (\psi_W^*\psi_W)\vert_B) & = 1 + m(1, (\chi_1)\vert_B) + m(1, \psi\vert_B) +  m(1, \psi'\vert_B) \\ & = 1 +  m(1, \psi\vert_B) +  m(1, \psi'\vert_B),\end{align*} and we may assume that $m(1, \psi) >   0$, whence $\psi_W = \psi$.

This implies that $\psi^* = \psi$. Then also $(\psi')^* = \psi'$, in view of the relation \eqref{main-rel}. Write $\psi' = g \psi$, $1 \neq g \in G(H)$.  Therefore $g^2 = 1$, and since $G(H)$ is cyclic $g = a \in G[\chi_1]$. But this implies that the irreducible characters $G[\chi_1], \chi_1, \psi, \psi'$ span a standard subalgebra which corresponds to a Hopf subalgebra of $H$  of dimension 24; hence $H$ is not simple in this case.

{\bf Case (b).}  Let $\psi_i$ be the character of $W_i$, $i  = 1, 2$. Then  $a\psi_1 = \psi_2 \neq \psi_1$, because $a W_1$ is not isomorphic to $W_1$ and it is contained in $H^{\co q}$. Moreover, $\psi_1 \in \{ \psi, \psi' \}$, since otherwise, we would have $\psi' = a\psi$, implying as before that $H$ contains a Hopf subalgebra of dimension 24.

Then we may assume that $\psi_1 = \psi$. Write $\psi' = g \psi$, $1 \neq g \in G(H)$. If $\psi^* = \psi$, then also $(\psi')^* = \psi'$ and therefore $g^2 = 1$, whence $g = a$ since $G(H)$ is cyclic.
If  $\psi^* = \psi'$, then $\psi' = \psi_2$, by Lemma \ref{idf-coinv}. Thus, in any case, $\psi' = a \psi$, implying as before, that $H$ contains a Hopf subalgebra of dimension 24. This proves that if $H$ is simple, then the group $G(H)$ is not cyclic, as claimed. \end{proof}

Suppose that $H$ has a Hopf subalgebra $A$ of dimension 12. By
Lemma \ref{qt12-cc}, $H^*$ has no cocommutative quotient of
dimension 12. Therefore, $A$ is not commutative. Since $G(H)$ is
not cyclic, then $A \simeq \mathcal A_0$ by \cite{fukuda}. By
Proposition \ref{twist-pqq} (i), there exists a 2-cocycle $\phi
\in kG(A)^{\otimes 2} = kG(H)^{\otimes 2}$ such that $\mathcal A_0
\simeq (kG)_{\phi}$. This establishes the lemma in this case.

  We may therefore assume that $H$ is of type $(1, 4; 2,
11)$ and contains no Hopf subalgebra of dimension 12, by Lemma
\ref{tw-2}. In particular, $H^*$ contains no Hopf subalgebra of
dimension 16,  by Theorem \ref{coinvariantes}. We may also assume that $H^*$ is not of type $(1, 4;
2, 2; 3, 4)$ nor $(1, 4; 2, 2; 6, 1)$, since otherwise we are done
by letting $\widetilde H = H^*$.

By Lemma \ref{qt16} (i), we may suppose that $H^*$ is also of type
$(1, 4; 2, 11)$. Then there is a quotient $q: H \to B$, where
$\dim B = 8$ and $B$ is not commutative. We may assume that
$H^{\co B} = k1 \oplus kt \oplus V_1 \oplus V_2$, where $t \in
G(H)$ and $V_i$ is an irreducible coideal of dimension 2, such
that $V_i$ is not isomorphic to $V_2$. Let $\tau_i \in X_2$ and
$C_i \subseteq H$ be the irreducible character and simple
subcoalgebra corresponding, respectively, to $V_i$, $i = 1, 2$.

Then $m(\tau_i, \Ind_{B^*}^{H^*}1) = 1$, and therefore $q(\tau_i)
= \tau_i\vert_{B^*} = 1 + g_i$, where $1 \neq g_i \in G(B)$.
Moreover, the subgroups $\langle g_1\rangle$, $\langle g_2
\rangle$ are of order at most 4. By construction, $q$ induces by
restriction a surjective Hopf algebra map $q: k[C_i] \to k\langle
g_i\rangle$; in particular, $k[C_i] \neq H$.

  If $k[C_i]^{\co B} = k1 \oplus V_i$, respectively $k1
\oplus kt \oplus  V_1 \oplus V_2$, then we have $\dim k[C_i, G(H)]
= 12$, respectively $\dim k[C_i] = 12$. Hence we are done in these
cases. If otherwise, $k[C_i]^{\co B} = k1 \oplus kt \oplus V_i$,
for all $i = 1, 2$, Then $tV_i = V_i = V_it$ and by Corollary
\ref{cor-mk}, $t$ is central in $k[C_i]$. Thus $t$ is central in
 $k[C_1, C_2]$;  and we may assume that $ k[C_1, C_2] \neq H$.
 Since $6$ divides $\dim k[C_1, C_2]$,  then we may even assume that $\dim k[C_1, C_2] = 12$
 and we are done in view of Lemma \ref{tw-2}. \end{proof}

\section{Main result} In this section we shall prove our main result in dimension 48.
For this we shall first study the normal Hopf subalgebras in a
semisimple Hopf algebra $K$ with $\vert G(K) \vert = 12$. In what
follows, we shall denote by $G$ the group considered in Lemma
\ref{red-tw}.

\begin{theorem} Let $H$ be a semisimple Hopf algebra of dimension 48. Then $H$ is not simple. \end{theorem}

\begin{proof} By Lemma \ref{red-tw}, we may assume that for $\widetilde
H = H$ or $H^*$, there exists an invertible normalized 2-cocycle
$\phi \in kG(\widetilde H) \otimes kG(\widetilde H)$ such that
$\widetilde H_{\phi}$ is of one of the types
\begin{equation*}(1, 12; 2, 9), \quad (1, 12; 3, 4),
\quad (1, 12; 6, 1),\end{equation*} and $G(\widetilde H_{\phi})
\simeq G$, where $G$ is as  in \ref{caso1}.

By Lemma \ref{12-no}, $\widetilde H_{\phi}$ is not simple. Note
that  $kG(\widetilde H) = (kG(\widetilde H))_{\phi}$ is a
cocommutative Hopf subalgebra of dimension 4 of $\widetilde
H_{\phi}$. Note also that, because $\widetilde H = \widetilde
H_{\phi}$ as algebras, we may assume that $|G(\widetilde
H_{\phi}^*)|$ is divisible by 4. On the other hand, $\phi^{-1} \in
kG(\widetilde H) \otimes kG(\widetilde H)$ is a normalized
2-cocycle for $\widetilde H_{\phi}$.

  Let $K$ be a semisimple Hopf algebra of dimension 48
such that $G(K) \simeq G$. In view of Lemmas \ref{prim}, \ref{seg}
and \ref{ter} below, at least one of the following conditions
hold:
\begin{flalign} \label{cond1} &K  \text{ has a nontrivial central group-like element; }& \\
\label{cond2} &G(K) \text{ is contained in a normal Hopf
subalgebra of } K.& \end{flalign}

Thus, if $J \in kG(K) \otimes kG(K)$ is an invertible 2-cocycle,
the twisted Hopf algebra $K_J$ is not simple, by Lemma
\ref{hfsb-tw} and Corollary \ref{andrus}. This implies the
theorem, applied to $K = \widetilde H_{\phi}$ and $J = \phi^{-1}$.
\end{proof}

\bigbreak In the rest of this section $K$ will be a semisimple
Hopf algebra of dimension 48 such that $G(K) \simeq G$; that is,
$G(K) \simeq G = F \rtimes \Gamma$, and thus $G$ contains a unique
abelian subgroup $M$ of order 6, $M = F \times Z$, where $Z =
Z(G)$ is of order 2.

Also,  $A \subseteq H$ will be a proper normal Hopf subalgebra.
Our aim is to show that at least one of the conditions
\eqref{cond1} or \eqref{cond2} hold.

  Hence, we may assume in what follows that $\dim A \neq
2$. Indeed, by Corollary \ref{kob-mas}, if $\dim A = 2$, then $A
\subseteq kG(K)\cap Z(K)$. Thus, in this case, $K$ verifies
condition \eqref{cond1}.

  On the other hand, note that we cannot have $A \cap
kG(K)$ of dimension 4: this would imply that $A \cap kG(K) = kL$,
where $L$ is a \emph{normal} subgroup of order $4$ in $G$ (because
$A \cap kG(K)$ is invariant under the adjoint action of $G$),
which contradicts the structure of $G$.

Hence $\dim A = 4, 8,  16$ are impossible.

  Therefore, the possibilities for $\dim A$ are 3, 6, 12
and 24. Moreover, if $\dim A = 12$ or 24, we may assume that $K$
is of type $(1, 12; 3, 4)$. Otherwise, (that is, if $K$ is of type
$(1, 12; 2, 9)$ or $(1, 12; 6, 1)$), since $\dim A \cap kG(K) \neq
4$, necessarily $G(K) \subseteq A$. See Lemma \ref{conteo24}.

\begin{lemma}\label{prim}
Assume $K$ is of type $(1, 12; 2, 9)$ as a coalgebra. Then at
least one of the conditions \eqref{cond1} or \eqref{cond2} hold.
\end{lemma}

\begin{proof} By Proposition \ref{cociente8}, $K$ contains a Hopf subalgebra
$K_0$ of dimension 8. We may assume that $k[K_0, F] = K$; if not,
$\dim k[K_0, F] = 24$ and $k[K_0, F]$ is normal in $K$, whence we
are done because $G(K) = G \subseteq k[K_0, F]$.

For all $g \in F$, $gK_0g^{-1}$ is an 8-dimensional Hopf
subalgebra of $K$. If $gK_0g^{-1} = K_0$ for some $1 \neq g \in
F$, then $G(K_0)$ would be a normal subgroup of order $4$ of $G$,
which is a contradiction. Thus the conjugation action of $F$ gives
rise to 8-dimensional Hopf subalgebras $K_0, K_1, K_2$ and
$G(K_i)$, $0\leq i\leq 2$, are distinct subgroups of order 4 of
$G$.

We may further assume that $F \subseteq k[K_0, K_1, K_2]$ and
therefore that $k[K_0, K_1, K_2] = K$. Otherwise $G(k[K_0, K_1,
K_2])$ would be of order $4$ and since $k[K_0, K_1, K_2]$ is
normal in $K = k[K_0, K_1, K_2, F]$, then this group would be
normal in $G(K)$, which is not possible. Note that for each $i =
0, 1, 2$, $G(K) = \langle G(K_i), F\rangle$.

We  may assume $K$ contains a normal Hopf subalgebra $A$ with
$\dim A = 3$ or $6$. If $\dim A = 3$, then $A = kF$, and in this
case, since $kG(K_i)$ is normal in $K_i$, then $\ad_{K_i}kG(K) =
kG(K)$. Thus $kG(K)$ is normal in $K = k[K_0, K_1, K_2]$. If on
the other hand, $\dim A = 6$, then because $\dim K/KA^+ = 8$, we
have $kF \subseteq A$. Hence $A$ is cocommutative and arguing as
before we get that $kG(K)$ is normal in $K$. This finishes the
proof of the lemma.
\end{proof}

\begin{lemma}\label{seg}
Assume $K$ is of type $(1, 12; 3, 4)$ as a coalgebra. Then at
least one of the conditions \eqref{cond1} or \eqref{cond2} hold.
\end{lemma}

\begin{proof} In this case may assume that $\dim A = 3, 6, 12$ or 24.
We may further assume that $K$ contains a normal Hopf subalgebra
of dimension 3; that is, $\dim A = 3$. To see this, we argue as
follows. Suppose that $\dim A = 12$. Then $A$ is of type $(1, 3;
3, 1)$ as a coalgebra, and $A$ is commutative. Also, $G(A)$ is the
unique (normal) subgroup of order 3 of $G(K) = G$. Thus $kG(A)$ is
a 3-dimensional normal Hopf subalgebra in $k[A, G] = K$.

If $\dim A = 6$, then  $A = kN$, where $N$ is a normal subgroup of
order 6 of $G(K)$. Consider the Hopf algebra quotient $K \to
K/KA^+$, since $\dim K/KA^+ = 8$, there is a Hopf algebra quotient
$K \to K/KA^+ \to B$, where $\dim B = 4$. We have $A \subseteq
K^{\co B}$, and unless $kG(K)$ is normal in $K$, $K^{\co B} = A
\oplus U \oplus V$, where $U$ and $V$ are irreducible left
coideals of $K$ of dimension 3; since $N K^{\co B} = K^{\co B}$,
we get that $U$ and $V$ are not isomorphic.

In addition, $xV \simeq V \simeq Vx$ and $xU \simeq U \simeq Ux$,
for all $x \in F$, because $G[\chi_U]$ and $G[\chi_V]$ are
necessarily of order 3.  Let $C_U$ and $C_V$ be the simple
subcoalgebras of $K$ containing $U$ and $V$, respectively. In view
of Corollary \ref{cor-mk}, $kF$ is normal in $k[C_V, C_U]$.
Therefore, since it is also normal in $kG(K)$, then $kF$ is normal
in $k[C_V, C_U, G(K)] = K$.

Finally, if $\dim A = 24$, we may assume that $A$ is of type $(1,
6; 3, 2)$ as a coalgebra; see Lemma \ref{conteo24}. In particular,
$G(A)$ is a normal subgroup of $G$ of order 6.

Let $A_0 \subseteq A$ be a normal Hopf subalgebra. If $\dim A_0 =
12$, then $kF \subseteq Z(A_0)$, and it follows that $kF$ is
normal in $K = k[A_0, G(K)]$. If, on the other hand, $A_0$ is
cocommutative, then $\dim A_0 = 2, 3$ or 6.

In the first case, $A_0 = kZ$ is the only subgroup of order 2
contained in $G(A)$. Hence $kZ$ is central in $A$ and thus also in
$k[A, G(K)] = K$. If $\dim A_0 = 3$, then $A_0 = kF$ and $kF$ is
normal in $k[A, G(K)] = K$. If $\dim A_0 = 6$, then $A_0 = kN$ is
normal in $k[A, G(K)] = K$ and $K$ has a normal Hopf subalgebra of
dimension 3, be the above.

\smallbreak We may therefore assume that $\dim A = 3$, as claimed.
Then $A = kF$, since $F$ is the only subgroup of order 3 in
$G(K)$. In particular, $F = G[\chi]$ for all irreducible character
$\chi$ of degree 3. Hence, the quotient Hopf algebra $K/KA^+$ is
cocommutative; see Remark \ref{particular}. Say $K/KA^+ \simeq
kN$, where $N$ is a group of order 16. Since $kF \simeq k^F$, then
$K$ fits into the abelian extension
\begin{equation}\label{3-16}1 \to k^F \to K \to kN \to 1.\end{equation}
We shall prove in what follows that in this case, $G(K) \cap Z(K) \neq 1$. Hence $K$ satisfies condition \eqref{cond1}.

  Consider the matched pair  $(F, N)$ associated to the
extension, with the actions $\triangleleft : F \times N \to F$ and
$\triangleright : F \times N \to N$. By \cite[Lemma 1.1.7]{pqq},
the exact sequence dual to \eqref{3-16} gives rise to an exact
sequence of groups $$1 \to \widehat F \to G(K) \to N^F,$$ where
$N^F$ is the subgroup of invariants in $N$ under the action
$\triangleright$. Since $|G(K)| = 12$, $N^F$ contains a subgroup
$N_1$ of order 4 of $N$.

On the other hand, the action $\triangleleft$ permutes the set
$F\backslash \{ 1 \}$; hence, there is a subgroup $N_0 \subseteq
N$ such that $N_0$ acts trivially on $F$ and $|N_0| = 8$. We have
$N_0 \cap N_1 \neq 1$. Let $1 \neq x \in N_0 \cap N_1$.

As a Hopf algebra, $K \simeq k^F {}^{\tau}\#_{\sigma}kN$ is a
bicrossed product, corresponding to the actions $kN \otimes k^F
\to k^F$ obtained from $\triangleleft$, and $kF \otimes k^N \to
k^N$ obtained from $\triangleright$, and certain 2-cocycles
$\sigma: N \times N \to (k^F)^{\times}$ and $\tau: F \times F \to
(k^N)^{\times}$. Moreover, by \cite[Lemma 1.2.5]{pqq}, we may
assume that $\tau = 1$.

The compatibility conditions between $\triangleleft$ and
$\triangleright$ imply that $F$ acts on $N_0$ by group
automorphisms through $\triangleright$. Hence $A_0 = k^F \otimes
kN_0$ is a Hopf subalgebra of $K$ of dimension 24.

  Since $x \in N_1$, then $F$ acts trivially on $x$ and $x
\in G(K)$. Also, since $x \in N_0$, then $x$ acts trivially in
$F$, and thus $x$ commutes with $k^F$ \cite{Maext}.

Since $F$ acts on $N_0$ by group automorphisms, this action
preserves the center of $N_0$. If $N_0$ is not abelian, then
$|Z(N_0)|$ is of order 2, and $F$ acts trivially on $Z(N_0)$.
Hence we may assume that $x \in Z(N_0)$. Therefore, $1 \neq x \in
G(A_0) \cap Z(A_0)$ is an element of order 2. Clearly, the same
conclusion holds if $N_0$ is abelian, since in this case any
$x \in N_1 \cap N_0$ commutes with $N_0$.

On the other hand, we may assume that $A_0$ is of type $(1, 6; 3,
2)$ as a coalgebra, whence $G(A_0) \cap Z(A_0) \supseteq Z(G(K))$,
because $A_0$ is normal in $K$. Therefore $Z(G(K))$ is central in
$k[A_0, G(K)] = K$. This proves the claim. The proof of the lemma
is now complete. \end{proof}

\begin{lemma}\label{ter} Assume $K$ is of type $(1, 12; 6, 1)$ as a coalgebra.
Then at least one of the conditions \eqref{cond1} or \eqref{cond2}
hold. \end{lemma}

\begin{proof} In this case may assume that $\dim A = 3$ or 6.
We claim that if  $K$ has a Hopf algebra quotient $K \to B$ with
$\dim B = 16$, then $kG(K)$ is normal in $K$. Hence the lemma
follows when $\dim A = 3$.

Observe that the quotient $q: K \to B$ is necessarily normal,
since $kF = K^{\co B}$ by \cite{NZ}. In order to establish the
claim, we shall follow the lines of the proof of Lemma
\ref{nq-16}. We have $q(kG(K))$ is a four-dimensional Hopf
subalgebra of $B$.  Let $\psi$ be the unique irreducible character
of degree 6. Then we have $\psi^2 = \sum_{g \in G(K)}g + 4 \psi$.
This implies, as in the proof of Lemma \ref{nq-16}, that $q(\psi)
= \lambda_1 + \lambda_2 + \lambda_3$, where $\lambda_i$ are
pairwise distinct irreducible characters of degree 2 in $B$. Hence
$B$ must be of type $(1, 4; 2, 3)$ as a coalgebra. Hence,
$q(kG(K)) = kG(B)$.

We know that $kG(B)$ is normal in $B$. Consider the sequence of
surjective Hopf algebra maps $K \overset{q}\to B \overset{q'}\to
B'$, where $B' = B / B (kG(B))^+$.  Since $q(kG(K)) = kG(B)$, then
we have  $kG(K) \subseteq K^{\co q'q}$. Thus, $kG(K) = K^{\co
q'q}$ is a normal Hopf subalgebra, as claimed.

\bigbreak Suppose next that $\dim A = 6$, so that $A = kS$, where
$S$ is a subgroup of order 6 in $G$. The quotient Hopf algebra
$K/KA^+$ is cocommutative  by Remark \ref{particular}. Say $K/KA^+
= k\Gamma$, where $|\Gamma| = 8$. Then $K^* = k^{\Gamma}
\#_{\sigma}k^S$ is a crossed product. Since there exists an
irreducible $K^*$-module of dimension 6, then $S = M$ is abelian,
by \cite[Proof of Theorem 2.1]{MoW}. Thus $K$ fits into the
abelian extension $1 \to k^M \to K \to kN \to 1$, where $N$ is a
group of order 8 such that $ K/KA^+ \simeq kN$. Dualizing, we get
an abelian extension $1 \to k^N \to K^* \to kM \to 1$. Let
$\triangleleft: N\times M \to N$, $\triangleright: N\times M \to
M$, be the associated matched pair.

Since the action $\triangleright$ fixes $1 \in M$, and because $N$
is of order 8, the set of fixed points $M^N$ has at least 2
elements. Moreover, by the compatibility condition
\cite[(4.10)]{Maext}, $M^N$ is a subgroup of $M$. It follows also
from the formulas \cite[(4.2) and (4.5)]{Maext} for the
multiplication and comultiplication of $K$ that the subspace $B: =
k^N \otimes k(M^N)$ is a Hopf subalgebra of $K^*$. If $B = K^*$, then
then the action $\triangleright $ is trivial. This implies that
$kM$ is a central Hopf subalgebra of $K$, which contradicts the
assumption on the structure of $G(K)$.

Therefore $\dim B = 8 |M^N| = 16$ or 24. If $\dim B = 16$, the
lemma follows from the proof in the case $\dim A = 3$. If $\dim B
= 24$, then $G(K) \cap Z(K) \neq 1$ and the lemma follows as well.
\end{proof}

\chapter{Dimension $54$}\label{54}

\section{First reduction} Let $H$ be a nontrivial semisimple Hopf algebra of dimension
$54$.

\begin{lemma}\label{conteo54} The order of $G(H^*)$ is
either  $2$, $6$, $9$, $18$ or $27$ and as an algebra $H$ is of
one of the following types:
\begin{itemize} \item $(1, 2; 2, 4; 6, 1)$, \item $(1, 2; 2, 4; 3, 4)$, \item $(1, 2; 2, 13)$, \item $(1, 6; 2, 3; 6, 1)$, \item $(1, 6; 2, 12)$, \item $(1, 6; 2, 3; 3, 4)$, \item $(1, 9; 3, 1; 6, 1)$, \item $(1, 9; 3, 5)$, \item $(1, 18; 2, 9)$, \item $(1, 18; 3, 4)$, \item $(1, 18; 6, 1)$, \item $(1, 27; 3, 3)$. \end{itemize} \end{lemma}

\begin{proof} It follows from \ref{alg-struct} and Corollary \ref{cor-g-1} that $|G(H^*)| \neq 1$.
If $|G(H^*)| = 2$, then the possible algebra types  for $H$ other
than the prescribed ones are
\begin{align*}& (1, 2; 2, 1; 4, 3), (1, 2; 3, 4;
4, 1), (1, 2; 4, 1; 6, 1), \\ & (1, 2; 2, 9; 4, 1), (1, 2; 2, 5;
4, 2).\end{align*} In the last two cases, by Theorems
\ref{thm-nr} and \ref{corolario} $H$ has a quotient Hopf algebra
of algebra types $(1, 2; 2, 9)$ or $(1, 2; 2, 5)$, respectively,
which contradicts \cite{NZ}. In the third case, $H$ has one
irreducible character $\chi$ of degree $4$, one irreducible
character $\psi$ of degree $6$ and the remaining two are of degree
$1$. The character $\chi$ is necessarily stable under left
multiplication by $G(H^*)$ and therefore we have a decomposition
$\chi^2 = \chi \chi^* = \epsilon + g + n \chi + m \psi$, implying
that $n = 2$ and $m = 1$. Then we have $1 = m = m(\chi, \psi
\chi)$ and $\psi \chi = \chi + l \psi$, $l \in \mathbb Z$, which
is not possible. The type $(1, 2; 2, 1; 4, 3)$ is discarded by
Lemma \ref{alg-type}.

Suppose that $H$ is of type $(1, 2; 3, 4; 4, 1)$.  Write $G(H) =
\{ \epsilon, g \}$ and let $\tau$ be the unique irreducible
character of degree $4$. We must have $\tau \tau^* = 1 + g + 2
\tau + \psi_1 + \psi_2$, where $\psi_1 \neq \psi_2$ are
irreducible characters of degree $3$. Then $m(\tau, \psi_1 \tau) =
1$ and therefore, $\psi_1 \tau = \tau + \sum_j \psi_j$, where
$\psi_j$ are irreducible of degree $3$. Taking degrees we get a
contradiction. This discards this possibility.

  The possibility $|G(H^*)| = 3$ is discarded using Theorem
\ref{thm-nr} and \ref{alg-struct}.  As to $|G(H^*)| = 6$, we need to discard the type $(1, 6; 4, 3)$, which is done as follows: let $\chi$ be a irreducible character of degree $4$, and write $\chi \chi^* = \sum_{g \in G[\chi]} g + \sum_{i}m_i \chi_i$, where $\chi_i$ are irreducible of degree $4$; it turns out that $4$ divides the order of $G[\chi]$ which is impossible.

Suppose that $|G(H^*)| = 9$. Then $H$ has
no irreducible characters of degree $2$ by Theorem \ref{thm-nr} and the
only possibilities are the prescribed ones. The rest of the lemma
follows also from \ref{alg-struct}. \end{proof}

\begin{lemma}\label{g9-54} Suppose that $H$ is of type $(1, 9; 3, 5)$ as a coalgebra. Then $H$ is not simple. \end{lemma}

\begin{proof} We claim that there exist irreducible characters $\psi \neq \psi'$
of degree $3$, which commute with $G(H)$, and such that $$\psi \psi^* = \psi' (\psi')^* = \sum_{g \in G(H)} g.$$
Indeed, the group $G(H)$, being abelian, acts by left multiplication on the set
$X_3' : = \{ \psi \in X_3: g\psi = \psi g, \forall g \in G(H) \}$,
and we have $|X_3'| = 2$ or $5$. It follows that there are at
least $2$ stable elements in $X_3'$, which therefore satisfy the
claimed equations.

The claim implies that there are two subcoalgebras $C_1$ and $C_2$
of  dimension $9$, such that $gC_i = C_i = C_i g$, for all $g \in
G (H)$. Moreover, we may assume that $k[C_1, C_2] = H$ since it is a Hopf subalgebra of dimension bigger than $18$, which divides $\dim H$: if $\dim k[C_1, C_2] = 27$, then $H$ is not simple by Corollary \ref{kob-mas}.

Fix  $i = 1, 2$. By Proposition \ref{dima=dimc},
$kG(H)$ is normal in $H =  k[C_1, C_2]$. Therefore $H$
is not simple in this case. \end{proof}

\begin{lemma}\label{g9-54-2} Suppose that $H$ is of type $(1, 9; 3, 1; 6, 1)$ as a coalgebra. Then $H$ is not simple. \end{lemma}

\begin{proof}  Let $D(H)$  be the Drinfeld double of $H$. By
\cite[Remark 2.2.4]{pqq} $|GD(H)^*| \neq 1$. By \cite{R}, the
group-like elements in $D(H)^*$ are of the form $g \otimes \eta$,
where $g \in G(H)$, $\eta \in G(H^*)$, are such that $\eta \otimes
g$ is central in $D(H)$. Observe that if the orders of $g$ and
$\eta$ are different, say $n = \vert \eta \vert < \vert g \vert$,
then the element $1 \otimes \epsilon \neq g^n \otimes \epsilon  =
(g \otimes \eta )^n$ would be such that $\epsilon \otimes g^n$ is
central in $D(H)$, and {\it a fortiori}, $1 \neq g^n \in G(H) \cap
Z(H)$.

Then we may assume that there is an
element $g \otimes \eta \in GD(H)^*$, where $g \in G(H)$ and $\eta
\in G(H^*)$ are of the same order. In particular, $|G(H^*)|$ is
divisible by $3$.

On the other hand, the irreducible  characters of degree $1$ and
$3$ span a standard subalgebra, which corresponds to a Hopf
subalgebra $A$ of dimension $18$. Consider the projection $H^* \to
A^*$; we may assume that $(H^*)^{\co A^*} = k1 \oplus V$, where
$V$ is an irreducible left coideal of dimension $2$. Then, by Theorem \ref{coinvariantes}, $|G(H^*)|$ is even and there is a
quotient Hopf algebra $H \to B$, with $\dim B = 3 |G(H^*)|$. Hence either
$\dim B = 18$ and necessarily $H^{\co B} \subseteq kG(H)$,  or else $H^*$ is of type $(1, 18; 2, 9)$ as a coalgebra and $H^{\co k^{G(H^*)}} \subseteq kG(H)$. This implies that
$H^*$ contains  a normal Hopf subalgebra, and hence $H$ is not
simple. \end{proof}

\begin{remark}\label{rmk-54} (i) If $H$ is of type $(1, 27; 3, 3)$, then $H$ is not simple by Corollary \ref{kob-mas}.

(ii) Suppose that $H$ is simple. It follows from Lemmas \ref{conteo54}, \ref{g9-54} and \ref{g9-54-2} that there exist subgroups $F \subseteq G(H)$ and $F' \subseteq G(H^*)$ such that $|F| = |F'| = 2$. We have a projection $q: H \to k^{F'}$ such that $\dim H^{\co q} = 27$. Hence, $F \cap H^{\co q} = 1$ and $H = R \# kF$ is a biproduct. By Proposition \ref{som}, $R$ is not commutative and not cocommutative.

(iii) Suppose that $H$ is of type $(1, 6; 2, 12)$ or $(1, 18; 2, 9)$. Then $H$ is not simple.

\begin{proof} Keep the notation in part (ii).
We claim that $G(H)$ contains a central (hence unique) subgroup of order $2$, which must stabilize all simple subcoalgebras of dimension $4$. This implies, in view of Remark \ref{particular}, that the braided Hopf algebra $R$ is a cocommutative  coalgebra and we are done.

To prove the claim we distinguish both cases. Suppose first that $H$ is of type
$(1, 6; 2, 12)$. Then $G(H)$ is abelian by Proposition \ref{nab-pq}, and the claim follows in this case.

In the case of type $(1, 18; 2, 9)$ we argue as follows. First note that by Lemma \ref{index-3}, $D(H)$ has an irreducible module of dimension $2$. It follows from Theorem \ref{thm-nr} that $G(D(H)^*)$ has an element $g \otimes \eta$ of order 2 or 3. If $g \otimes \eta$ is of order 2, we may assume that the center of $G(H)$ has an element of order $2$, as claimed.

Finally suppose that $g \otimes \eta$ is of order 3. By \cite[Corollary 2.3.2]{pqq} there is an exact sequence of Hopf algebras
 $1 \to kG \to D(H) \to K \to 1$, and also $G(K^*) \neq 1$ by \cite[Remark 2.2.4]{pqq}; moreover, we may assume that $g \otimes \eta \in G(K^*) \subseteq G(D(H)^*)$. Therefore, by part (ii) in \cite[Corollary 2.3.2]{pqq}, we get $\langle \eta,  g \rangle = 1$.

Consider the natural projection $q: H \to k^{\langle \eta \rangle}$; we have shown that $\langle g \rangle \subseteq H^{\co q}$. Also, it is clear that $H^{\co q}$ contains every 2-subgoup of $G(H)$. Since $H^{\co q}$ is of dimension 18, and we may assume does not coincide with $kG(H)$, then $G(H) \cap H^{\co q}$ is a normal subgroup of order 2 or 6 in $G(H)$. In particular, since we may assume that $1 \neq g \in Z(G(H))$, then  the group $G(H) \cap H^{\co q}$ is abelian. But $G(H) \cap H^{\co q}$ contains a unique subgroup of order 2, then so does $G(H)$. This finishes the proof of the claim.  \end{proof}

(iv) By Theorem \ref{thm-nr}, since $\dim H$ is not divisible
by $4$,  for every irreducible character $\chi$ of degree $2$ we
have $G[\chi] \neq 1$.

(v)  If  $H$ is of type $(1, 2; 2, 4; 6, 1)$, $(1, 2; 2, 4; 3,
4)$, $(1, 6; 2, 3; 6, 1)$ or $(1, 6; 2, 3; 3, 4)$ as a coalgebra,
then $G(H)$ and $X_2$ span a standard subalgebra of $R(H^*)$,
which corresponds to a non-cocommutative Hopf subalgebra $A$ of
dimension $18$. We may thus conclude  that if $H$ is simple, then
$H$ contains a Hopf subalgebra of dimension $18$. \end{remark}

\begin{lemma}\label{des-54} Assume that $H$ is of type $(1, 2;
2, 13)$ as a coalgebra. Then $H$ is commutative. \end{lemma}

\begin{proof} By Theorem \ref{1-2-2-n}, we may assume that the order of $G(H^*)$ is odd; thus, by Lemma \ref{conteo54}, $|G(H^*)| = 9$ or $27$. Also, by Corollary \ref{sek-rf}, we may assume that  $G(H) \cap Z(H) = 1$ and, in particular, $|G(H^*)| \neq 27$.

By \cite[Remark 2.2.4]{pqq}, we have $G(D(H)^*) \neq 1$. Since
$|G(H^*)|$ and $|G(H)|$ are relatively prime, the description of
$G(D(H)^*)$ in \cite{R} implies that there is an isomorphism
$$G(D(H)^*) \simeq \left(Z(H) \cap G(H)\right) \times \left(Z(H^*)
\cap G(H^*)\right) = Z(H^*) \cap G(H^*).$$ Therefore, $n =
|G(D(H)^*)| = 3$ or $9$, and $H^*$ contains a central group-like
element of order 3. Hence $H$  fits into an abelian cocentral
extension  $$1 \to k^G \to  H \to k\mathbb Z_3 \to 0,$$ where $G$
is a group of order 18. Then the extension induces the trivial
action  $\triangleright : G \times \mathbb Z_3 \to \mathbb Z_3$,
and an action by group automorphisms $\triangleleft: G \times
\mathbb Z_3   \to G$. This implies that the transpose action
$\rightharpoonup : k\mathbb Z_3 \otimes k^G \to k^G$ is by Hopf
algebra automorphisms, and therefore that $kG(H) = k\widehat G$ is
invariant under this action. Hence this action is trivial on
$kG(H)$. By \cite[Proposition 1.2.6]{pqq} $H$ is isomorphic as an
algebra to the crossed product $H \simeq k^G
\#_{\rightharpoonup}k\mathbb Z_3$. So it turns out that $kG(H)$ is
central in $H$. This implies that $H$ is commutative.  \end{proof}

\begin{lemma}\label{cor-54} Assume that $|G(H)| = 18$. Then $H$ is not simple. \end{lemma}

This discards the following possibilities for the coalgebra type of $H$:
$$(1, 18; 2, 9), \quad (1, 18; 3, 4), \quad (1, 18; 6, 1)$$

\begin{proof} In view of Remark \ref{rmk-54}
(iii) and Lemma \ref{conteo54},  we may assume that $H$ contains
no irreducible left coideal of dimension $2$. On the other hand,
we know from Remark \ref{rmk-54} (v) that there is a quotient Hopf
algebra $H \to B$, where $\dim B = 18$. Then necessarily $H^{\co
B} \subseteq kG(H)$ is a Hopf subalgebra, and thus $H$ is not
simple. \end{proof}

\begin{lemma} Suppose that $H$ is simple. Then $H$ contains a Hopf subalgebra
$A \subseteq H$, such that $A \simeq k^{\Gamma}$,  where $\Gamma$
is a non-abelian group of order $18$. In particular, the dimension
of an irreducible left $H$-module is at most $3$.  \end{lemma}

\begin{proof} By Remark \ref{rmk-54} (v) and Lemma \ref{cor-54}  we may assume that $H$ contains a non-cocommutative Hopf subalgebra $A$ of dimension
$18$.
By \cite{masuoka-further} there are two isomorphism classes of nontrivial semisimple Hopf algebras of dimension $18$, which are dual to each another: $\mathcal B_0$ and $\mathcal B_1 = {\mathcal B_0}^*$. We have $|G(\mathcal B_0)| = 6$ and  $|G(\mathcal B_1)| = 9$.

If $A$ is not commutative, and since $G(A) \subseteq G(H)$, we may
assume that $|G(A)| = 6$. By Remark \ref{rmk-54} (ii), $H = R \#
k\mathbb Z_2$, and since $k\mathbb Z_2 \subseteq A$ and $A$ is not
contained in $R = H^{\co \mathbb Z_2}$, then also $A \simeq R' \#
k\mathbb Z_2$. In particular, the order of $G(A^*)$ is also
divisible by $2$, which implies that $A$ is commutative as
claimed. Finally, applying \cite[Corollary 3.9]{harmonic} to the
inclusion $k^{\Gamma} \subseteq H$, we find that $\dim V \leq 3$,
for all irreducible $H$-module $V$. \end{proof}

\section{Main result} After what we have already shown until now, we can
conclude that if  $H$ is  simple, then the possible (co)algebra
types for $H$ are $(1, 2; 2, 4; 3, 4)$ and $(1, 6; 2, 3; 3, 4)$.

\begin{theorem} Let $H$ be a semisimple Hopf algebra of dimension
$54$ over $k$. Then $H$ is not simple. \end{theorem}

\begin{proof} Let $A \subseteq H$ be the commutative Hopf subalgebra of
dimension $18$. Note that $G(H) \subseteq A$.

By Lemma \ref{index-3} and Proposition \ref{doble}, we have
$G(D(H)^*) \neq 1$. Hence, we may assume that there exists a
nontrivial one-dimensional Yetter-Drinfeld module for $H$. By
Lemma \ref{yd-1}, this has necessarily  the form $V_{g, \eta}$,
for some $1 \neq g \in G(H)$ and $1 \neq \eta \in G(H^*)$.

Consider the projection $q: H \to k^{\langle \eta \rangle}$,
obtained by transposing the inclusion $k{\langle \eta \rangle}
\subseteq H^*$.
 Since $A$ is
commutative, $g^{-1}ag = a$, for all $a \in A^{\co q}$. By
Theorem \ref{A-comm-g}, we must have that $A^{\co q}$ is a Hopf
subalgebra of $A$. But this is not possible since $\dim A^{\co q}
= 9$, implying that $A^{\co q}$ is a cocommutative Hopf
subalgebra, while $G(H)$ is of order $2$ or $6$.
 This contradiction finishes the proof of the theorem.
\end{proof}

\chapter{Dimension $56$}\label{56}

\section{First reduction} Let $H$ be a nontrivial semisimple Hopf algebra of
dimension $56$.

\begin{lemma}\label{conteo56} The order of $G(H^*)$ is either $4$, $7$, $8$  or $28$, and as an algebra $H$ is of one of the following types:
\begin{itemize} \item $(1, 4; 2, 13)$,  \item $(1, 4; 2, 9; 4, 1)$,  \item $(1, 4; 2, 5; 4, 2)$, \item $(1, 4; 2, 1; 4, 3)$, \item $(1, 7; 7, 1)$, \item $(1, 8; 4, 3)$, \item $(1, 8; 2, 4; 4, 2)$, \item $(1, 8; 2, 8; 4, 1)$, \item $(1, 8; 2, 12)$, \item $(1, 28; 2, 7)$. \end{itemize} \end{lemma}

\begin{proof} We have $G(H^*) \neq 1$  by Corollary \ref{cor-g-1}; indeed, counting arguments show that the assumption $G(H^*) = 1$ implies that $H$ must have an irreducible character of degree 2.
Suppose that  $|G(H^*)| = 2$, and $H$ is of type $(1, 2; 2, n_2;
3, n_3; \dots )$. If $n_2 \neq 0$, by Theorems \ref{thm-nr} and
\ref{corolario}, there is a quotient Hopf algebra of dimension $2
(2n_2+1)$; thus $n_2 = 3$ by \cite{NZ}, and using \ref{alg-struct} we find a
contradiction. The possibilities with $n_2 = 0$ are the types
$(1, 2; 3, 2; 6, 1)$ and $(1, 2; 3, 6)$: these cases are seen to
be impossible, after trying to decompose the product $\psi \psi^*$
into irreducibles, where $\psi$ is an irreducible character of
degree $3$.

In the case $|G(H^*)| = 4$,  apart from the types
listed in the lemma, we have the following possibilities: $$(1, 4; 2,
4; 6, 1), (1, 4; 2, 4; 3, 4), (1, 4; 3, 4; 4, 1), (1, 4; 4, 1; 6, 1).$$
In the first two cases, $H$ has a
quotient Hopf algebra of dimension $20$, which is impossible. In
the third case, $H$ has four degree $1$ representations, four
irreducible characters $\chi_1, \dots, \chi_4$ of degree $3$, and
one irreducible character $\psi$ of degree $4$. Then we must have
$m(\chi_i, \psi^2) = m(\chi_i, \psi \psi^*) > 0$ for some $i$;
also, since $G(H^*)$ must permute transitively the degree $3$
characters under left or right multiplication and since $g \psi =
\psi$, for all $g \in G(H^*)$, then $m(\chi_i, \psi^2) > 0$ for
all  $i = 1, \dots, 4$, implying that $m(\psi, \chi_i \psi) =
m(\chi_i, \psi^2) = 1$ for all  $i = 1, \dots, 4$. Thus $\chi_1
\psi = \psi + \sum_{i} m_i \chi_i$,  where $m_i$ are non-negative
integers not all them equal to zero. Again since $\psi g = \psi$
for all $g$, we find that $m_i \neq 0$, for all $i = 1, \dots, 4$.
Taking degrees we find a contradiction. This discards the type
$(1, 4; 3, 4; 4, 1)$.

The type $(1, 4; 4, 1; 6, 1)$ is discarded as follows: let $\chi$ and $\psi$ be the unique irreducible characters of degrees $4$ and $6$, respectively. Write $\chi^2 = \sum_{g \in G(H^*)}g + n \chi + m \psi$; thus $m \neq 0$, since this would give a quotient Hopf algebra of dimension $20$, which is not possible. Hence $m = 2$, and since $m = m (\chi, \psi \chi)$, we have $\psi \chi = 2 \chi + l \psi$. Taking degrees we get a contradiction.

It is easy to see that $|G(H^*)| = 14$ is not
possible. The rest of  the lemma follows from \ref{alg-struct}.
\end{proof}

\begin{remark}\label{rmk-56} (i) If $H$ is of type $(1, 28; 2, 7)$, then $H$ is not simple, by Corollary \ref{kob-mas}.

  (ii) It follows from Proposition  \ref{cociente8} that,
except for the coalgebra type $(1, 7; 7, 1)$, $H$ has a Hopf subalgebra of dimension $8$.

  (iii) By \cite{pqq2},  if $G(H) \cap Z(H) = 1 = G(H^*)
\cap Z(H^*)$, then $G(H)$ and $G(H^*)$ are not both of order $8$.
\end{remark}

\begin{lemma} Suppose that $H$ is of type $(1, 7; 7, 1)$ as a coalgebra. Then $H$ is not simple. \end{lemma}

\begin{proof} We may assume that $H^*$ is also of type $(1, 7; 7, 1)$ as a coalgebra;
otherwise,  there is a quotient Hopf algebra $H \to B$, where
$\dim B = 8$ and necessarily $H^{\co B} = kG(H)$ is a normal Hopf
subalgebra.  Hence $H = R \# kG(H)$ is a biproduct, where $\dim R =
8$ and moreover, by Remark \ref{particular}, $R$ is cocommutative. Then
the lemma follows from Proposition \ref{som}. \end{proof}

\begin{lemma}\label{red56}  Suppose that $H$ is simple.

Then $H = R \# A$, where $A$ is a  semisimple Hopf algebra of
dimension $8$ and $R$ is a Yetter-Drinfeld Hopf algebra over $A$
of dimension $7$.  \end{lemma}

\begin{proof} We may assume that $H$  has a  Hopf subalgebra $A$ and a Hopf algebra
quotient $q: H \to B$, such that $\dim A = \dim B = 8$.

Since $\dim H^{\co B} = 7$,  we may assume that $kG(H) \cap H^{\co
B} = k1$. In particular, the lemma follows in the case where $A$
is cocommutative.

\smallbreak If $A$ is not cocommutative,  then $A$ contains a
unique $4$-dimensional simple coalgebra $C$ of $H$. Note that $C$
is not contained in $H^{\co B}$ since $G(H) \subseteq C^2$.

If the restriction $q\vert_A: A \to B$ is  an isomorphism, then we
are done. Otherwise, $H^{\co B}$ contains a $2$-dimensional
irreducible left coideal $V$ of $A$. But since $H^{\co B}$ cannot
contain $C$, we must have $H^{\co B} \cap A = k1 \oplus V$. This
contradicts Lemma \ref{restriction}. The proof of the lemma is complete.
\end{proof}

\section{Main result} In this section, we apply  Lemma \ref{red56} to show that a
semisimple Hopf algebra of dimension $56$ is not simple as a Hopf algebra.

\begin{theorem} Let $H$ be a semisimple  Hopf algebra of dimension $56$.
Then $H$ is not simple. \end{theorem}

\begin{proof} Suppose on the contrary that $H$ is simple.
Keep the notation in Lemma \ref{red56}.  After dualizing if
necessary, we may assume that $|G(H)| = 4$ (see Remark
\ref{rmk-56} (ii)); so that $A$ is not cocommutative. As a left
coideal of $H$, $R$ must be of one of the following types: $$(1) \
k1 \oplus V_1 \oplus V_2 \oplus V_3, \quad (2) \ k1 \oplus V
\oplus W,$$ where $\dim V_j = 2 = \dim V$, for all $j = 1, 2,  3$, and $\dim W = 4$. In particular, the type $(1, 4; 2, 1; 4, 3)$
is impossible.

Consider first the case (1). By Lemma \ref{bip-2-dim}
$\rho(V_j) \subseteq kG(A) \otimes V_j$, for all $j$; thus $\rho
(R) \subseteq kG(A) \otimes R$ and therefore $R \# kG(A)$ is a
Hopf subalgebra of $H$. This implies that $H$ is not simple, since
$G(A)$ has index $2$ in $A$.

Suppose finally that $R$ is as in case (2). By
Lemma \ref{bip-2-dim} $\rho(V) \subseteq kG(A) \otimes V$.  Moreover,
$kG(A) . V \subseteq V$, since $g V g^{-1}$ is a $2$-dimensional
irreducible left coideal of $H$ contained in $R$, for all $g \in
G(H)$.

By Proposition \ref{R-coideal}, $V$ is an $A$-subcomodule
subcoalgebra of $R$. Then, by Lemma \ref{util}, $B = k[V] \#
kG(A)$ is a Hopf subalgebra of $H$. Since $\dim B$ is divisible by
$4$ and $\dim k[V] \geq 3$, we must have $\dim B = 28$ and $kV =
R$. As before, we get that $R \# kG(A)$ is a Hopf subalgebra of
index $2$ in $H$, and $H$ is not simple. \end{proof}

\appendix

\chapter{Drinfeld Double of $H_8$}

We shall denote by $H_8$ the unique nontrivial 8-dimensional semisimple
Hopf algebra over $k$ \cite{k-p, ma-6-8}. We shall use the notation $D_4$ and $Q$ to indicate, respectively, the dihedral and quaternionic groups of order $8$. For a Hopf algebra $A$, $D(A)$ denotes the Drinfeld double of $A$.

By \cite{ma-6-8} the only non-commutative semisimple Hopf algebras of dimension $8$ are $H_8$, $kD_4$ and $kQ$.

\section{Structure of $D(H_8)$}
Tambara and Yamagami show in \cite{ty} that the categories of representations of these three Hopf algebras are not equivalent as monoidal categories. The comparison of Schur indicators implies  that the representation theory of $H_8$ is in some sense closer to that of $kD$ than to that of $kQ$; see \cite{Mo-8}.
On the other hand, note that $H_8$ fits into an exact sequence $1 \to k^{\Gamma} \to H_8 \to kF \to 1$, where $\Gamma \simeq \mathbb Z_2 \times \mathbb Z_2$ and $F \simeq \mathbb Z_2$, such that the associated product group associated to the extension is $D_4$.

The main results in this appendix are the following:

\begin{theorem}\label{th1}  (i) $D(H_8)$ fits into an abelian central extension
$$0 \to k^G \to D(H_8) \to kG \to 1,$$
where $G = G(D(H_8)^*) \simeq \mathbb Z_2 \times \mathbb Z_2 \times \mathbb Z_2$.

(ii) $D(H_8)$ is of type $(1, 8; 2, 14)$ as an algebra, and as a coalgebra it is of type $(1, 16; 2, 8; 4, 1)$.
\end{theorem}

\begin{theorem}\label{th2} $D(H_8)$ has no quotient Hopf algebra isomorphic to $kQ$. \end{theorem}

Observe that since $H_8$ admits quasitriangular structures \cite{suzuki}, then $D(H_8)$ also has quotient Hopf algebras isomorphic to $H_8$.

\begin{remark} Let $G$ be a finite group. Then $D(G)$ is a semidirect product $D(G) = k^G \rtimes kG$, with respect to the action comming from the adjoint action of $G$ on itself.
Hence, the irreducible representations of $D(G)$, viewed as irreducible Yetter-Drinfeld modules over $kG$, are classified by the modules $V_{g, \rho} : = kG \otimes_{Z_G(g)} \rho$, where $g$ runs over a system of representatives of the conjugacy classes in $G$ and $\rho$ runs over the irreducible representations of the centralizer $Z_G(g)$.

This implies that if $|G| = p^3$, $p$ prime, then $D(G)$ has exactly $p^3$ one-dimensional representations and the remaining irreducible representations are of dimension $p$; that is, $D(G)$ is of type $(1, p^3; p, p(p^3-1))$ as an algebra.

In particular, the Drinfeld doubles of the three non-commutative $8$-dimensional semisimple Hopf algebras have the same algebra structure. \end{remark}

\section[Proof]{Proof of Theorem \ref{th1}}

Let $A$ be a finite-dimensional Hopf algebra and let $g \in G(A)$, $\eta \in G(A^*)$. Let $V_{g, \eta}$ denote the vector space $k1$ endowed with the action $h.1 = \eta (h)1$, $h \in H$, and the coaction $1 \mapsto g \otimes 1$.

By Lemma \ref{yd-1}, the one-dimensional Yetter-Drinfeld modules over $A$ are exactly of the form $V_{g, \eta}$, where $g \in G(A)$ and $\eta \in G(A^*)$ are such that $(\eta \rightharpoonup h) g = g (h \leftharpoonup \eta)$, for all $h \in A$.

Note that if the condition $(\eta \rightharpoonup h) g = g (h \leftharpoonup \eta)$ holds for all $h$ in a set of generators of $A$, then it holds for all $h \in A$.

It turns out that the elements  of the form $\eta \otimes g$, where $g$ and $\eta$ satisfy the condition in Lemma \ref{yd-1} are exactly the central group-like elements in $D(A)$ \cite{R}.
In particular, $V_{g, \epsilon}$ (respectively, $V_{1, \eta}$) is a Yetter-Drinfeld module if and only if $g \in Z(A)$ (respectively, $\eta \in Z(A^*)$).

As in \cite{ma-6-8}, $H_8$ can be presented by generators $x$, $y$, $z$ with relations
\begin{align*}& x^2 = y^2 = 1, \\ & xy = yx, \quad zx = yz, \quad zy = xz, \\ & z^2 = \dfrac{1}{2}(1 + x + y - xy). \end{align*}
The coalgebra structure is determined by
\begin{align*}\Delta (x) & = x \otimes x,  \qquad \Delta (y) = y \otimes y, \\ \Delta(z) & = \dfrac{1}{2} \left( (1+y) \otimes 1 + (1-y) \otimes x \right) (z \otimes z). \end{align*}

In particular, $z \in H_8^{\times}$, $\epsilon(z) = 1$ and $\mathcal S(z) = z^{-1}$. We have in addition $H_8 \simeq H_8^*$,  $G(H_8) = \{ 1, x, y, xy \} \simeq \mathbb  Z_2 \times \mathbb  Z_2$, and $G(H_8) \cap Z(H_8) = \{ 1, xy \}$.

\begin{lemma} Let $g \in G(H_8) \backslash Z(H_8)$ and $\eta \in G(H_8^*) \backslash Z(H_8^*)$. Then $V_{g, \eta}$ is a Yetter-Drinfeld module of $H_8$. \end{lemma}

This implies that $D(H_8)$ has exactly $8$ distinct one-dimensional representations. In other words, we have
\begin{equation}\label{explicita}G(D(H_8)^*) = \langle xy \otimes \epsilon \rangle \oplus \langle 1\otimes \alpha \beta \rangle \oplus \langle x \otimes \alpha \rangle,
\end{equation}
where $G(H_8^*) = \{ \epsilon,  \alpha,  \beta,  \alpha \beta \}$ and $G(H_8^*) \cap Z(H_8^*) = \{ \epsilon,  \alpha \beta \}$.

\begin{proof} We shall show that if $1, xy \neq g$, then every element $\eta \in G(H_8^*) \backslash Z(H_8^*)$ satisfies the condition in Lemma \ref{yd-1}.

Notice that for any one-dimensional representation $\eta: H_8 \to k$ we must have $\eta (x) = \eta (y)$, because of the relation $zx = yz$ and the fact that $z \in H_8^{\times}$. Moreover, if $\eta$ is not central in $H_8^*$, then $\eta (x) = \eta (y) = -1$: indeed, if $\eta\vert_{G(H)} = 1$, then $k\langle \eta \rangle = (H_8^*)^{\co G(H_8)}$ and thus $\eta \in Z(H_8^*)$.

Also, it is enough to see that the condition in Lemma \ref{yd-1} is satisfied for $h = z$, since it is always satisfied for $h = x, y$ and these generate $H_8$ as an algebra.

We compute
\begin{align} \eta \rightharpoonup z & = \eta(z_2) z_1 = \dfrac{1}{2} \left( 1 + y + \eta(x) (1-y) \right) \eta(z) z, \\
z \leftharpoonup \eta  & = \eta(z_1) z_2 = \dfrac{1}{2} \left( \eta(1 + y)1 + \eta(1-y) x \right) \eta(z) z. \end{align}

Replacing in this identities $g = y$, the condition $(\eta \rightharpoonup z) y = y (z \leftharpoonup \eta)$ is equivalent to the equation $xy + x + \eta (x) (x - xy) = \eta (1 + y)  y + \eta (1-y) xy$.
This is always satisfied for $\eta \in G(H_8^*) \backslash Z(H_8^*)$.

The argument for $g = x$ is similar.  \end{proof}

\begin{corollary}\label{gp-lk} $G(D(H_8)^*) \simeq \mathbb  Z_2 \times \mathbb  Z_2 \times \mathbb  Z_2$. \end{corollary}

\begin{proof} We have $|G(D(H_8)^*)| = 8$ and $G(D(H_8)^*)$ is isomorphic to a subgroup of $G(D(H_8)) = G(H_8^*) \times G(H_8)$.  \end{proof}

{\it Proof of Theorem \ref{th1}.}  (i) In virtue of \cite[Corollary 2.3.2]{pqq} and Corollary \ref{gp-lk},  there is a central extension
$$\begin{CD}0 @>>> kG @>{\iota}>> D(H_8) @>{\pi}>> K @>>> 1, \end{CD}$$
where $G \simeq G(D(H_8)^*) \simeq \mathbb Z_2 \times \mathbb  Z_2 \times \mathbb  Z_2$, and the map $\iota$ is determined by $\iota(g \otimes \eta) = \eta \otimes g$. In particular, $\dim K = 8$.

We have $kG \simeq k^G$. Identify $K$ with a Hopf subalgebra of $D(H_8)^*$. Part (i) will be established if we show that $G(K^*) = G$. To see this, we observe that a group-like element $g \otimes \eta$ belongs to $G(K^*)$  exactly when it is a one dimensional representation of $D(H_8)$ which factorizes through $K$; that is, $g \otimes \eta$ belongs to $G(K^*)$
 if and only if
\begin{equation}\langle \eta, g'\rangle \langle \eta', g\rangle = 1,\end{equation} for all group-like elements  $g' \otimes \eta' \in G(D(H_8)^*)$. See \cite[Corollary 2.3.2]{pqq}.
Finally, using the description of the elements in $G(D(H_8)^*)$ given in \eqref{explicita}, one sees that $G(K^*) = G$, thus proving part (i).

(ii) As a coalgebra, $D(H_8)$ is a tensor product: $D(H_8) = (H_8^*)^{\cop} \otimes H_8$. This proves the statement corresponding to the coalgebra structure.

Combining part (i) with the description in \cite[Theorem 3.3]{KMM} (the action being trivial in our situation), we find that the simple $D(H_8)$-modules are of the form $p_g \otimes V$, where $g \in G$ and $V$ is an irreducible $k_{\sigma_g}G$-irreducible module, for some $2$-cocycle $\sigma_g: G \times G \to k^{\times}$.
The dimensions of the irreducible $k_{\sigma_g}G$-modules are either $1$ or $2$ (since they divide the order of $G$ and their square is less than $8$). This finishes the proof of (ii). \qed

\section[Proof]{Proof of Theorem \ref{th2}}

 Suppose that $\pi: D(H_8) \to kM$ is a surjective Hopf algebra map, where $M$ is a group of order $8$. Then $kM = \pi(H_8^*) \pi(H_8)$ and $\dim \pi(H_8^*)$ and $\dim \pi(H_8)$ divide  $4$ (because $H_8 \simeq H_8^*$ is not isomorphic to $kM$). Then some of them, say $\pi (H_8)$ is of dimension $4$. Then $\pi(H_8) \simeq \mathbb Z_2 \times \mathbb Z_2$, since $H_8^*$ has no Hopf subalgebra isomorphic to $\mathbb Z_4$. Therefore $M \neq Q$.  \qed

\backmatter

\bibliographystyle{amsalpha}

\end{document}